\documentstyle{amsppt}
\voffset=-1.5cm
\parskip 6pt
\NoBlackBoxes
\magnification=1200
\def\exam#1#2{\noindent\hangindent=3em\hangafter1
                \hbox to3em{#1\hfil\quad}#2}

\def\exama#1#2#3{\dimen1=2.25em \dimen2=4.5em
                \noindent\hangindent=\dimen2\hangafter1
                \hbox to\dimen1{#1\hfil~~}\hbox to\dimen1{#2\hfil~~}#3}  
\def\card{\text{card}}
\def\var{\text{var}}
\def\Var{\text{Var}}

\def\rbx{\hfill{\vbox{\hrule\hbox{\vrule\kern6pt\vbox{\kern6pt}\vrule}
          \hrule}}}

\def\bx{\vbox{\hrule\hbox{\vrule\kern6pt\vbox{\kern6pt}\vrule}\hrule}}
\rightheadtext{Concentration of Measure and Isoperimetric 
Inequalities}
\leftheadtext{M. Talagrand}

\topmatter
\title
Concentration of Measure and Isoperimetric\\
Inequalities in Product Spaces
\endtitle
\author
Michel Talagrand$^{(\ast)}$
\endauthor
\address
C.N.R.S. and The Ohio State University
\endaddress
\date (February `94)
\enddate
\thanks
AMS Classification numbers:  Primary 60E15, 28A35, 60G99; Secondary 
60G15, 68C15.
\endthanks
\thanks
$(\ast )$ Work partially supported by an NSF grant.
\endthanks
\abstract
The concentration of measure prenomenon roughly states that, if a set 
$A$ in a product $\Omega^N$ of probability spaces has measure at least 
one half, ``most'' of the points of $\Omega^N$ are ``close'' to $A$.  
We proceed to a systematic exploration of this phenomenon.  The 
meaning of the word ``most'' is made rigorous by isoperimetric-type 
inequalities that bound the measure of the exceptional sets.  The 
meaning of the work ``close'' is defined in three main ways, each of 
them giving rise to related, but different inequalities.  The 
inequalities are all proved through a common scheme of proof.  
Remarkably, this simple approach not only yields qualitatively optimal 
results, but, in many cases, captures near optimal numerical 
constants.  A large number of applications are given, in particular 
in Percolation, Geometric Probability, Probability in Banach Spaces, 
to demonstrate in concrete situations the extremely wide range of 
application of the abstract tools.
\endabstract
\endtopmatter

\vfil\eject

\centerline{\bf Table of Contents}
\bigskip

\exam{I.}{Introduction}

\exam{1.}{Inequalities }

\exam{2.}{Control by one point}

\exama{~}{2.1}{Basic principle}
\smallskip
\exama{~}{2.2}{Sharpening}
\smallskip
\exama{~}{2.3}{Two point space}
\smallskip
\exama{~}{2.4}{Penalties, I}
\smallskip
\exama{~}{2.5}{Penalties, II}
\smallskip
\exama{~}{2.6}{Penalties, III}
\smallskip
\exama{~}{2.7}{Penalties, IV}

\exam{3.}{Control by $q$ points}

\exama{~}{3.1}{Basic result}
\smallskip
\exama{~}{3.2}{Sharpening}
\smallskip
\exama{~}{3.3}{Penalties}
\smallskip
\exama{~}{3.4}{Interpolation}

\exam{4.}{Convex hull}

\exama{~}{4.1}{Basic results}
\smallskip
\exama{~}{4.2}{Sharpening}
\smallskip
\exama{~}{4.3}{Two point space}
\smallskip
\exama{~}{4.4}{Penalties}
\smallskip
\exama{~}{4.5}{Interpolation}

\exam{5.}{The Symmetric group}

\exam{II.}{Applications}

\exam{6.}{Bin Packing}

\exam{7.}{Subsequences}

\exama{~}{7.1}{Longest increasing subsequence}
\smallskip
\exama{~}{7.2}{Longest common subsequence}

\exam{8.}{Percolation}

\exama{~}{8.1}{Basic results}
\smallskip
\exama{~}{8.2}{General moments}
\smallskip
\exama{~}{8.3}{First time passage in Percolation}

\exam{9.}{Chromatic number of random graphs}

\exam{10.}{The Assignment Problem}

\exam{11.}{Geometric Probability}

\exama{~}{11.1}{Irregularities of the Poisson point process}
\smallskip
\exama{~}{11.2}{The traveling salesman problem}
\smallskip
\exama{~}{11.3}{The minimum spanning tree}
\smallskip
\exama{~}{11.4}{The Gabriel graph}
\smallskip
\exama{~}{11.5}{Simple matching}

\exam{12.}{The free energy of spin glasses at high temperature}

\exam{13.}{Sums of vector valued independent random variables}
\vfil\eject

\noindent {\bf 1.~~Introduction}

Upon reading the words ``isoperimetric inequality'' the average reader 
is likely to think to the classical statement:

(1.1)  Upon the bodies of a given volume in $\Bbb R^N$, the ball is 
the one with the smallest surface area.

This formulation, that needs the notion of surface area, is not very 
appropriate for generalization in abstract setting.  A less known 
(equivalent) formulation is as follows:

(1.2)  Among the bodies $A$ of a given volume in $\Bbb R^N$, the one 
for which the set $A_t$ of points within Euclidean distance $t$ of 
$A$ has minimum volume is the Euclidean ball.

It should be intuitive, taking $t\to 0$, that (1.2) implies (1.1).  We 
will, however, rather be interested in large values of $t$.  A first 
sight, this is uninteresting; but this first impression is created 
only by our deficient intuition, that functions correctly only for 
$N\le 3$, and lamely fails for the large values of $N$ that are of 
interest here.

For our point of view, the main feature of (1.2) is that it gives a 
lower bound on the volume of $A_t$ that depends only on $t$ and the 
volume of $A$.

From now on, all the measures considered will be probabilities (i.e. 
of total mass one).  The concept of concentration of measure 
phenomenon largely arose through the work of V. Milman on Dvoretzky's 
theorem on almost Euclidean sections of convex bodies.  Following 
[G-M], [Mi-S], the basic ideas may be described in the following way.  
Consider a (Polish) metric space $(X,d)$.  For a subset $A$ of $X$, 
consider the $d$-ball $A_t$ centered on $A$, i.e.
$$\leqalignno{A_t &= \{ x\in X\colon d(x,A)\le t\}\,. &(1.3)\cr}$$
Consider now a Borel probability measure $P$ on $X$.  The 
concentration function $\alpha (P,t)$ is defined as
$$\alpha (P,t) = \sup\{ 1-P(A_t)\colon P(A)\ge {1 \over 2}\,,~A\subset
X\,,~A~\text{Borel}\}\,.$$
In other words
$$\leqalignno{P(A) &\ge {1 \over 2}\quad\Rightarrow\quad P(A_t)\ge 
1-\alpha (P,t)\,.&(1.4)\cr}$$
It turns out that in many situations the function $\alpha (P,t)$ 
becomes extremely small when $t$ grows.  In rough words, if one starts 
with any set $A$ of measure $\ge 1/2$, $A_t$ is almost the 
entire space.  This is the concentration of measure phenomenon, that 
was promoted most vigorously by V. Milman.  It plays an important 
role in local theory of Banach spaces, and has become the 
central concept of the area of probability known as Probability in 
Banach spaces.  (See the book [L-T2], and subsequent work such as 
[T6], [T7]).

A prime example of space where concentration of measure holds is the 
Euclidean sphere $S_N$ of $\Bbb R^{N+1}$ equipped with its geodesic 
distance $d$ and normalized Haar measure $P_N$, for which it can be 
shown that
$$\leqalignno{\alpha (P_N,t) &\le \left({\pi \over 8}\right)^{1/2}
\exp\left( -{(N-1) \over 2}t^2\right)\,. &(1.5)\cr}$$
Closely related, and more in line with the topic of the present paper 
is the case $X=\Bbb R^N$, equipped with the Euclidean distance and 
the canonical Gaussian measure $\gamma_N$ (whose covariance is the 
Euclidean dot product).  In that case
$$\leqalignno{\alpha (\gamma_N,t) &\le \int^\infty_t {1 \over 
\sqrt{2\pi}}e^{-u^2/2}du\le {1\over 2} e^{-t^2/2}\,. &(1.6)\cr}$$

It should be pointed out that more is known.  The Gaussian 
isoperimetric inequality states that
$$\leqalignno{\gamma_N(A) &= \gamma_1 ((-\infty ,a])\Rightarrow
\gamma_N (A_t)\ge\gamma_1((-\infty ,a+t]) &(1.7)\cr}$$
which implies (1.6) when $a=0$.  However, it is sufficient for many 
applications to know (1.6) or even the weaker inequality
$$\leqalignno{\alpha (\gamma_N,t) &\le Ke^{-t^2/K} &(1.8)\cr}$$
where $K$ is a universal constant.

In the present work we perform a systematic investigation of the 
concentration of measure phenomenon in product spaces.  Thus with the 
terminology above, $X$ will be a product of probability spaces, and 
$P$ a product measure.  The statements will have the form (1.4).  
However, the set $A_t$, which consists of points close in a certain 
sense to $A$, (and that, for convenience, we will call the 
$t$-fattening of $A$) will not always have the form (1.3).  Indeed, 
it turns out that it is extremely fruitful to consider various notions 
of fattening.  We will define three rather distinct notions
of fattening.  These notions are studied respectively in Chapters 2 
to 4.  Each of these notions can be studied with various level of 
sophistication, and they are at times closely connected.  Discussing 
the whole theory in this introduction would require too much 
repetition and is inappropriate for an article of the present length.  
Thereby, we have decided to mention here only the main new 
theme (that did not appear in this author's previous work) as well 
as a simple result that appears to have a considerable potential for 
applications.

Assume that $X=\Omega^N$ is a product of probability spaces, and that 
$P=\mu^N$ is a product probability.  We recall that the Hamming 
distance $d$ on $X$ is given by
$$\leqalignno{d(x,y) &= \text{card}\{ i\le N\colon x_i\not= y_i\}\,.
 &(1.9)\cr}$$

When $A_t$ is given by (1.3), where $d$ is the Hamming distance, an 
important result, proved in special cases in [Mi-S] (with a proof 
that extends verbatim to the general situation) is that the 
concentration function $\alpha (P,t)$ satisfies 
$$\leqalignno{\alpha (P,t) &\le K\exp \left(-{t^2 \over K}\right)\,, 
&(1.10)\cr}$$
where $K$ is a universal constant.

One could interpret (1.9) by saying that we put a penalty $1$ for 
each coordinate $i$ where $x_i\not= y_i$.  One recurring theme of the 
present paper is the investigation of what happens when, instead, we 
put a penalty $h(x_i,y_i)$, where $h(x,y)$ is a non-negative function 
on $\Omega^2$.  A striking and unexpected finding is that in several 
instances there is a high disymmetry between the roles of $x$ and $y$.  
For example, in one of the main results of the paper (Theorem 4.4.1) 
if one requires that $h(x,y)$ should depend on $x$ only, it has to be 
bounded; but, if it depends on $y$ only, weak integrability conditions 
suffice.

Suppose now that $(\alpha_i)_{i\le N}$ are positive numbers, and let 
us replace the distance (1.9) by
$$d_\alpha (x,y)=\sum_{i\le N}\alpha_i 1_{\{x_i\not= y_i\}}\,.$$
It is then shown in [Mi-S] that (1.10) can be extended into
$$\leqalignno{\alpha (P,t) &\le K\exp\left( -{t^2 \over K\sum_{i\le
N}\alpha^2_i}\right)\,. 
&(1.11)\cr}$$
One way to spell out this result is as follows:

Given $A\subset\Omega^N$, with $P(A)\ge{1 \over 2}$, then, for all 
numbers $(\alpha_i)_{i\le N}$, $\alpha_i\ge 0$, $\sum\limits_{i\le N}
\alpha^2_i=1$, we have
$$\leqalignno{P(A_{t,\alpha}) &\ge 1-K\exp\left( -{t^2 \over K}\right)
 &(1.12)\cr}$$
where
$$A_{t,\alpha} =\left\{ x\in\Omega^N\colon \exists y\in A\,,~\sum_{i
\le N}\alpha_i1_{\{ x_i\not= y_i\}}\le t\right\}\,.$$

The first result of Chapter 4 states that (1.12) can be improved into
$$\leqalignno{P\left(\bigcap\limits_\alpha A_{t,\alpha}\right) &\ge 
1-K\exp\left( -{t^2 \over K}\right) &(1.14)\cr}$$
where the intersection is over all families $\alpha =(\alpha_i)_{i\le 
N}$ as above.  The power of this principle (that will be considerably 
perfected in Chapter 4) is by no means obvious at first sight, but 
will be demonstrated repeatedly through Chapters 6 to 9 (the easiest
applications being in Chapter 6 and 7).

We have explained in terms of sets what is the concentration of 
measure phenomenon.  However, rather than sets, one is more often 
interested in {\it functions}.  In that case, the concentration 
of measure phenomenon takes the following form:  if a function $f$ 
on $X$ is sufficiently regular, it is very concentrated around its 
median (hence around its mean).  If $M_f$ is a median of $f$, this is 
expressed by a (fast decreasing) bound on $P(\vert f-M_f\vert >t)$.  
For a simple example, (1.4) implies that if $f$ has a Lipschitz 
constant $1$ with respect to the underlying distance
$$\leqalignno{P(\vert f-M_f\vert\ge t) &\le 2\alpha (P,t)\,. 
&(1.15)\cr}$$

Despite the fact that functions are potentially more important than 
sets, all our concentration of measure results are stated in terms of 
sets.  (This is done in Part I.)  The essential reason for this choice 
is that the power and the generality of these results largely arise 
from the fact that they require only minimal structure (a condition 
better achieved by considering sets only).  A secondary reason is 
that much of the progress reported on the present paper (including on 
some rather concrete questions presented in Part II) has been 
permitted, or at least helped by the abstract point of view; and 
thereby, it seems worthwhile to promote this approach.  Nevertheless, 
the natural domain of application of the tools of Part I is the 
obtention of bounds on $P(\vert f-M_f\vert\ge t)$ when $f$ is a 
function defined on a product of measure spaces.  We will, however, 
give no abstract statement of this type.  We prefer instead to 
analyze a number of specific situations, reducing each time to 
statements about sets (the great variety of situations encountered 
indicates that this is possibly a clever choice).  This is the 
purpose of Part II, where we will demonstrate the efficiency of the 
tools of Part I.  It must be said that these specific situations have 
been of considerable help in pointing out the directions in which the 
abstract theory should be developed.  Most of the abstract results 
are indeed directly motivated by applications.

Certainly there is a considerable number of situations where occur 
naturally functions that are defined on a product of many measure 
spaces, or equivalently that depend on many independent random 
variables.  The examples presented here are certainly influenced by 
the past interests of the author.  Their boundary, however, is likely 
to reflect the limited knowledge of this author rather than the limit 
of the power of abstract tools of Part I.  (Should a reader be aware 
of another potential domain of application, he is urged to introduce 
it to this author.)  Quite logically, several of the examples we 
present have an ``applied'' flavor.  This is simply because stochastic 
models occur in physics (such as Percolation and spin glasses) and 
Computer Science (bin packing, assignment problem, geometric 
probability).  The reason for the later is that these stochastic 
models do shed some light on the behavior of computationally 
intractable problems, and, for this reason, are widely studied today; 
see e.g., [C-L], [~~~].  No previous knowledge whatsoever of these 
problems is required for reading the material of Part II, that we 
briefly describe now.

Each of the examples of Part II studies the deviation of a specific 
function $f$ of many independent random variables from its mean.  In 
each example, the function $f$ is obtained as the solution of an 
optimization problem.  This is not a coincidence, but rather reflects 
the fact that such situations are well adapted to the use of our 
methods.  In Chapter 6, we apply (4.1.3) to stochastic bin packing.  
This simple application is presented first since it is while 
considering this problem (while proctoring a Calculus exam) that the 
power of (4.1.3) beyond Probabilities in Banach spaces was first 
realized.  The application is not really typical.  More typical is 
the application of Chapter 7, to the length of the longest increasing 
subsequence of a random permutation.  This application puts forward 
the fact that when one studies the size of substructures whose 
existence is determined by a comparatively small number of random 
variables, rather than by the whole collection of random variables, 
inequality (4.1.3) fully takes advantage of that feature.  This 
characteristic occurs again in Chapter 8, where it is presented as a 
general result, that allows, as a rather weak and special corollary, 
to improve upon H. Kesten's recent results on first time passage in 
Percolation [K2].  In Chapter 9, we show how (4.1.3) again provides a 
natural approach to questions on random graphs.  The challenge of the 
Assignment problem considered in Section 10 is that the objective 
function $f$ considered there is very small; it is of order one, while 
depending on $N^2$ independent variables of order one, each of them 
with a potentially disastrous influence on the objective function.  
In Chapter 11, we consider situations where the objective function $f$ is 
defined in a geometrical manner from a random set of $N$ points in the 
unit square.   The common objective is to prove that $f$ has 
Gaussian-like tails. However, the richness of the situation is 
unsuspected beforehand; apparently similar definitions require rather 
different levels of sophistication.  In Chapter 12, we provide a 
simple derivation of the free energy of spin glasses at high 
temperature.  Finally, in Chapter 13, we discuss how the study of 
sums of vector-valued independent random variables motivated the 
approach of this paper, and we discuss a few new specific results.

We now comment on the methods of Part I, their history, and compare 
them with competing methods.

There is a general method, that is becoming increasingly popular, to 
prove deviation inequalities for $\vert f-Ef\vert$.  (That the mean 
rather than the median is involved is very much irrelevant).  It is 
to decompose $f$ as the sum of a martingale difference sequence $f=
\sum d_i$, and to use martingale inequalities.  The generality of the 
method stems from the fact that such a decomposition is easy, simply 
writing $d_i=E(f\mid {\Cal F}_i)-E(F\mid {\Cal F}_{i-1})$ for any 
increasing filtration $({\Cal F}_i)$.  This method was used in 
Probability in Banach Spaces (under the name of ``Yurinski's method'') 
for the study of $f=\left\|\sum\limits_{i\le N} X_i\right\|$, where 
$X_i$ are independent Banach space random variables (r.v.).  The 
generality of the method was discovered by B. Maurey [M1], and it was 
further developed in [Mi-S].  It soon became apparent, however, that 
this method would not always yield optimal results; this is what 
prompted the invention of the isoperimetric inequality of [T2] (more 
details on history are given in Chapter 12).  An inequality  very
similar to the inequality of [T2], but with a much simpler proof, 
appears in the present paper as Theorem 3.1.1.  The phenomenon 
described by this inequality was completely new at that time, and had 
a major impact in Probability in Banach spaces (prompting, in 
particular, the writing of the book [L-T2]).  One could reasonably 
hope that this inequality would find applications to other domains; 
but as of today, this has not been the case.  Another inequality that 
was discovered in relation with Probability in a Banach space is a 
predecessor of (4.1.3) [T1].  The inequality of [T1] did not, however, 
play a crucial role in that theory, because, for most applications, 
it could be replaced by the Gaussian isoperimetric inequality (1.6) 
to which it is related.  For this reason, the discovery that (4.1.3) 
was the direction to pursue for applications outside Probability in 
Banach spaces was delayed until very recently.  It does not seem 
possible to prove either (4.1.3), or even some of its most interesting 
consequences we will present in Part II through the martingale method.  
This should not be so surprising, since the inequalities of the 
present paper have been developed precisely to achieve what 
martingales seem unable to attain.  Among the results of Chapters 2 
to 5, apparently only those of Sections 2.1, 2.2 can be obtained using 
martingales; and the only reason why these are included here is that 
they provide an excellent and very simple setting to introduce our 
basic scheme of proof.  A major thesis of the present paper is that, 
while in principle the martingale method has a wider range of 
applications, in many situations the abstract inequalities of Part I 
are not only more powerful, but require considerably less ingenuity 
to apply.  In all the examples we examined, only in some rare 
situations, where the martingale is close to a sum of independent 
r.v., and where the value of numerical constants is crucial (such as 
[M-H]) did our methods fail to supersede martingales.

We now comment on the method of proof of the inequalities of Part I.  
Isoperimetric inequalities such as (1.5) or (1.7) are often proved 
via rearrangements.  That is, one produces a (simple if possible) 
way to transform the set $A$ in a set $T(A)$, of the same measure, 
but more regular, so that the measure of $T(A)_t$ is not more than 
the measure of $A_t$.  The procedure is then iterated, in a way 
that the iterates of $A$ converge to the ``extremal case''.  
Rearrangements are the only known technique to obtain perfect 
inequalities such as (1.5), (1.6).  The inequality of [T2], that 
started the present line of work was proved using rearrangements.  
The difficult proof requires different types of transformations, some 
of which prevent from obtaining the external sets.

Despite considerable efforts, rearrangements did not yield a proof of 
the inequality of [T1].  (As pointed out to me by N. Alon, the reason 
could be the complicated nature of the extremal sets.)  A completely 
new method was developed in [T2].  The main discovery there was that 
of a formulation that allows an easy proof by induction upon the 
number of coordinates.  The wide applicability of the method became 
apparent only gradually.  This method and its variations provide a 
unified scheme of proof of all our inequalities, that, in its simplest 
occurrence, is described in great detail in Section 2.1.  Ironically 
enough, this method is, in its principle, rather similar to the
martingale method; the extra power is gained from the possibility of 
abstract manipulations in product spaces.  A considerable advantage 
of the method is that, proving the induction hypothesis reduces to 
proving certain statements involving only functions on $\Omega$.  At 
times this is extremely easy; sometimes it is a bit harder.  But 
certainly the nature of the statements that have to be decided is
such that they are bound to yield to sufficient effort.  What on the 
other hand, is not entirely clear, is why this simple procedure seems 
so miraculously sharp; in the situations where explicit computations 
of the best possible constants given by the method has been possible, 
these constants have proved very close to the optimal.  In the cases 
where only less precise estimates have been possible, these estimates 
appear nonetheless to capture, up to a constant, the exact order of 
what really happens, and this, in every single situation that has 
been investigated.

The paper has been written to be read without any knowledge of this 
author's previous work or of the topic in general.  For the sake of 
completeness, the only previous result of the author that has not 
been either vastly generalized or considerably simplified has been 
reproduced (as Theorem 4.2.4).  Significant effort has been made in 
writing the paper in an easily accessible form.  For example, it turns 
out in several situations that the simplest occurrence of a new 
principle is also the most frequently used.  In these cases, we have 
taken care to give a separate proof for this most important case.  
These (short) proofs also serve as an introduction to the more 
complicated proofs of subsequent more specialized results.

During the preparation of this paper, I asked a number of people 
whether they were aware of recent or potential uses of the martingale 
method.  I am pleased to thank D. Aldous, E. Bolthausen, A. Frieze, 
C. McDiarmid, B. Pittel, M. Steele, W. Szpankowski for their precious 
suggestions.  Special thanks are due to H. Kesten, who communicated 
to me preprints of his recent work on percolation [K].  Analysis of 
his results pointed the way to several of the major developments that 
are presented in the present paper.  The material of Chapter 5 was 
directly motivated by questions of G. Schechtman concerning 
the ``correct form'' of the concentration of measure on the symmetric 
group.  A. Frieze, S. Janson and J. Wehr most helpfully contributed 
to literally hundreds of improvements upon the easy version of this 
work.  Finally, it must be acknowledged that this paper would not have 
been written if Professor Milman had not, over the years, convinced 
this author of the central importance of the concentration of measure 
phenomenon.
\vfil\eject

\noindent {\bf 2.~~Control by one point}

{\bf 2.1.~~The basic principle}

Throughout the paper we will consider a probability space $(\Omega ,
\Sigma ,\mu )$ and the product $(\Omega^N,\mu^N)$.  The product 
probability $\mu^N$ will be denoted simply by $P$.

Consider a subset $A$ of $\Omega^N$.  For $x\in \Omega^N$, we measure 
how far $x$ is from $A$ by
$$\leqalignno{f(A,x) &= \min\{\card\{ i\le N\,;\,x_i\not= y_i\}\,;\,y
\in A\}\,. &(2.1.1)\cr}$$
This is simply the Hamming distance from $x$ to $A$.  The reason that 
we use a different notation is that at later stages, we will introduce 
different ways to measure how far $x$ is from $A$.  These ways will 
not necessarily arise from a distance.

It should be observed that the function $f(A,\cdot )$ need not be 
measurable even when $A$ is measurable.  This is the reason for the 
upper integral and outer probability in Proposition 2.1.1. below.  On 
the other hand, measurability questions are simply irrelevant in the 
study of inequalities.  Simple and standard approximation arguments 
show that none of the results of this paper would lose any power if 
one should assume that $\Omega$ is Polish, $\mu$ is a Borel measure, 
and that one studies only compact sets.  It would be distracting to 
devote space and energy to these routine considerations.  Therefore, 
we have felt that it would be better to simply ignore all 
measurability questions, and treat all sets and functions as if they 
were measurable.  The reader will keep in mind that in the sequel, 
when measurability problems do arise, certain integrals (resp.
probabilities) have to be replaced by upper integrals (resp. outer 
probabilities) just as in the statement of Proposition 2.1.1.   (The 
reader who desires to have a proof of our statements without 
measurability assumption should be warned that it does not work to 
try to extend the proofs we give by putting outer integrals rather 
integrals -- the reason being that Fubini theorem fails for outer 
integrals.  Rather one has to derive the general result from the 
special case of well behaved sets by approximation.)

\proclaim{Proposition 2.1.1}  For $t>0$, we have
$$\leqalignno{\int^\ast e^{tf(A,x)}dP(x) &\le {1 \over P(A)}\left(
{1 \over 2} +{e^t+e^{-t}\over 4}\right)^N &(2.1.2)\cr
&\le {1 \over P(A)} e^{t^2N/4}\,.\cr}$$

In particular,
$$\leqalignno{P^\ast (\{ f(A,\cdot )\ge k\}) &\le {1 \over P(A)} 
e^{-k^2/N}\,. &(2.1.3)\cr}$$
\endproclaim

As was pointed out in the introduction, the power of our approach 
largely rests upon the fact that it reduces the proof of an inequality 
in $\Omega^N$ such as (2.1.2) to the proof of a much simpler fact 
about functions on $\Omega$.  In the present case, the meat of 
Proposition 2.1.1 is as follows.

\proclaim{Lemma 2.1.2}  Consider a (measurable) function $g$ on 
$\Omega$.  Assume $0\le g\le 1$.  Then we have
$$\leqalignno{\int_\Omega \min \left( e^t,{1 \over g(\omega )}\right) 
d\mu (\omega)\int_\Omega g(\omega )d\mu (\omega ) &\le a(t) 
&(2.1.4)\cr}$$
where we have set $a(t) = \left({1 \over 2} +{e^t+e^{-t} \over 4}
\right)$.
\endproclaim

\demo{Proof}  If we replace $g$ by $\max (g,e^{-t})$, this does not 
change the first integral, but increases the second.  Thus it suffices 
to prove that if $e^{-t}\le g\le 1$, we have
$$\int_\Omega {1 \over g} d\mu \int_\Omega gd\mu \le a(t)\,.$$

Consider the convex set ${\Cal C}$ of measurable functions $g$ on 
$\Omega$ for which $e^{-t}\le g\le 1$.  On ${\Cal C}$, the functional 
$g\to\int_\Omega g^{-1}d\mu$ is convex.  On the subset ${\Cal C}_b$ 
of ${\Cal C}$ that consists of the functions with integral $b$, this 
functional attains its maximum on an extreme point.  There is no loss 
of generality to assume that $\mu$ has no atoms; then it is well known 
that an extreme point of ${\Cal C}$ takes only the values $e^{-t}$ and 
$1$.  Thereby it suffices to show that for $0\le u\le 1$ we have
$$(1-u+ue^t)(1-u+ue^{-t})\le a(t)\,.$$
But the left hand side is invariant by changing $u$ into $1-u$, so 
that the maximum is obtained at $u=1/2$, and is $a(t)$.\rbx
\enddemo

The proof of Proposition 2.1.1 goes by induction over $N$.  The case 
$N=1$ follows from the application of (2.1.4) to $g=1_A$.

Suppose now that the result has been proved for $N$, and let us prove 
it for $N+1$.  Consider $A\subset\Omega^{N+1}=\Omega^N\times\Omega$.  
For $\omega\in\Omega$, we set
$$\leqalignno{A(\omega ) &= \{ x\in\Omega^N\,;\,(x,\omega )\in A\}\,. 
&(2.1.5)\cr}$$
and
$$B=\{ x\in\Omega^N\,;\,\exists\omega\in\Omega\,,\, (x,\omega )\in A
\}\,.$$

With obvious notations, we have
$$f(A,(x,\omega ))\le f(A(\omega ),x)\,.$$
Indeed, if $y\in A(\omega )$, then $(y,\omega )\in A$, and the number 
of coordinates where $(y,\omega )$ and $(x,\omega )$ differ is the 
number of coordinates where $x$ and $y$ differ.  Thus, by induction 
hypothesis, we have
$$\leqalignno{\int_{\Omega^N}\exp (tf(A,(x,\omega )))dP(x) &\le {a(t
)^N \over P(A(\omega))}\,. &(2.1.6)\cr}$$

We also observe that
$$f(A,(x,\omega )) \le f(B,x)+1$$
so that, by induction hypothesis, we have
$$\int_{\Omega^N} e^{tf(A,(x,\omega ))} dP(x)\le {e^t a(t)^N \over 
P(B)}\,,$$
and combining with (2.1.6) we get
$$\int_{\Omega^N} e^{tf(A,(x,\omega ))}dP(x) \le a(t)^N\min\left({e^t 
\over P(B)}\,,\,{1 \over P(A(\omega ))}\right)\,.$$
Integrating in $\omega$, we have 
$$\int_{\Omega^{N+1}} e^{tf(A,(x,\omega ))}dP(x)d\mu (\omega )\le
a(t)^N\int_\Omega\min \left({e^t 
\over P(B)}\,,\,{1 \over P(A(\omega ))}\right) d\mu (\omega )\,.$$

To complete the induction, it suffices to show, by Fubini theorem, 
that
$$\int_\Omega \min\left({e^t \over P(B)}\,,\,{1 \over P(A(\omega ))}
\right) d\mu (\omega )\le {a(t) \over P\otimes\mu (A)} = {a(t) \over 
\int_\Omega P(A(\omega ))d\mu (\omega )}\,.$$
But this follows from (2.1.4) applied to the function $g(\omega )=P(A
(\omega ))/P(B)$.

We now finish the proof of Proposition 2.1.1.  We note that
$$a(t)=1+\sum_{n\ge 1} {t^{2n} \over 2(2n)!}\,.$$
Now $2(2n)!\ge 4^n n!$.  Indeed, this holds for $n=1$, $n=2$, while 
if $n+1\ge 4$, we have
$${(2n)! \over n!} = (n+1)\cdots (2n)\ge 4^n\,.$$
Thus
$$a(t)\le 1+\sum_{n\ge 1}t^{2n}/4^nn! = \exp (t^2/4)\,.$$
Finally, (2.1.3) follows from Chebyshev inequality
$$\eqalign{P(\{ f(A,\cdot )\ge k\} ) &\le e^{-tk}\int e^{tf(A,x)}dP
(x)\cr
&\le {1 \over P(A)} e^{-tk+Nt^2/4}\cr}$$
for $t=2k/N$.\rbx

\remark{Remark 2.1.3}  Consider a sequence $(a_i)_{i\le N}$ of 
positive numbers.  If we now replace (2.1.1) by
$$\leqalignno{f(A,x) &= \inf\{\sum\{ a_i\colon i\le N\,;\, x_i\not= 
y_i\}\colon y\in A\} &(2.1.7)\cr}$$
the proof of Proposition 2.1.1. shows that
$$\leqalignno{\int e^{tf(A,x)}dP(x) &\le {1 \over P(A)} e^{t^2
\sum_{i\le N} a^2_i/4} &(2.1.8)\cr}$$
and, by Chebyshev inequality
$$\leqalignno{P(\{ f(A,\cdot )\ge u\} ) &\le {1 \over P(A)} 
e^{-u^2/\sum_{i\le N} a^2_i}\,. &(2.1.9)\cr}$$

A number of inequalities presented in Chapters 2 to 5 allow extensions 
that parallel the way Remark 2.1.3 expands Proposition 2.1.1.  These 
extensions are immediate, and will not be stated.  It should be 
pointed out, on the other hand, that no gain of generality would be 
obtained in Proposition 2.1.1. by replacing the product $\Omega^N$, 
$\mu^N$ by a product $\prod\limits_{i\le N}\Omega_i$, 
$\bigotimes\limits_{i\le N}\mu_i$.  This comment also applies to many 
inequalities that we will subsequently prove.
\endremark

{\bf 2.2.~~Sharpening}

Having proved (2.1.2), it is natural to wonder whether this could be 
improved by allowing another type of dependence of the right-hand side 
as a function of $P(A)$.  The most obvious choice is to replace $P(A
)^{-1}$ by $P(A)^{-\alpha}$ for some $\alpha >0$.

\proclaim{Proposition 2.2.1}  For $t\ge 0$, we have
$$\leqalignno{\int e^{tf(A,x)} dP(x) &\le {a(\alpha ,t)^N \over P
(A)^\alpha} &(2.2.1)\cr}$$
where
$$\leqalignno{a(\alpha ,t) &={\alpha^\alpha \over (\alpha +1)^{\alpha
+1}}~{(e^t-e^{-t/\alpha}
)^{1+\alpha} \over (1-e^{-t/\alpha})(e^t-1)^\alpha}\,. &(2.2.2)\cr}$$
\endproclaim

\demo{Proof}  Following the scheme of proof of Proposition 2.1.1, 
(2.2.1) holds provided, for each function $0\le g\le 1$ on $\Omega$, 
we have
$$\int_\Omega \min \left( e^t\,,\,{1 \over g^\alpha}\right) d\mu 
\left(\int_\Omega gd\mu \right)^\alpha\le a(\alpha ,t)\,.$$
Following the proof of Lemma 2.1.2, we see that we can take
$$\leqalignno{a(\alpha ,t) &=\sup_{0\le u\le 1}(1+u(e^t-1))(1-u(1
-e^{-t/\alpha}))^\alpha\,, &(2.2.3)\cr}$$
from which (2.2.2) follows by calculus.\rbx
\enddemo

Certainly neither the author nor the reader are enthusiastic about the 
prospect of using (2.2.1) and optimizing in Chebyshev inequality.  The 
purpose of the next result is to obtain a more manageable bound, that 
also makes clearer the gain obtained by taking large values of $\alpha$.

\proclaim{Lemma 2.2.2}
$$a(\alpha ,t)\le \exp {t^2 \over 8}\left( 1+{1 \over \alpha}\right)$$
\endproclaim

\demo{Proof}  Interestingly, rather than using (2.2.2), it seems 
simpler to go back to (2.2.3) and to show that, whenever $0\le u\le 
1$, we have
$$(1+u(e^t-1))(1-u(1-e^{-t/\alpha}))^\alpha\le\exp{t^2 \over 8}\left( 
1+{1 \over \alpha}\right)\,,$$
or, equivalently
$$\leqalignno{\log (1+u(e^t-1))+\alpha\log (1-u(1-e^{-t/\alpha})) 
&\le {t^2 \over 8}\left( 1+ {1 \over \alpha}\right)\,. &(2.2.4)\cr}$$

Since (2.2.4) holds for $t=0$, it suffices to show that the derivative 
of the left-hand side is bounded by the derivative of the right-hand 
side for $t\ge 0$, i.e.,
$$t\ge 0\Rightarrow {ue^t \over 1+u(e^t-1)} - {ue^{-t/\alpha} \over
1-u(1-e^{-t/\alpha})}\le {t \over 4}\left( 1+{1 \over \alpha}
\right)\,,$$
or, equivalently
$$\leqalignno{t\ge 0 &\Rightarrow {u-1 \over 1+u(e^t-1)} - {u-1 \over 
1-u(1-e^{-t/\alpha})}\le {t \over 4}\left( 1+{1 \over \alpha}
\right)\,. &(2.2.5)\cr}$$

Again (2.2.5) holds for $t=0$.  So it suffices to show that for $t\ge 
0$, the derivative of the left-hand side of (2.2.5) is bounded by the 
derivative of the right-hand side; or, equivalently, that
$$u(1-u)\left[{e^t \over (1-u+ue^t)^2} +{1 \over \alpha}~{e^{-t/
\alpha} \over (1-u+ue^{-t/\alpha})^2}\right]\le{1 \over 4}+{1 \over 
4\alpha}\,.$$
Now, using the inequality $4ab\le (a+b)^2$, we see that
$$\eqalignno{{u(1-u)e^t \over (1-u+ue^t)^2} &\le {1 \over 4}\,;
\qquad{u(1-u)e^{-t/\alpha} \over (1-u+ue^{-t/\alpha})^2}\le {1 
\over 4}\,. &\bx\cr}$$
\enddemo

\proclaim{Corollary 2.2.3}  For $t\ge 0$, we have
$$\leqalignno {\int e^{tf(A,x)}dP(x) &\le {1 \over P(A)^\alpha} \exp 
N {t^2 \over 8}\left( 1+{1 \over \alpha}\right)\,. &(2.2.6)\cr}$$
In particular, for $k\ge \sqrt{{N \over 2}\log {1 \over P(A)}}$, we 
have
$$\leqalignno{P(\{ f(A,\cdot )\ge k\} ) &\le \exp\left( -{2 \over N}
\left( k-\sqrt{{N \over 2}\log {1 \over P(A)}}\right)^2\right)\,. 
&(2.2.7)\cr}$$
\endproclaim

\demo{Proof}  Certainly (2.2.6) follows from (2.2.1) and Lemma 2.2.2.  
Optimization over $t$ in Chebyshev inequality yields
$$P(\{ f(A,\cdot )\ge k\} )\le {1 \over P(A)^\alpha} \exp\left( -{2
k^2 \over N}~{\alpha \over \alpha +1}\right)\,.$$
For $k\ge\sqrt{{N \over 2}\log{1 \over P(A)}}$, making the (optimal) 
choice
$$\alpha =-1+\sqrt{{2k^2 \over N\log{1 \over P(A)}}}$$
yields (2.2.7).\rbx
\enddemo

It is an interesting fact that (2.2.7) is exactly the best bound that 
has been proved on $P(\{f(A,\cdot )\ge k\} )$ using martingales (see 
[McD]).  It is a natural question to wonder whether, when 
$P(A)\ge 1/2$, one indeed has
$$P(\{ f(A,\cdot )\ge k\} )\le K\exp \left( -{2k^2 \over N}\right)$$
for some universal constant $K$.  More or less standard arguments 
(e.g., those contained in [T2]) show that it suffices to consider the 
case where $\Omega =\{ 0,1\}$, where $P$ is the product of measures 
$(\mu_i)_{i\le N}$ on $\Omega$, and where $A$ is even ``hereditary''.  
The case where $\mu_i(\{ 1\})=1/2$ for each $i\le N$ is known, as a 
consequence of more precise results, such as Harper's inequality.  
Intuitively, this is the worst case.

Having obtained (2.2.6), one must wonder whether further improvements 
upon (2.2.6) are possible by considering yet other general 
dependencies of the right-hand side as a function of $P(A)$. The 
reader who wishes to truly penetrate this paper will convince himself 
that this is not the case.

{\bf 2.3.~~Two point space}

Let us now consider the case where $\Omega =\{ 0,1\}$, and set $p=\mu 
(\{ 1\})$, so that $\mu (\{ 0\})=1-p$.

\proclaim{Proposition 2.3.1}  For $t\ge 0$, $\alpha \ge 1$, we have
$$\leqalignno{\int e^{tf(A,x)}dP(x) &\le {b(\alpha ,t,p)^N \over P
(A)^{\alpha}}\,, &(2.3.1)\cr}$$
where, for $p\ge 1/2$, we have set
$$\leqalignno{b(\alpha ,t,p) &= ((1-p)e^t+p)(p+(1-p)e^{-t/\alpha}
)^\alpha\,, &(2.3.2)\cr}$$
and, for $p\le 1/2$,
$$\leqalignno{b(\alpha, t,p) &= b(\alpha,t,1-p)=((1-p)e^{-t}+p)(p+
(1-p)e^{t/\alpha})^\alpha\,. &(2.3.3)\cr}$$
\endproclaim

\demo{Proof}  Following the proofs of Propositions 2.1.1. and 2.2.1 
it suffices to show that for any function $0\le g\le 1$ on $\Omega$ 
we have
$$\int_\Omega \min\left( e^t,{1 \over g^\alpha}\right)d\mu \left(\int 
g\,d\mu \right)^\alpha\le b(\alpha ,t,p)\,.$$

As in the proof of Lemma 2.1.2., we reduce to the case where $g\ge 
e^{-t/\alpha}$.  Setting $a=g(0)$, $b=g(1)$, it suffices to show 
that, for $e^{-t/\alpha}\le a$, $b\le 1$ we have
$$\left( (1-p){1 \over a^\alpha} +{p \over b^\alpha}\right)((1-p)a+p
b)^\alpha\le b(\alpha ,t,p)\,.$$
Setting $x=b/a$, it suffices to show that
$$e^{-t/\alpha}\le x\le e^{t/\alpha}\Rightarrow \varphi (x)\le b(
\alpha ,t,p)$$
where we have set
$$\varphi (x) = ((1-p)x^\alpha +p)\left({1-p \over x}+p
\right)^{\alpha}\,.$$

Now,
$$\varphi '(x)=\alpha p(1-p)\left( x^{\alpha -1}-{1 \over x^2}
\right)\left({1-p \over x} +p \right)^{\alpha -1}$$
so that $\varphi$ decreases for $x\le 1$, increases for $x\ge 1$.

Also, we have
$$\varphi '(x)-\left(\varphi \left({1 \over x}\right)\right) ' = 
\alpha p(1-p)\left( 1-{1 \over x^{\alpha +1}}\right)((1-p)+px
)^{\alpha -1} - ((1-p)x+p)^{\alpha -1})\,,$$
so that, for $x\ge 1$, this has the sign of $2p-1$.  Thus for $p\le 
1/2$, $\varphi$ attains its maximum on the interval $[e^{-t/\alpha},
e^{t/\alpha}]$ at the right end of this interval, while for 
$p\ge 1/2$ it attains its maximum at the left end.  (One should 
observe that changing $x$ in $1/x$ and $p$ in $1-p$ leave $\varphi$ 
invariant.)\rbx
\enddemo

A particularly important example is when
$$A=\left\{ x=(x_i)\in\{ 0,1\}^N\,;\,\sum_{i\le N}x_i\le k
\right\}\,.$$
The use of (2.3.1) for this set and of Chebyshev inequality will in 
particular produce bounds for the tails of the binomial law.  Thereby, 
it is not surprising that the computations involved in the use of 
(2.3.1) do run into the same type of difficulties as those involving 
the tails of the binomial law.  We now show how, nonetheless, some 
simple and reasonably sharp results can be deduced (for general sets 
$A$) from (2.3.1).  The reader will observe that the bound (2.3.1) is 
(of course) invariant when $p$ is replaced by $1-p$, so that there is 
no loss of generality to assume $p\ge 1/2$.  Let us fix $p$, $\alpha
\ge 1$, and consider
$$f(t)=\log b(\alpha ,t,p)=\log ((1-p)e^t+p)+\alpha\log (p+(1-p)e^{-t/
\alpha})\,.$$
Thus $f(0)=0$, and
$$f'(t)=(1-p)\left({1 \over (1-p)+pe^{-t}} - {1 \over pe^{t/\alpha}+
(1-p)}\right)\,.$$
Thus $f'(0)=0$, and
$$f''(t)=p(1-p)(h(e^{-t})+{1 \over \alpha} h(e^{t/\alpha}))$$
where
$$h(x)={x \over (px+1-p)^2} = {x \over (1-(1-x)p)^2}\,.$$

Simple computations show that when $x\ge 1/e$, we have $\vert h(x)-1
\vert\le K\vert x-1\vert$ for some universal constant $K$.  It follows 
that
$$t\le 1\Rightarrow f''(t)\le p(1-p)\left(\left( 1+{1 \over \alpha}
\right) +4Kt\right)$$
and, by integration, that
$$t\le 1\Rightarrow f(t)\le p(1-p)\left(\left( 1+{1 \over \alpha}
\right) {t^2 \over 2}+Kt^3\right)\,.$$

Thus, we have shown the first half of the following.

\proclaim{Corollary 2.3.2}  For $\alpha\ge 1$, $0\le t\le 1$, we have 
$$\leqalignno{\int e^{tf(A,x)}dP(x) &\le {1 \over P(A)^\alpha}\exp N
\left[ p(1-p)\left( 1+{1\over \alpha}\right) {t^2 \over 2}+Kt^3
\right]\,. &(2.3.4)\cr}$$
In particular, for
$$\left( 4p(1-p)N\log {1 \over P(A)}\right)^{1/2}\le k\le p(1-p)N$$
we have
$$\leqalignno{&~~~P(\{ f(A,x)\ge k\} ) &(2.3.5)\cr
&~~~~~\le \exp\left( -{1 \over 2p(1-p)N}\left( k-\sqrt {2p(1-p)N\log 
{1 \over P(A)}}\right)^2 +{K k^3 \over (p(1-p))^3 N^2}\right)\,.\cr}$$
\endproclaim

To obtain (2.3.5), one proceeds as in the proof of (2.2.7), using 
first Chebyshev inequality for $t={k \over p(1-p)}~{\alpha \over (1+
\alpha )N}$, then taking
$$\alpha =-1+\sqrt{{k^2 \over 2p(1-p)N\log {1 \over P(A)}}}\,.$$

It is of interest to compare the bound (2.3.5) with the isoperimetric 
inequalities obtained in [Lea]; these isoperimetric inequalities are 
optimal, but apply only to special sets (the so called hereditary 
sets).  The bound (2.3.5) is more general, and provides estimates of 
essentially the same quality.

We now turn to a rather different situation.  Beside the measure 
$\mu$, we consider on $\Omega$ another probability $\mu_1$, with $p_1
=\mu_1 (\{ 1\})>p$, and we set $P_1=\mu^N_1$.

\proclaim{Theorem 2.3.4}  For a subset $A$ of $\Omega^N$, and $x\in 
\Omega^N$, we consider
$$f(A,x)=\min\{\card\{ i\le N\,;\,x_i=1\,,\, y_i=0\}\,;\,y\in A\}\,.$$
Then, for $t\ge 0$,
$$\leqalignno{\int e^{tf(A,x)}dP(x) &\le {a(\alpha ,t)^N \over P_1(
A)^\alpha} &(2.3.6)\cr}$$
where 
$$a(\alpha ,t)=\max (1,(1-p+pe^t)(p_1e^{-t/\alpha}+1-p_1)^\alpha )
\,.$$
\endproclaim

\remark{Comment}  The really new phenomenon here is that for small 
$t$, one has $a(\alpha ,t)=1$.  In particular, if $\alpha =1$, this 
occurs whenever $e^t\le p_1(1-p)/p(1-p_1)$ so that one has
$$\leqalignno{\int\left({p_1(1-p) \over p(1-p_1)}\right)^{f(A,x)}dP(x)
&\le {1 \over P_1(A)}\,. &(2.3.7)\cr}$$

The remarkable feature about this statement is that it is independent 
of $N$ (and so is in essence an infinite dimensional statement).  This 
is the first of the results we present that apparently cannot be 
obtained via martingales (so it deserves to be called a theorem rather 
than a proposition).  The reader that would like to gain intuition 
about the phenomenon captured by Theorem 2.3.4 should consider the 
case where $A=\{ x\in\{ 0,1\}^N\,;\,\sum\limits_{i\le N}x_i\le n\}$.  
In order to have $P_1(A)$ of order $1/2$, one takes $n$ equal to 
$Np_1$, assuming for simplicity that this number is an integer.  
Observing that
$$f(A,x)>k\Leftrightarrow \sum_{i\le N} x_i> n+k=Np_1+k=Np+(k+N(p_1-p)
)$$
the quantity $P(\{ f(A,x)>k\})$ can be estimated through the tails of 
the binomial law; the most interesting values of $N$ are such that 
$N(p_1-p)\sim k$.
\endremark

The induction scheme of Proposition 2.1.1. will reduce Theorem 2.3.4. 
to an elementary two-point inequality, that is the object of the next 
lemma.

\proclaim{Lemma 2.3.5}  If $a\le b\le 1$, we have
$$\leqalignno{{1-p \over b^\alpha}+p\min\left({1 \over a^\alpha}\,,
\,{e^t \over b^\alpha}\right) &\le {a(\alpha ,t) \over (ap_1+b(1-p_1)
)^\alpha}\,. &(2.3.8)\cr}$$
\endproclaim

\demo{Proof}  If we set $x=\min\left({b \over a}\,,\,e^{t/\alpha}
\right)$, we are reduced to show that
$$1\le x\le e^{t/\alpha}\Rightarrow \varphi (x)\le a(\alpha ,t)$$
where
$$\varphi (x) = (1-p+px^\alpha )\left({p_1 \over x}+(1-p_1)
\right)^\alpha\le a(\alpha ,t)\,.$$
But $\varphi '(x)$ has the sign of $p(1-p_1)x^{\alpha +1}-p_1(1-p)$, 
so it is negative for values of $x$ close to one, and then, possibly, 
becomes positive.  Thus $\varphi$ attains its maximum on the interval 
$[1,e^{t/\alpha}]$ at one of the endpoints.\rbx
\enddemo

\demo{Proof of Theorem 2.3.4}  We proceed by induction over $N$.  For 
$N=1$, since $f(A,\omega )\equiv 0$ when $1\in A$, it suffices to 
consider the case $A=\{ 0\}$, in which case the result follows 
from (2.3.8) with $a=0$, $b=1$.

Assuming now that the theorem has been proved for $N$, we prove it for 
$N+1$.  Consider $A\subset \Omega^{N+1}$, and set $A_1=\{x\in\Omega^N
\,;\,(x,1)\in A\}$.  Consider the projection $B$ of $A$ on $\Omega^N$.  
We observe that
$$\eqalign{f(A,(x,\omega )) &\le 1+f(B,x)\cr
f(A,(x,\omega )) &\le f(A_1,x)\cr}$$
so that setting $a=P_1(A_1)$, $b=P_1(B)$ and using the induction 
hypothesis, the result follows from (2.3.8).\rbx
\enddemo

{\bf 2.4.~~Penalties, I.}

A (somewhat imprecise) way to reformulate (2.1.1) is that we measure 
how far $x$ is from $A$ by simply counting the smallest number of 
coordinates of $x$ that cannot be captured by a point of $A$.  
Rather than just giving a penalty of $1$ for each coordinate we miss, 
it is natural to consider, given a non-negative function $h$ on 
$\Omega\times\Omega$, the quantity
$$\leqalignno{f_h(A,x) &= \inf\left\{\sum_{i\le N}h(x_i,y_i)1_{\{x_i
\not= y_i\}}\,;\,y\in A \right\}\,. &(2.4.1)\cr}$$

To simplify the notations, we will assume
$$\leqalignno{\forall x\in\Omega\,,\qquad h(x,x)=0 &&(2.4.2)\cr}$$
so that (2.4.1) becomes
$$\leqalignno{f_h(A,x) &= \inf\left\{\sum_{i\le N}h(x_i,y_i)\,;\,y
\in A\right\}\,. &(2.4.3)\cr}$$

Concerning (2.4.2), we should point out that we will let $x$, $y$ 
denote points in $\Omega^N$ as well as points in $\Omega$; when there 
is too much danger of confusion, however, points of $\Omega$ will be 
denoted by $\omega$, $\omega '$.

The function $h$ will always be assumed to be measurable.  The 
following simple result is already useful, as will be demonstrated in 
Chapter 11.

\proclaim{Theorem 2.4.1}  For each measurable subset $A$ of $\Omega^N$, 
and each $t>0$ for which\hfil\break $\iint\exp th(x,y) d\mu (x)d\mu 
(y)<\infty$, we have, setting $v(\omega ,\omega ')=\max (h(\omega ,
\omega '),h(\omega ',\omega ))$, that
$$\leqalignno{\int_{\Omega^N}e^{tf_h(A,x)}dP(x) &\le {1 \over P(A)}
\left({1 \over 2}\int_{\Omega^2}(e^{tv(\omega ,\omega ')}+e^{-tv(
\omega ,\omega ')}) d\mu (\omega )d\mu (\omega ')\right)^N\,. 
&(2.4.4)\cr}$$
\endproclaim

The crucial point of Theorem 2.4.1 is as follows.

\proclaim{Proposition 2.4.2}  Consider a function $g\ge 0$ on $\Omega$, 
and set
$$\leqalignno{\widehat{g}(x) &=\inf_{y\in\Omega}(g(y)+th(x,y))\,. 
&(2.4.5)\cr}$$
Then
$$\leqalignno{\int^\ast e^{\widehat g}d\mu\int e^{-g}d\mu &\le {1 
\over 2}\int_{\Omega^2} (e^{tv(\omega ,\omega ')}+e^{-tv(\omega ,
\omega ')})d\mu (\omega )d\mu (\omega ')\,. &(2.4.6)\cr}$$
\endproclaim

Indeed, a simple truncation argument shows that Proposition 2.4.2 
remains true if one allows (using obvious conventions) $g$ to take 
values in $\Bbb R^+\cup\{\infty\}$.  To prove Theorem 2.4.1. by 
induction over $N$, considering a subset $A$ of $\Omega^{N+1}$, for
$\omega\in\Omega$ we set
$$A(\omega )=\{x\in\Omega^N\,;\, (x,\omega )\in A\}\,,$$
and we define $g$ by $P(A(\omega ))=e^{-g(\omega )}$.  It should then 
be clear that (2.4.6) is exactly what is needed to make the induction 
work.

\demo{Proof of Proposition 2.4.2}  For simplicity we assume $\widehat 
g$ measurable.  Then the left-hand side of (2.4.6) coincides with
$$\iint_{\Omega^2}e^{\widehat{g}(x)-g(y)}d\mu (x)d\mu (y)\,.$$
We set $u(x,y)=\widehat{g}(x)-g(y)$.  By definition of $\widehat{g}$, 
we have $\widehat{g}(x)\le g(y)+th(x,y)$.  Since $h(x,x)=0$, we also 
have $\widehat{g}(x)\le g(x)$.  Hence
$$\leqalignno{u(x,y)\le th(x,y)\,;\qquad u(x,y)\le g(x)-g(y) 
&&(2.4.7)\cr}$$

We now observe that for two numbers $a,b$, if $a+b\le 0$, then
$$e^a+e^b\le e^{\max (a,b,0)}+e^{-\max (a,b,0)}\,.$$
Thereby, by (2.4.7), we have
$$e^{u(x,y)}+e^{u(y,x)}\le e^{tv(x,y)}+e^{-tv(x,y)}\,.$$
The result follows by integration.\bx
\enddemo

It is of interest to get simpler bounds.  Let us observe the following 
elementary fact (that is obvious on power series expansions)
$$\leqalignno{\text{The~function}~x^{-2}(e^x+e^{-x}-2)~
\text{increases~for}~ x\ge 0\,. &&(2.4.8)\cr}$$

Thus, for $t\le 1$
$${e^{tv} +e^{-tv}-2 \over t^2v^2} \le {e^v+e^{-v}-2 \over v^2}$$
and hence
$$\leqalignno{{e^{tv}+e^{-tv} \over 2} &\le 1+{t^2 \over 2}( e^v +
e^{-v}-2)\,. &(2.4.9)\cr}$$

We note that, for an increasing function $\varphi$,
$$\varphi (\max (a,b))\le \max (\varphi (a),\varphi (b))\le \varphi 
(a) +\varphi (b)\,.$$

Using this for $\varphi (x)=e^x+e^{-x}-2$, $a=h(\omega ,\omega ')$, 
$b=h(\omega ',\omega )$, using then (2.4.9) and integrating, we get 
the following from (2.4.4).

\proclaim{Theorem 2.4.3} If
$$\leqalignno{\iint \exp h(x,y)d\mu (x)d\mu (y) &<\infty\,, 
&(2.4.10)\cr}$$
we have for $t\le 1$,
$$\leqalignno{&\int_{\Omega^N}e^{tf_h(A,x)}dP(x) &(2.4.11)\cr
&\qquad \le {1 \over P(A)}\exp\left( Nt^2\iint (e^{h(\omega ,
\omega ')}+e^{-h(\omega ,\omega ')}-2)d\mu (\omega )d\mu (\omega ')
\right)\,. \cr}$$
\endproclaim

\proclaim{Corollary 2.4.4}  Assume
$$\leqalignno{\iint \exp h(x,y)d\mu (x)d\mu (y) &\le 2 &(2.4.12)\cr}$$
Then for all $u\le 2N$ we have
$$\leqalignno{P(\{ f_h(A,\cdot )\ge u\} ) &\le {1 \over P(A)} e^{-u^2
/4N}\,. &(2.4.13)\cr}$$
\endproclaim

\demo{Proof}  Since $e^{-h}\le 1$, under (2.4.12), the right-hand side 
of (2.4.11) becomes bounded by $P(A)^{-1}\exp Nt^2$, from which 
(2.4.13) follows by Chebyshev inequality.\bx
\enddemo

The following resembles Bernstein's inequality.

\proclaim{Corollary 2.4.5}  Assume that $\Vert h\Vert_\infty =
\sup\limits_{x,y\in\Omega^2} h(x,y)$ 
is finite.  Then
$$\leqalignno{P(\{f_h(A,\cdot )\ge u\} ) &\le {1 \over P(A)}\exp 
\left( -\min\left({u^2 \over 8N\Vert h\Vert^2_2}\,,\,{u \over 2\Vert 
h\Vert_\infty}\right)\right) &(2.4.14)\cr}$$
where we have set $\Vert h\Vert_2 = (\iint_{\Omega^2} h^2(\omega ,
\omega ')d\mu (\omega )d\mu (\omega '))^{1/2}$.
\endproclaim

\demo{Proof}  By homogeneity, we can replace $h$ by $h'=h/\Vert h
\Vert_\infty$.  For $x\le 1$, by (2.4.8), we have $e^x+e^{-x}-2\le 
x^2(e+e^{-1}-2)\le 2x^2$.  Thereby the right hand side of (2.4.11) 
becomes bounded by $P(A)^{-1}\exp 2Nt^2\Vert h\Vert^2_2$, from which 
the result follows by Chebyshev inequality.
\enddemo

\remark{Remark}  The reader has possibly observed that we have made no 
special efforts to get sharp numerical constants (in contrast with the 
previous sections) and we have used the simplest estimates, however 
crude.  This feature will occur repeatedly.  For a number of the 
results we will present, the proofs do not seem adapted to the 
obtention of sharp constants.  Thereby, there is actually no point to 
track explicit values of the numerical constants involved.  Throughout 
the paper, $K$ will denote a universal constant, that may vary at 
each occurrence.
\endremark

{\bf 2.5.~~Penalties, II.}

It should be apparent from (2.4.1) that $f_h$ depends on $h$ only 
through the properties of the 
following functional, defined for subsets $B$ of $\Omega$
$$\leqalignno{h(\omega ,B) &= \inf \{ h(\omega ,\omega ')\,;\,\omega '
\in B\}\,. &(2.5.1)\cr}$$
(The reader should carefully compare this definition with (2.4.3) and 
note that in both cases the infimum is taken over the second variable.)

Thereby, one should expect that the exponential integrability of $h$ 
can be replaced in Theorem 2.4.1.  By a weaker condition on the 
functional $h(x,B)$.  This is indeed the case.

\proclaim{Theorem 2.5.1}  Assume that for each subset $B$ of $\Omega$ 
we have
$$\leqalignno{\int_\Omega \exp 2h(x,B)d\mu (x) &\le {e \over \mu (B)}
\,. &(2.5.2)\cr}$$
Then, for each subset $A$ of $\Omega^N$, and each $0\le t\le 1$, we 
have
$$\leqalignno{\int_{\Omega^N}e^{tf_h(A,x)}dP(x) &\le {e^{3t^2N} \over 
P(A)}\,. &(2.5.3)\cr}$$
\endproclaim

\remark{Discussion}  1)  It is good to observe and keep in mind that 
by H\"older's inequality, we have for $a\le 1$
$$\leqalignno{\int e^{ah}d\mu &\le (\int e^hd\mu )^a\,. &(2.5.4)\cr}$$
Thus, the precise value of constants such as the constants $2$, $e$ 
that occur in (2.5.2) is completely irrelevant.  Actually we will use 
the following consequence of (2.5.2):
$$\leqalignno{\int_\Omega \exp h(x,B)d\mu (x) &\le {\sqrt {e} \over 
\sqrt{\mu (B)}}\le {2 \over \sqrt{\mu (B)}}\,. &(2.5.5)\cr}$$ 

2)  It is very instructive to compare (2.5.2) with a condition such 
as (2.4.10).  Indeed, under (2.4.10), we have for all $x$
$$\leqalignno{\mu (B)\exp h(x,B) &= \mu (B)\inf_{y\in B}\exp h(x,y)
\le \int_{\Omega}\exp h(x,y)d\mu (y) \,. &(2.5.6)\cr}$$

Integrating in $x$ gives
$$\int \exp h(x,B)d\mu (x)\le {1 \over \mu (B)}\iint_{\Omega^2}\exp 
h(x,y)d\mu (x)d\mu (y)\,.$$
Thus (with the exception of the largely irrelevant factor 2), (2.4.10) 
appear stronger than (2.5.2).  It is indeed much stronger, a fact that 
is not surprising in view of the crudeness of (2.5.6).  To get a 
concrete example, consider the case where $\Omega$ is itself a product 
of $m$ spaces (and $\mu$ a product measure), and denote by $f(x,y)$ 
the Hamming distance in $\Omega$.  Then Proposition 2.1.1. asserts 
that the function $h=m^{-1/2}f$ satisfies (2.5.2).  On the other hand 
(except in trivial cases) the function $f/a$ will fail (2.4.12) unless 
$a$ is of order $m$.
\endremark

To prove Theorem 2.5.1, the induction method reduces to the proof of 
the following.

\proclaim{Proposition 2.5.2}  Consider $0\le t\le 1$, and a function 
$g\ge 0$ on $\Omega$.  For $s\ge 0$, we set $B_s=\{ g\le s\}$, and we 
consider
$$\leqalignno{\widehat{g} (x) &= \inf_{s>0}s+th(x,B_s)\,. 
&(2.5.7)\cr}$$
Then under (2.5.2) we have
$$\leqalignno{\int^\ast e^{\widehat{g}}d\mu\int e^{-g}d\mu &\le 
e^{3t^2}\,. &(2.5.8)\cr}$$
\endproclaim

\demo{Proof}  We observe that
$$\leqalignno{\widehat{g}(x) -g(y) &\le th(x,B_{g(y)})\,. 
&(2.5.9)\cr}$$

We then follow the argument of Proposition 2.4.2, using (2.5.9) rather 
than the first part of (2.4.7).  Combining with the argument of 
Theorem 2.4.3, we are led to show that
$$\iint_{\Omega^2}e^{h(x,B_{g(y)})}d\mu (x)d\mu (y)\le 4\,.$$
Using (2.5.5) and Fubini theorem, it suffices to show that
$$\leqalignno{\int_\Omega {1 \over \sqrt{\mu (B_{g(y)})}}d\mu (y) 
&\le 2\,. &(2.5.10)\cr}$$

The best way to prove this inequality is to observe that the left-hand 
side depends only on the function $s\to\mu (B_s)$.  Thus there is no 
loss of generality to assume that $\Omega =[0,1]$, that $\mu$ is 
Lebesgue measure, and that $g$ is nondecreasing.  But then $\mu (B_{g
(y)})\ge y$, and $\int^1_0 y^{-1/2}dy=2$.\rbx
\enddemo

As pointed out in the discussion, a natural application of Theorem 
2.5.1. is to the case where $\Omega$ is already a product space.  This 
will be used implicitly, but crucially in Section 11.5.  To formulate 
in words what happens, Proposition 2.1.1. states that if $A$ is a 
subset of a product $\Omega^N$ of $N$ spaces, of measure $1/2$, all 
but exceptional points $x$ of $\Omega^N$ are such that there is a 
point in $A$ that captures all but about $\sqrt{N}$ of their 
coordinates.  Suppose now that $N=N_1N_2$, and we think of the $N$ 
coordinates as $N_1$ blocks of $N_2$ coordinates.  Then, using Theorem 
2.5.1, we know that (for but exceptional points $x$) not only we will 
find a point in $A$ that misses only about $\sqrt{N}$ coordinates of 
$x$, but these coordinates will be concentrated in only about 
$\sqrt{N_1}$ blocks.  An interesting question would be to quantify 
precisely what can be said when, rather than considering only two 
``levels'', one considers a large number of levels.

{\bf 2.6.~~Penalties, III.}

In this section, we explore a new phenomenon, that will also be met 
in Sections 3.3.3 and 4.4.4.  The notations of the present section 
will be used throughout the paper.  Roughly speaking, what happens 
is that if, in (2.5.3), one allows a more general type of dependence 
on $P(A)$ of the right hand side, then a weaker condition than (2.5.2) 
will suffice; this will mean in practice weaker integrability 
requirements on $h$.

The dependence in $P(A)$ we will consider will be of the type 
$e^{\theta (P(A))}$.  Throughout the paper, $\theta$ will denote a 
convex decreasing function from $]0,1]$ to $\Bbb R^+$, such that
$\theta (1)=0$, $\lim\limits_{x\to 0}\theta (x)=\infty$.  The most 
important example is $\theta (x)=-\log x$, in which case $e^{\theta (
P(A))}$ is the familiar quantity $1/P(A)$.  We will always denote 
by $\xi$ the inverse function of $\theta$, so that $\xi$ is a convex 
function from $\Bbb R^+$ to $]0,1]$, with $\xi (0)=1$.  We will always 
assume the following
$$\leqalignno{\xi ''~\text{decreases}\,;\qquad \forall b>0\,,~~\xi 
''(b)\le\vert\xi '(b)\vert &&(2.6.1)\cr}$$
For $x\in\Bbb R$, we set $x^+=\max (x,0)$, and we will keep the 
following notation, for $x\in\Bbb R$, $b\in\Bbb R^+$
$$\leqalignno{\Xi (x,b) &= \xi (x^+)-\xi (b)-(x^+-b)\xi '(b) 
&(2.6.2)\cr}$$

We denote by $\lambda$ Lebesgue measure on $[0,1]$.  The measure of a 
Borel set $B$ is simply denoted by $\vert B\vert$.

Central to this section is the following technical condition, that 
relates $\xi$ and a function $w\ge 0$ defined on $[0,1]$.

\remark{Condition $H(\xi ,w)$}
$$\leqalignno{\forall b\ge 0\,,~\forall t\,,~0\le t\le 1\,,~
\int^1_0\Xi (b-tw(u),b)d\lambda (u) &\le t^2\vert\xi '(b)\vert\,. 
&(2.6.3)\cr}$$
First we will investigate conditions that imply (2.6.3) under two 
simple choices of $\xi$.  Then we will look at a rather general 
situation where the meaning of (2.6.3) can be considerably clarified; 
and before stating the main result (Theorem 2.6.5) we will prove a 
technical lemma that will explain the precise purpose of condition 
$H(\xi ,w)$.
\endremark

\proclaim{Proposition 2.6.1}  When $\xi (x)=e^{-x}$, condition $H(\xi 
,w)$ holds provided $\int e^wd\lambda\le 2$.
\endproclaim

\demo{Proof}  Indeed, we have
$$\eqalign{\Xi (x,b) &= e^{-x^+}-e^{-b}+(x^+-b)e^{-b}\cr
&\le e^{-x}-e^{-b} +(x-b)e^{-b}\cr
&= e^{-b}(e^{-(x-b)}+(x-b)-1)\,.\cr}$$
Thus (2.6.3) holds provided
$$t\le 1\Rightarrow \int (e^{tw}-tw-1)d\lambda \le t^2\,.$$
But, since the function $x^{-2}(e^{x}-x-1)$ increases for $x\ge 0$, 
we have $(e^{tw}-tw-1)\le t^2(e^w-w-1)$.\rbx
\enddemo

\proclaim{Proposition 2.6.2}  If $\xi (x)={1 \over x+1}$, then 
condition $H(\xi ,\omega )$ holds provided $\int w^2d\lambda\le 1$.
\endproclaim

\demo{Proof}  Setting $y=x^+$, we have
$$\eqalignno{\Xi (x,b) &= {1 \over y+1}-{1 \over b+1} +{y-b \over 
(b+1)^2}\cr
&= {(y-b)^2 \over (y+1)(b+1)^2} \le (y-b)^2\vert\xi '(b)\vert\,. 
&\bx\cr}$$
\enddemo

One obvious consequence of (2.6.3), taking $t=1$ is that
$$\leqalignno{\vert\{ w\ge b\}\vert \Xi (0,b) &\le \vert\xi '(b)
\vert\,. &(2.6.4)\cr}$$
In practice for $b$ large $\Xi (0,b)$ is of order $1$; so (2.6.4) is 
really a tail condition.  The next result shows that a condition of a 
similar nature is indeed sufficient, provided $\xi '$ varies smoothly 
(i.e., satisfies the $\Delta_{2^-}$ condition; which is not the case 
when $\xi (x)=e^{-x}$).

\proclaim{Proposition 2.6.3}  Assume that for a certain number $L>0$, 
we have
$$\leqalignno{\forall b>0\,,~\forall t\le 1\,,~\vert\xi '({b \over 2t} 
)\vert\le Lt^2\vert\xi '(b)\vert\,. &&(2.6.5)\cr}$$
Then (2.6.3) holds provided the following two conditions hold:
$$\leqalignno{\int w^2d\lambda \le {1 \over L} &&(2.6.6)\cr}$$
$$\leqalignno{\forall b>0\,,~\vert\{ w\ge b\}\vert\le {1 \over 2L}
\vert\xi '(b)\vert\,. &&(2.6.7)\cr}$$
\endproclaim

\demo{Proof}  We write
$$\leqalignno{\int\Xi (b-tw,b)d\lambda &\le \int_{\{tw\le b/2\}}\Xi 
(b-tw,b)d\lambda +\vert \{ tw\ge b/2\}\vert\,. &(2.6.8)\cr}$$
By Taylor's formula, since $\xi ''$ decreases, and $\xi ''(b)\le 
\vert\xi '(b)\vert$, we have, by (2.6.5)
$$\eqalign{x\ge b/2\Rightarrow\Xi (x,b) &\le {(x-b)^2 \over 2}\xi 
''\left({b \over 2}\right)\cr
&\le {(x-b)^2 \over 2}\left|\xi '\left({b \over 2}\right)\right|\cr
&\le L{(x-b)^2 \over 2} \vert\xi '(b)\vert\,.\cr}$$
Thus
$$\int_{\{ tw\le b/2\}}\Xi (b-tw,b)d\lambda\le {Lt^2 \over 2}\vert 
\xi '(b)\vert\int w^2d\lambda\,.$$

Also, by (2.6.7), (2.6.5)
$$\vert\{ tw\ge b/2\}\vert\le {1 \over 2L}\left|\xi '\left({b \over 
2t}\right)\right|\le {t^2 \over 2}\vert\xi '(b)\vert\,.$$
The result follows, combining with (2.6.8)\,.\rbx
\enddemo

The reader should observe that the functions $\xi (x)=(1+x)^{-\alpha}$ 
($\alpha\ge 1$) satisfy (2.6.5).

The following lemma explains the purpose of condition $H(\xi ,w)$.

\proclaim{Lemma 2.6.4}  Consider a function $f\ge 0$ on $\Omega$.  
Assume that for a certain $t$, $0\le t\le 1$ and all $s\le b$ we have 
$$\leqalignno{\mu (\{ f\le s\} )\le \vert \{ tw\ge b-s\}\vert &= 
\vert\{ b-tw\le s\}\vert\,. &(2.6.9)\cr}$$
Then under condition $H(\xi ,w)$, for each set $C$ we have
$$\leqalignno{&\int_C\xi (f)d\mu &(2.6.10)\cr 
&\quad \le \xi (b)\mu (C)+\xi '(b)\int_C (f-b)d\mu +t^2\vert\xi '
(b)\vert +{1 \over 2}\xi ''(b)\int_{C\cap\{ f\ge b\}}(f-b)^2d\mu\,
.\cr}$$
\endproclaim

\demo{Proof}  By definition of $\Xi$, (2.6.10) is equivalent to
$$\int_C\Xi (f,b)d\mu \le t^2\vert\xi '(b)\vert +{1 \over 2}\xi ''(b)
\int_{C\cap\{ f \ge b\}}(f-b)^2 d\mu\,.$$

By Taylor's formula, and since $\xi ''$ decreases, for $x>b$ we have
$$\Xi (x,b) \le {1 \over 2} (x-b)^2\xi ''(b)$$
and thus
$$\int_{C\cap\{ f\ge b\}}\Xi (f,b)d\mu \le {1 \over 2}\xi '' (b)
\int_{C\cap\{ f\ge b\}}(f-b)^2 d\mu\,.$$

If we remember that $\Xi\ge 0$, and if we use condition $H(\xi ,w)$, 
we then see that it suffices to show that
$$\leqalignno{\int_{\{f\le b\}}\Xi (f,b)d\mu &\le \int\Xi (b-tw,b)d
\lambda\,. &(2.6.11)\cr}$$

Now, (2.6.9) implies that for all $s<b$ we have
$$\mu (\{ f\le s\} ) \le \vert\{ b-tw\le s\}\vert\,.$$
Thus, since $\Xi (x,b)$ decreases for $x\le b$, we have, for all $z>0$
$$\mu (1_{\{ f\le b\}}\Xi (f,b)\ge z) \le \vert\{\Xi (b-tw,b)\ge z\}
\vert\,,$$
from which (2.6.11) follows.\rbx
\enddemo

\proclaim{Theorem 2.6.5}  Consider a function $h$ on $\Omega\times
\Omega$, and a nonincreasing function $w$ on $[0,1]$ such that $\int 
w^2d\lambda\le 1$.  Assume that for each subset $B$ of $\Omega$, we 
have
$$\leqalignno{\int_\Omega \exp h(x,B)d\mu (x) &\le \exp (w(\mu (B))) 
&(2.6.12)\cr}$$
where we keep the usual notation $h(x,B)=\inf\{ h(x,y)\,;\,y\in B\}$.  
Consider a function $\theta$ as usual, and assume that the condition 
$H(\zeta ,w)$ holds.

Then, for each subset $A$ of $\Omega^N$, and all $t\le 1$, we have, 
for all $t\le 1$
$$\int_{\Omega^N}e^{tf_h(A,x)}dP(x)\le \exp (4Nt^2+\theta (P(A)))\,.$$
\endproclaim

To understand better (2.6.12) it is of interest to specialize to the 
case where $h$ depends only on $x$ (resp. $y$).  If $h$ depends on $x$ 
only, (2.6.12) means that $\int_\Omega \exp h(x)d\mu (x)\le \exp w(0
)$.  If $h$ depends on $y$ only, then (2.6.12) becomes
$$\inf\{ h(y)\,;\,y\in B\}\le w(\mu (B))\,.$$
Taking $B=\{ h\ge s\}$, we get $s\le w(\mu (\{h\le s\} ))$ and, since 
$w$ in nonincreasing, this implies
$$\mu(\{ h\ge s\} )\le \vert\{ w\ge s\}\vert\,.$$
It is easy to see that, conversely, this implies (2.6.12) when $h$ 
depends upon $y$ only and when $w$ is left continuous.

To prove Theorem 2.6.5, it suffices, using induction over $N$, to 
prove the following.

\proclaim{Proposition 2.6.6}  Consider a function $g$ on $\Omega$, 
$0<g\le 1$, and set
$$\widehat{\theta} g(x) = \inf_{y\in\Omega}\{\theta (g(y))+th(x,y)\}
\,.$$
Then, under the conditions of Theorem 2.6.5, for $t\le 1$, we have
$$\int e^{\widehat{\theta}g}d\mu\le \exp (4t^2+\theta (\int gd\mu ))
\,.$$
\endproclaim

Clearly, this is equivalent to the following.

\proclaim{Proposition 2.6.7}  Consider a function $f$ on $\Omega$, $f
\ge 0$, and set
$$\widehat{f}(x)=\inf_{y\in\Omega}\{ f(y)+th(x,y)\}\,.$$
Then, under the conditions of Theorem 2.6.5, for $t\le 1$, we have
$$\leqalignno{\int e^{\widehat f} d\mu &\le \exp (5t^2+\theta (\int\xi 
(f)d\mu ))\,. &(2.6.13)\cr}$$
\endproclaim

\demo{Proof}  The problem is that we have on the right of (2.6.13) 
the quantity $\theta (\int\xi (f)d\mu)$ rather than the larger 
quantity $\int fd\mu$.  We consider $t$ as fixed through the proof.

{\bf Step 1.}  We set $B_s=\{ f\le s\}$ for $s\ge 0$, and
$$b = \inf_s\{s+tw(\mu (B_s))\}\,.$$
We note that $\widehat{f}(x)\le f(x)$.  We consider the function $f'$ 
given by
$$\eqalign{f'(x)=\widehat{f}(x)~\quad~ &\text{if}\quad \widehat{f}(x) 
>b\cr
f'(x)=b~\quad~ &\text{if}\quad \widehat{f}(x)\le b<f(x)\cr
f'(x)=f(x)~\quad~ &\text{if}\quad f(x)<b\,.\cr}$$
Since $\widehat{f}\le f$, we have $\widehat{f}\le f'\le f$.  Thus 
$\int e^{\widehat f}d\mu\le\int e^{f'}d\mu$, and $\xi (f)\le\xi (f')$, 
so $\int\xi (f) d\mu\le \int\xi (f')d\mu$ and $\theta (\int\xi (f')d
\mu )\le\theta (\int\xi (f)d\mu )$.  Thereby, it suffices to prove 
that
$$\leqalignno{\int e^{f'}d\mu &\le \exp (5t^2+\theta (\int\xi (f')d
\mu ))\,. &(2.6.14)\cr}$$

{\bf Step 2.}  By definition of $b$, for $s<b$, we have
$$tw(\mu (B_s))\ge b-s\,.$$
Since $w$ is nonincreasing, we have
$$\leqalignno{\vert\{ tw\ge b-s\}\vert &\ge\mu (B_s)\,. &(2.6.15)\cr}$$
Since $f'(x) =f(x)$ when $f(x)<b$, we see that (2.6.9) holds (for $f'$ 
rather than $f$).  Since $f'=\widehat f$ when $f'(x)>b$, by (2.6.10) 
used for $C=\Omega$, we get, since $\xi ''(b)\le \vert\xi '(b)\vert$
$$\leqalignno{\int\xi (f')d\mu &\le \xi (b)+\xi '(b)\int (f'-b)d\mu 
+\vert\xi '(b)\vert (t^2 +{1 \over 2}\int_{\{\widehat{f}>b\}}(
\widehat{f}-b)^2d\mu )\,. &(2.6.16)\cr}$$

{\bf Step 3.}  If $y\in B_s$, we have
$$\widehat{f}(x)\le f(y)+th(x,y)\le s+th(x,y)$$
so that
$$\widehat{f}(x)\le s+th(x,B_s)\,.$$
Thus, by (2.6.12) we have
$$\int e^{t^{-1}\widehat{f}}d\mu \le \exp (t^{-1}(s+tw(\mu (B_s))))
\,.$$
Taking the infimum over $s$ yields
$$\leqalignno{\int e^{t^{-1}(\widehat {f}-b)}d\mu &\le 1\,. 
&(2.6.17)\cr}$$
Since $e^{x^+}\le 1+e^x$, we get
$$\leqalignno{\int e^{t^{-1}(\widehat{f}-b)^+}d\mu &\le 2\,. 
&(2.6.18)\cr}$$

{\bf Step 4.}  The inequality $e^x\ge 1+x^2/2$ for $x\ge 0$, and 
(2.6.18) show that $\int ((\widehat{f}-b)^+)^2d\mu\le 2t^2$.  
Combining with (2.6.16), we get
$$\leqalignno{\int\xi (f')d\mu &\le \xi (b)+\xi '(b)\int (f'-b)d\mu 
+2t^2\vert\xi '(b)\vert \,. &(2.6.19)\cr}$$

The convexity of $\theta$ implies that $\theta (x)\ge\theta (y)+(x-y)
\theta '(y)$.  Also, since $\theta (\xi (x))=1$, we have $\theta '(
\xi (b))=1/\xi '(b)$.  Thus (2.6.19) implies
$$\leqalignno{\theta (\int\xi (f')d\mu ) &\ge b+\int (f'-b)d\mu 
-2t^2 &(2.6.20)\cr
&=\int f'd\mu -2t^2\,.\cr}$$

{\bf Step 5.}  To finish the proof, it is thereby sufficient to show 
that
$$\leqalignno{\int e^{f'}d\mu &\le \exp (2t^2+\int f'd\mu )\,. 
&(2.6.21)\cr}$$
Consider the function $r(x)=e^x-x-1$, so that
$$\eqalign {\int e^{f'-b}d\mu &= 1+\int (f'-b)d\mu +\int r(f'-b)d
\mu\cr
&\le \exp (\int (f'-b)d\mu +\int r(f'-b)d\mu )\cr}$$
and thus it suffices to show that $\int r(f'-b)d\mu\le 2t^2$.  We 
observe by (2.6.18) that
$$\int_{\{ f'>b\}}r(t^{-1}(f'-b))d\mu\le 1$$
and, since as already observed, the function $x^{-2}r(x)$ increases 
for $x>0$, this implies
$$\int_{\{ f'>b\}} r(f'-b)d\mu\le t^2\,.$$
Also, it is elementary to see that $r(x)\le x^2/2$ for $x<0$.  Now, 
by (2.6.15), we have
$$\eqalignno{\int_{\{ f'<b\}}(f'-b)^2d\mu &\le t^2\int w^2d\mu\le 
t^2\,. &\bx\cr}$$
\enddemo

{\bf 2.7.~~Penalties, IV.}

This section is devoted to remarkable fact that if (2.5.2) is suitably 
reinforced, the term $\exp t^2N$ can be removed in (2.5.3).

To express conveniently the conditions we need, we introduce the 
function $c(a,t)$, defined for $0<a<1$, $t>0$, as follows ($c$ stands 
for concentration):  if $\nu_1$ is the measure on $\Bbb R$ of 
density ${1 \over 2}e^{-\vert x\vert}$ with respect to Lebesgue 
measure, we have $c(a,t)=\nu_1((-\infty ,b+t])\,,$ where $b$ is given 
by $a=\nu_1((-\infty ,b])$.  Simple considerations show that
$$\eqalign{a\ge {1 \over 2} &\Rightarrow c(a,t)=1-e^{-t}(1-a)\cr
a\le {1 \over 2}\,,\quad e^ta\le {1 \over 2} &\Rightarrow c(a,t)=
e^ta\cr
a\le {1 \over 2}\,,\quad e^ta\ge {1 \over 2} &\Rightarrow c(a,t)=1-{1 
\over e^ta}\,.\cr}$$

\proclaim{Theorem 2.7.1}  Assume that for each subset $B$ of $\Omega$ 
we have
$$\leqalignno{t\le 1 &\Rightarrow\mu (\{ h(\cdot ,B)\le t^2\} )\ge c(
\mu (B),t) &(2.7.1)\cr
t\ge 1 &\Rightarrow \mu (\{ h(\cdot ,B)\le t\} )\ge c(\mu (B),t)\,. 
&(2.7.2)\cr}$$
Then, for each subset $A$ of $\Omega^N$, we have
$$\leqalignno{\int_{\Omega^N} e^{K^{-1}f_h(A,x)} dP(x) &\le {1 \over 
P(A)} &(2.7.3)\cr}$$
where $f_h$ is given by (2.4.3) and where $K$ is universal.
\endproclaim

Our first task should be to give natural examples of situations where 
(2.7.1), (2.7.2) occur.

\proclaim{Proposition 2.7.2}  Consider the probability $\nu_1$ on 
$\Bbb R$, of density ${1\over 2}e^{-\vert x\vert}$ with respect to 
Lebesgue measure.  Then the function $h(x,y)=\min (\vert x-y\vert ,
\vert x-y\vert^2)$ satisfies (2.7.1), (2.7.2) (for $\nu_1$, rather 
than $\mu$).
\endproclaim

\demo{Proof}  For a subset $C$ of $\Bbb R$, and $t>0$, let us set 
$C_t=\{ x\in\Bbb R\,;\,d(x,C)\le t\}$.  To prove (2.7.1), (2.7.2), it 
suffices to show that
$$\nu_1(C_t)\ge c(\nu_1(C),t)\,.$$
This is proved in [T4] using rearrangements.

We sketch below a simpler alternative argument to prove the weaker 
result
$$\leqalignno{\nu_1(C_t) &\ge c(\nu_1(C),t/2)\,. &(2.7.4)\cr}$$
(The reader should observe that this suffices to prove that $h/4$ 
satisfies (2.7.1), (2.7.2).) 

First, we reduce to the case where $C$ is a finite union of intervals.  
Setting
$$u(t)=\inf\{\vert x\vert\,;\,x\in\overline{C_t}/C_t\}\,,$$
it should be clear that
$$\leqalignno{{d\nu_1 \over dt}(C_t) &\ge {1 \over 2}\exp (-u(t))\,. 
&(2.7.5)\cr}$$

By definition of $u(t)$, we see that the interval $[-u(t), u(t)]$ is 
either contained in the closure of $C_t$, or else it does not meet 
$C_t$.  Thereby, we have either
$$\nu_1 (C_t)\ge 1-2\nu_1 ([u(t),\infty ))=1-e^{-u(t)}$$
or else
$$\nu_1(C_t)\le 2\nu_1([u(t),\infty ))=e^{-u(t)}$$
so that, in any case
$$e^{-u(t)}\ge \min (\nu_1(C_t),1-\nu_1(C_t))\,.$$
Combining with (2.7.5) shows that as long as $\nu_1(C_t)\le 1/2$, we 
have ${d \over dt}(\log\nu_1(C_t))\ge 1/2$, so that $\nu_1(C_t)\ge 
e^{t/2}\nu_1(C)$.  Similar considerations complete the proof.\rbx
\enddemo

Other examples can be generated using Proposition 2.7.2 and the 
following simple observation.

\proclaim{Proposition 2.7.3}  Consider a probability space $(\Omega ,
\mu )$, a function $h$ on $\Omega^2$, that satisfies (2.7.1), (2.7.2).  
Consider a measurable map $\eta$ from $\Omega$ to a measured space 
$\Omega '$, and the measure $\mu '=\eta (\mu )$ on $\Omega '$.  
Consider a function $h'$ on $\Omega^{\prime 2}$ such that
$$\leqalignno{\forall x,y\in\Omega\,,\,h'(\eta (x),\eta (y)) &\le 
h(x,y)\,. &(2.7.6)\cr}$$
Then $h',\mu '$ satisfy (2.7.1), (2.7.2).
\endproclaim

\demo{Proof}  This is obvious using the relations $\mu (\eta^{-1}(B))
=\mu '(B)$, $h(x,\eta^{-1}(B))\ge h'(\eta (x),B)$.

The use of Propositions 2.7.2 and 2.7.3 will allow the construction 
of a wide class of examples.
\enddemo

\proclaim{Proposition 2.7.4}  Consider a convex symmetric function 
$\psi\ge 0$ on $\Bbb R$, with $\lim\limits_{x\to\infty}\psi '(x)=
\infty$, and the probability $\nu_\psi$ of density $a_\psi e^{-\psi 
(x)}$ with respect to Lebesgue measure, where $a_\psi$ is the 
normalizing constant.  Then there is a constant $K(\psi )$ depending 
on $\psi$ only such that the function $h(x,y)$ on $\Bbb R^2$ given by
$$\leqalignno{\vert x-y\vert\le 1 &\Rightarrow h(x,y)={1 \over K(
\psi )}\vert x-y\vert^2 &(2.7.7)\cr
\vert x-y\vert\ge 1 &\Rightarrow h(x,y)={1 \over K(\psi )}\psi\left(
{1 \over K(\psi )}\vert x-y\vert\right) &(2.7.8)\cr}$$
satisfies (2.7.1), (2.7.2) with respect to $\nu_\psi$.
\endproclaim

\demo{Proof of Proposition 2.7.4}  Consider the nondecreasing map 
$\eta$ from $\Bbb R$ to $\Bbb R$ that transports $\nu_1$ to 
$\nu_\psi$.  Thus
$$\leqalignno{\int^\infty_{\eta (x)}a_\psi e^{-\psi (t)}d\lambda (t) 
&=\int^\infty_x {1 \over 2}e^{-\vert t\vert} d\lambda (t)\,. 
&(2.7.9)\cr}$$

By Propositions 2.7.2, 2.7.3, it suffices to show that
$$\leqalignno{h(\eta (x),\eta (y)) &\le \min (\vert x-y\vert ,\vert 
x-y\vert^2)\,. &(2.7.10)\cr}$$

It is simple to see that (2.7.10) will follow from (2.7.7), (2.7.8) 
(with a suitable choice of the constant there) provided we can show 
that
$$\leqalignno{\vert\eta (x)-\eta (y)\vert &\le K(\psi )\vert x-y\vert 
&(2.7.11)\cr
\psi\left({1 \over K(\psi )}\vert\eta (x)-\eta (y)\vert\right) &\le 
\vert x-y\vert\,. &(2.7.12)\cr}$$
There, as in the rest of this proof, $K(\psi )$ denotes a constant 
depending on $\psi$ only, that may vary at each occurrence.

To prove (2.7.11), it suffices to prove that $\eta '(x)$ is bounded 
when $x>0$.  Differentiating (2.7.9), we get $a_\psi\eta '(x)e^{-\psi 
(\eta (x))} = e^{-x}/2$, and plugging back in (2.7.9), we get
$$\eta '(x)=\int^\infty_{\eta (x)}e^{\psi (\eta (x))-\psi (t)}dt\,.$$
Thereby, it suffices to show that
$$\sup_{u\ge 0}\int^\infty_{u} e^{-\psi (t)+\psi (u)}dt<\infty\,.$$
Given $u_0>0$, the supremum for $u\le u_0$ is certainly bounded.  On 
the other hand, for $u\ge u_0$, 
by convexity of $\psi$ we have $\psi (t)-\psi (u)\ge (t-u)\psi '(u_0)$, 
so it suffices to choose $u_0$ with $\psi '(u_0)>0$.

We now turn to the proof of (2.7.12).  It suffices to prove that, for 
$y>x\ge 0$, we have $\psi ((\eta (y)-\eta (x))/K(\psi ))\le y-x$.  
Setting $a=\psi^{-1}(y-x)$, it suffices to show that $\eta 
(x+\psi (a))\le\eta (x)+K(\psi )a$, i.e., that
$$\leqalignno{a_\psi\int^\infty_{\eta (x)+K(\psi )a}e^{-\psi (t)}dt 
&\le {1 \over 2}e^{-x-\psi (a)}\,. &(2.7.13)\cr}$$

First, we note that, since $\psi (t)\ge \psi (y)+(t-y)\psi '(y)$, we 
have, for $y>0$, that
$$\int^\infty_y e^{-\psi (t)}dt\le {1 \over \psi '(y)}e^{-\psi (y)}$$
so that
$$\leqalignno{a_\psi \int_{\eta (x)+2a}e^{-\psi (t)}dt &\le {a_\psi 
\over \psi '(\eta (x)+2a)}e^{-\psi (\eta (x)+2a)}\,. &(2.7.14)\cr}$$

Also,
$$\leqalignno{{1 \over 2} e^{-x} &= a_\psi\int^\infty_{\eta (x)}
e^{-\psi (t)}dt \ge a_\psi ae^{-\psi (\eta(x)+a)}\,. &(2.7.15)\cr}$$

Since $\psi '(y)$ increases for $y>0$, we have $\psi (\eta (x)+2a)\ge
\psi (\eta (x)+a)+\psi (a)$.  
Thus, from (2.7.14), (2.7.15) we see that (2.7.13) holds provided 
$K(\psi )\ge 2$, $a\psi '(\eta (x)+2a)\ge a_\psi$.

On the other hand, using again convexity, we see that
$$\eqalign{\int^\infty_{\eta (x)+a}e^{-\psi (t)}dt=\int^\infty_{\eta 
(x)}e^{-\psi (v+a)}dv &\le e^{-a\psi '(\eta (x))}\int^\infty_{\eta 
(x)}e^{-\psi (v)}dv\cr
&\le {1 \over 2a_\psi} e^{-a\psi '(\eta (x))-x}\,. \cr}$$

Thereby, if $K(\psi )\ge 1$, (2.7.13) will hold provided $a\psi '(\eta 
(x))\ge \psi (a)$, and in particular if $\eta (x)\ge a$.

Thus we can assume $\eta (x)\le a$, $a\psi '(\eta (x)+2a)\le a_\psi$.  
This means that $a$ and $x$ stay bounded; but the conclusion is 
obvious then.\rbx
\enddemo

It is of particular interest to consider the case where $\psi (x)=
x^2$, so that $\nu_\psi$ is Gaussian.  In this case, Proposition 2.7.4 
shows that one can take $h(x,y)=K^{-1}(x-y)^2$.  This recovers the 
concentration of measure for the Gauss space, as expressed by (1.7).  
There is, however, a big loss of information in (2.7.10); and the 
result obtained by taking
$$h(x,y)=\min (\vert\eta^{-1}(x)-\eta^{-1}(y)\vert ,(\eta^{-1}(x)-
\eta^{-1}(y))^2)\,.$$
is rather more precise than (1.7).

The induction step of the proof of Theorem 2.7.1 reduces to the 
following.

\proclaim{Proposition 2.7.5}  There exists a universal constant $L$ 
with the following property.  Consider a function $g$ on $\Omega$, and 
define
$$\leqalignno{\widehat{g}(x) &= \inf_{y\in\Omega}g(y)+{1 \over L}h
(x,y)\,. &(2.7.16)\cr}$$

Then, under (2.7.1), (2.7.2), we have
$$\leqalignno{\int_\Omega e^{\widehat g} d\mu\int_\Omega e^{-g}d\mu 
&\le 1\,. &(2.7.17)\cr}$$
\endproclaim

Let us recall that we denote by $\nu_1$ the probability measure on 
$\Bbb R$ of density $e^{-\vert x\vert}/2$ with respect to Lebesgue 
measure.  During the end of this section, for $x\in\Bbb R$ we 
set $\varphi (x)=\min (\vert x\vert ,x^2)$.

The proof of Proposition 2.7.5 is considerably simplified by the 
following observation.

\proclaim{Proposition 2.7.6}  Consider a function $g$ on $\Omega$, 
and $\widehat g$ given by (2.7.16).  Then we can find two 
nonincreasing functions $g_1$, $g_2$ on $\Bbb R$ with the
following properties
$$\leqalignno{\int_\Omega e^{\widehat g}d\mu &=\int_{\Bbb R}e^{g_2}
d\nu_1\,;~~\int_\Omega e^{-g} d\mu = \int_{\Bbb R}e^{-g_1}d\nu_1 
&(2.7.18)\cr}$$
$$\leqalignno{\forall x\in\Bbb R\,,~~g_2(x) &\le \inf_{y\in\Bbb R}
g_1(y)+{1 \over L}\varphi (\vert x-y\vert )\,. &(2.7.19)\cr}$$
\endproclaim

In particular, this implies that we have reduced to the case $\mu = 
\nu_1$, $g$ nonincreasing, $h(x,y)=\varphi (\vert x-y\vert )$.

\demo{Proof}  We define, for $y\in\Bbb R$
$$\leqalignno{g_1(y) &= \inf\{ t\,;\,\mu (\{ g\le t\})\ge \nu_1([y,
\infty ))\} &(2.7.20)\cr
g_2(x) &= \sup\{ u\,;\,\mu (\{\widehat {g}\ge u\})\ge \nu_1((-\infty ,
x])\}\,. &(2.7.21)\cr}$$
Thereby both $g_1$, $g_2$ are nonincreasing; it should be obvious that 
(2.7.18) holds.  We prove (2.7.19).  Consider $x\le y$.  By (2.7.20), 
we have $\mu (B)\ge \nu_1 ([y,\infty ))$, where $B=\{ g\le g_1(y)\}$.  
By (2.7.1), (2.7.2), we have
$$\eqalign{\mu (\{ h(\cdot ,B)\le\varphi (t)\} ) &\ge c(\mu (B),t)\cr
&\ge c(\nu_1 ([y,\infty )),t)\cr
&=\nu_1 ([y-t,\infty ))\,.\cr}$$
Since $\widehat{g}\le g_1(y)+\varphi (t)/L$ on the set $\{h(\cdot ,B)
\le \varphi (t)\}$, we get
$$\mu\left(\left\{\widehat{g}\le g_1(y)+{\varphi (t) \over L}\right\}
\right) \ge \nu_1([y-t,\infty ))\,.$$

On the other hand, by (2.7.21) we have
$$\mu (\{\widehat{g}\ge g_2(x)\} )\ge\nu_1((-\infty ,x])\,.$$
Thus, if $t>y-x$, we have $g_2(x)<g_1(y)+\varphi (t)/L$.  Thus $g(x)
\le g(y)+\varphi (y-x)/L$, and (2.7.19) follows.\rbx
\enddemo

We next show that we have reduced the proof of Proposition 2.7.5 to 
the following.

\proclaim{Proposition 2.7.7}  There exists a universal constant $L$ 
with the following property.  Consider a nonincreasing function $f$ on 
$\Bbb R$, with $f(0)=0$.  Define
$$\leqalignno{\widehat{f}(x) &= \inf_{y\in\Bbb R}f(y)+{1 \over L}
\varphi (\vert x-y\vert )\,.&(2.7.22)\cr}$$
Then, if $f$ has a Lipschitz constant $\le 2/L$ we have
$$\int e^{\widehat f}d\nu_1\int e^{-f}d\nu_1\le 1\,.$$
\endproclaim

We prove the claim stated before Proposition 2.7.7.  In view of 
Proposition 2.7.6 and (2.7.19), it suffices to prove that $\int 
e^{\widehat{g}_1}d\nu_1\int e^{-g_1}d\nu_1\le 1$, where $\widehat{g}_1$ 
is given by the right-hand side of (2.7.19).  Define now
$$\leqalignno{f(y) &= \sup_{x\in\Bbb R}\widehat{g}_1(x)-{1 \over L}
\varphi (\vert x-y\vert )\,. &(2.7.23)\cr}$$
Since for all $x$ and $y$ we have $\widehat{g}_1(x)\le g_1(y)+{1 \over 
L}\varphi (\vert x-y\vert )$, we see that $f(y)\le g_1(y)$.  Thus, 
$\int e^{-f}d\nu_1\ge \int e^{-g_1}d\nu_1$.  Also, by (2.7.23), we 
have $\widehat{g}_1(x)\le f(y)+{1 \over L}\varphi (\vert x-y\vert )$ 
for all $x$, $y$, so that $\widehat{g}_1\le\widehat{f}$.  Thereby it 
suffices to prove that $\int e^{\widehat f} d\nu_1\int e^{-f}d\nu_1
\le 1$.  The condition $f(0)=0$ is certainly not restrictive, and $f$ 
has a lipschitz constant $\le 2/L$ by (2.7.23) since $\varphi$ has a 
lipschitz constant $\le 2$.

Upon seeing the result of [T4] exposed in a seminar, B. Maurey 
produced a rather magic proof of Proposition 2.7.7.  The proof we will 
give is more in the spirit of the arguments of the present paper, 
and is likely to be more instructive as it prepares for the 
considerably more delicate results to be presented in Chapter 4.  We 
start by a simple lemma.

\proclaim{Lemma 2.7.8}  Consider a nonincreasing function $u$ on 
$\Bbb R$, such that $u(0)=0$.  Then
$$\int_{\Bbb R^-}u^2d\nu_1\le K\sum_{k\ge 1}(u(-k)-u(-k+1))^2e^{-k}
\,.$$
\endproclaim

\demo{Proof}  For simplicity we set $u_k=u(-k)$.  Thus
$$\int_{\Bbb R^-}u^2d\nu_1\le S=:\sum_{k\ge 1}u^2_k e^{-k+1}\,.$$
Since $u^2_k\le 2u^2_{k-1}+2(u_k-u_{k-1})^2$, we have
$$S\le 2\sum_{k\ge 1} u^2_{k-1}e^{-k+1}+2\sum_{k\ge 1}(u_k-u_{k-1})^2
e^{-k+1}\,.$$
But since $u_0=0$, the first sum is exactly $2S/e$, so that
$$\eqalignno{S\left(1-{2 \over e}\right) &\le 2e\sum_{k\ge 1} (u_k-
u_{k-1})^2e^{-k}\,. &\bx\cr}$$
\enddemo

During the proof of Proposition 2.7.7, we will consider another number 
$1\le M\le L$.  The numbers $M$, $L$ will be chosen later.  The 
crucial part of the proof of Proposition 2.7.7 is as follows.

\proclaim{Proposition 2.7.9}  Consider a non-increasing function $u$ 
on $\Bbb R$, with $u(0)=0$.  Assume that $\vert u\vert \le 1/M$, and 
set $\widehat{u}(x)=\inf\limits_y u(y)+\varphi(\vert x-y\vert )/L$.  
Then, if $L\ge KM$, we have 
$$\leqalignno{\int_{\Bbb R} (u-\widehat{u})d\nu_1 &\ge {M \over K}
\int_{\Bbb R} u^2d\nu_1\,. &(2.7.23)\cr}$$
Moreover, if $M\ge K$, we have
$$\leqalignno{\int_{\Bbb R}e^{\widehat u} d\nu_1+\int_{\Bbb R}
e^{-u}d\nu_1 &\le 2-{M\over K}\int_{\Bbb R}u^2d\nu_1\,. &(2.7.24)\cr}$$
\endproclaim

\demo{Proof}  To prove (2.7.23), it suffices to prove it when the 
right-hand side is replaced by $\int_{\{u\ge 0\}}u^2d\mu$ (resp. 
$\int_{\{ u\le 0\}}u^2d\mu$).  The arguments for these two cases are 
similar so we treat the first case only.  We set $u_k=u(-k)$, so that 
$u_k\le M^{-1}$, and $M(u_k-u_{k-1})\le 1$.  We set $N_k=[2/M(u_k-
u_{k-1})]$.  Thus we have $N_k\ge 2$ and
$$\leqalignno{{1 \over 2MN_k} &\le u_k-u_{k-1}\le {1 \over MN_k}\,. 
&(2.7.25)\cr}$$
For $k\ge 1$, $\ell\ge 0$ we set $a_{k,\ell}=-k+1-\ell /N_k$, and 
$u_{k,\ell}=u(a_{k,\ell})$.  Thus $u_{k,0}=u_{k-1}$, $u_{k,N_k}=u_k$.  
For $k\ge 1$, $1\le\ell\le N_k$, we consider the subset $R_{k,
\ell}$ of $\Bbb R^2$ given by
$$R_{k,\ell}=]a_{k,\ell+1},a_{k,\ell}[\times ]\min (u_{k,\ell},u_{k,
\ell -1}+{4 \over LN_k^2}),u_{k,\ell}[\,.$$
We observe that no point belongs to more than two intervals $]a_{k,
\ell +1},a_{k,\ell}[$, for $1\le\ell\le N_k$, $k\ge 1$ so that the 
rectangles $R_{k,\ell}$ have the same property.  Since $u(x)\ge u_{k,
\ell}$ for $x\le a_{k,\ell}$, $R_{k,\ell}$ is below the graph of $u$; 
but, since $u(a_{k,\ell -1})=u_{k,\ell -1}$, we have $\widehat{u}(x)
\le u_{k,\ell -1}+4/LN^2_k$ on $[a_{k,\ell +1},a_{k,\ell}]$.  Thus 
$R_{k,\ell}$ is above the graph of $\widehat{u}$, and hence
$$\int_{\Bbb R} (u-\widehat{u})d\nu_1\ge {1 \over 2}\sum_{k\ge 1}~
\sum_{1\le\ell\le N_k}\nu_1([a_{k,\ell +1},a_{k,\ell}])\left( u_{k,
\ell}-\min\left( u_{k,\ell},u_{k,\ell -1}+{4 \over LN^2_k}\right)
\right)\,.$$
Since $\nu_1([a_{k,\ell +1},a_{k,\ell}])\ge e^{-k}/KN_k$, we have
$$\eqalign{\int_{\Bbb R} (u-\widehat{u})d\nu_1 &\ge {1 \over K}
\sum_{k\ge 1} {e^{-k}\over N_k}\sum_{1\le\ell\le N_k}\left( u_{k,\ell}
-u_{k,\ell -1}-{4 \over LN^2_k}\right)\cr
&= {1 \over K}\sum_{k\ge 1}{e^{-k} \over N_k} \left( u_{k,N_k}-u_{k,0} 
-{4 \over LN_k}\right)\cr
&={1 \over K}\sum_{k\ge 1}{e^{-k} \over N_k}\left( u_k-u_{k-1} -{4 
\over LN_k}\right)\cr
&\ge {1 \over K}\sum_{k\ge 1} e^{-k}(u_k-u_{k-1})^2\cr}$$
by (2.7.25), and provided $L\ge 16M$.  Thus, (2.7.23) follows from 
Lemma 2.7.8.

To prove (2.7.24), we use that $e^x\le 1+x+x^2$ for $\vert x\vert\le 
1$.  Thus
$$\int_{\Bbb R} e^{\widehat{u}} d\nu_1 +\int_{\Bbb R}e^{-u}d\nu_1\le 
2+\int_{\Bbb R}(\widehat{u} -u)d\nu_1 +\int_{\Bbb R} u^2d\nu_1+
\int_{\Bbb R}\widehat{u}^2d\nu_1\,.$$
Now,
$$\widehat{u}^2\le 2u^2+2(u-\widehat{u})^2\le 2u^2+{2 \over M} (u-
\widehat{u})\le 2u^2+{1 \over 2}(u-\widehat{u})$$
provided $M\ge 4$, and thus
$$\int_{\Bbb R}e^{\widehat u}d\nu_1+\int_{\Bbb R} e^{-u}d\nu_1\le 2+
3\int_{\Bbb R}u^2d\nu_1-{1 \over 2}\int_{\Bbb R}(u-\widehat{u})d\nu_1$$
and the result follows from (2.7.23).\rbx
\enddemo

\demo{Proof of Proposition 2.7.7}  We observe that, for $a\in\Bbb R$, 
we have $a(2-a)\le 1$.  Thus it suffices to show that
$$\leqalignno{\int_{\Bbb R}e^{\widehat f}d\nu_1+\int_{\Bbb R}e^{-f}d
\nu_1 &\le 2\,. &(2.7.26)\cr}$$
We set $u=\min (1/M,\max (f,-1/M))$.  Thus
$$\leqalignno{\int_{\Bbb R} e^{-f}d\nu_1 &\le \int_{\Bbb R} e^{-u}d
\nu_1+\int^\infty_b (e^{-f}-e^{1/M})d\nu_1 &(2.7.27)\cr}$$
where $f(b)=-1/M$.  We observe that if $\widehat{u}(x)<1/M$, then 
$\widehat{f}(x)\le\widehat{u}(x)$.  Indeed if $\widehat{u}(x)<1/M$, 
then given $\varepsilon$ with $\widehat{u}(x)<\varepsilon <1/M$, 
there exists $y$ with $u(y)+L^{-1}\varphi (\vert x-y\vert )<
\varepsilon$.  Thus $u(y)<1/M$, and thus $f(y)\le u(y)$, so that 
$\widehat{f}(x)<\varepsilon$.  Then, if $c$ is the largest so that 
$\widehat{f}(c)=1/M$, we have
$$\leqalignno{\int_{\Bbb R} e^{\widehat f}d\nu_1 &\le \int_{\Bbb R}
e^{\widehat u}d\nu_1+\int^c_{-\infty}(e^{\widehat f}-e^{1/M})d\nu_1
\,. &(2.7.28)\cr}$$
Since $f(0)=0$, we have $c<0<b$.  Since $f$ has a Lipschitz constant 
$\le 2/L$, for $x\ge b$ we have
$$-f(x)\le {1 \over M} +{2 \over L} (x-b)$$
and thus
$$\eqalign{\int^\infty_b e^{-f} d\nu_1 &\le \int^\infty_b{1 \over 2}
e^{1/M+2(x-b)/L}e^{-x}dx\cr
&= {e^{1/M} \over 1-2/L}\nu_1 ([b,\infty ))\,.\cr}$$
Hence we have
$$\eqalign{\int^\infty_b (e^{-f}-e^{1/M})d\nu_1 &\le e^{1/M}\left(
{1 \over 1-2/L} - 1\right)\nu_1 ([b,\infty ])\cr
&\le {K \over L}\nu_1 ([b,\infty [)\cr
&\le {KM^2 \over L}\int_{\Bbb R} u^2d\nu_1\cr}$$
since $u(x)= -1/M$ for $x\ge b$.  
Using (2.7.27), (2.7.30), and making a similar computation for 
$\int^c_{-\infty} (e^{\widehat f} - e^{1/M})d\nu_1$ yields
$$\int_{\Bbb R}e^{\widehat f}d\nu_1 +\int_{\Bbb R} e^{-f}d\nu_1\le
\int_{\Bbb R} e^{\widehat u} d\nu_1+\int_{\Bbb R} e^{-u}d\nu_1+{KM^2 
\over L}\int_{\Bbb R} u^2d\nu_1\,.$$
It then follows from (2.7.24) that (2.7.26) holds provided $M\ge K$, 
$L\ge KM^2$.\rbx
\enddemo

It would be of interest to understand exactly which are the functions 
$\varphi$ such that, if one sets
$$\widehat{f}(x) = \inf_{y\in\Bbb R} f(y)+\varphi (x-y)\,,$$
then $\int e^{\widehat f} d\nu_1\int e^{-f} d\nu_1\le 1$.  On the 
other hand, the situation is considerably clearer if one considers 
the standard Gaussian density $\gamma_1$ rather than $\nu_1$.  
In that case, the obvious adaptation of Maurey's argument shows that 
if $\alpha\ge 1$, and if $\widehat{f}(x) = \inf\limits_{y\in\Bbb R}
\alpha f(y)+{\alpha \over 2(\alpha +1)} (x-y)^2$, then 
$\int e^{\widehat f} d\gamma_1(\int e^{-f}d\gamma_1)^\alpha\le 1$.  
Thereby, by induction, and with the notations of (1.7), we get
$$\gamma_N(A_t)\ge 1-{1 \over \gamma_N(A)^\alpha} e^{-{\alpha t^2 
\over 2(\alpha +1)}}$$
so that, by optimization over $\alpha$, for $t\ge\sqrt{2\log (1/
\gamma_N(A))}$, we get
$$\gamma_N(A_t)\ge 1-\exp -{1 \over 2} (t-\sqrt{2\log (1/\gamma_N(A)
)}\,)^2$$
which is not so far from (1.7).
\vfil\eject

\noindent{\bf 3.~~Control by $q$ points}

{\bf 3.1.~~Basic result}

Consider an integer $q\ge 2$.  For subsets $A_1,\dots ,A_q$ of 
$\Omega^N$, and $x\in\Omega^N$, we set
$$\leqalignno{f(A_1,\dots A_q,x) = &\inf\{k\,;\,\exists y^1\in A_1
\,,\,y^2\in A_2,\dots ,y^q\in A_q\,; &(3.1.1)\cr
&\card\{ i\le N\,;\,x_i\not\in\{ y^1_i,\dots ,y^q_i\}\}\le k\}\,.\cr}$$

\proclaim{Theorem 3.1.1}  We have
$$\leqalignno{\int q^{f(A_1,\dots ,A_q,x)}dP(x) &\le {1 \over 
\prod\limits_{i\le q}P (A_i)}\,. &(3.1.2)\cr}$$
In particular we have
$$\leqalignno{P(\{ f(A,\dots A,x)\ge k\}) &\le {1 \over q^kP(A)^q}
\,. &(3.1.3)\cr}$$
\endproclaim

The induction method will reduce this statement to a simple fact about 
functions.

\proclaim{Lemma 3.1.2}  Consider a function $g$ on $\Omega$, such that 
$1/q\le g\le 1$.  Then
$$\leqalignno{\int_\Omega {1 \over g}d\mu (\int_\Omega gd\mu )^q 
&\le 1\,. &(3.1.4)\cr}$$
\endproclaim

\demo{Proof}  We could use the extreme point argument of Lemma 2.1.2.  
One alternative method is as follows.  Observing that $\log x\le x-1$, 
to prove that $ab^q\le 1$ it suffices to show that $a+qb\le 
q+1$.  Thus, it suffices to show that
$$\int_\Omega {1 \over g} d\mu +q\int_\Omega gd\mu\le q+1\,.$$
But this is obvious since $x^{-1}+qx\le q+1$ for $q^{-1}\le x\le 
1$.\rbx
\enddemo

\proclaim{Corollary 3.1.3}  Consider functions $g_i$ on $\Omega$, 
$g_i\le 1$.  Then
$$\leqalignno{\int_\Omega \min_{i\le q} \left( q,{1 \over g_i}\right) 
d\mu\prod_{i\le q}\int g_id\mu &\le 1\,. &(3.1.5)\cr}$$
\endproclaim

\demo{Proof}  Set $g=(\min\limits_{i\le q}(q,g^{-1}_i))^{-1}$, 
observe that $g_i\le g$, and use (3.1.4).

We now prove Theorem 3.1.1 by induction over $N$.  For $N=1$, the 
result follows from (3.1.5), taking $g_i=1_{A_i}$.

We assume now that Theorem 3.1.1 has been proved for $N$, and we 
prove it for $N+1$.  Consider sets $A_1,\dots ,A_q$ of $\Omega^{N+1}$.  
For $\omega\in\Omega$, we define the sets $A_i(\omega )$ 
as in (2.1.5) and we consider the projection $B_i$ of $A_i$ on 
$\Omega^N$.  The basic observation is that
$$\leqalignno{f(A_1,\dots ,A_q,(x,\omega )) &\le 1+f(B_1,\dots B_q,x) 
&(3.1.6)\cr}$$
and that, if, $j\le q$
$$f(A_1,\dots , A_q,(x,\omega )) \le f(C_1,\dots ,C_q,x)$$
where $C_i=B_i$ for $i\not= j$, $C_j=A_j(\omega )$.

If we set $g_i(\omega )=P(A_i(\omega ))/P(B_i)$ using Fubini theorem 
and induction hypothesis, we are reduced to show that
$$\int_\Omega\min \left( q,\min\limits_{i\le q} {1 \over g_i(\omega )}
\right)\le {1 \over \prod\limits_{i\le q}\int_\Omega g_id\mu}$$
which is (3.1.5).\rbx
\enddemo

{\bf 3.2.~~Sharpening.}

Given $\alpha >1$, we can now, in the spirit of Proposition 2.2.1, 
look for the largest number $a=a(q,\alpha )$ for which we can prove 
that
$$\leqalignno{\int^\ast a(q,\alpha )^{f(A_1,\dots ,A_q,x)}dP(x) 
&\le {1 \over \prod\limits_{i\le q}P (A_i)^\alpha}\,. &(3.2.1)\cr}$$
Following the proof of Theorem 3.1.1, we see that we can take for 
$a(q,\alpha )$ the unique number $x>1$ such that
$$\leqalignno{x+q\alpha x^{-1/\alpha} &= 1+q\alpha\,. &(3.2.2)\cr}$$

It then follows from (3.2.1) that
$$\leqalignno{P(\{ f(A,\dots ,A,x)\ge k\}) &\le \inf_{\alpha\ge 1}
{a(q,\alpha )^{-k} \over P(A)^{q\alpha}}\,. &(3.2.3)\cr}$$

There is no obvious way to compute the right-hand side of (3.2.3).  
However, for large $q$, we have the following, that improves upon 
(3.1.2) for large values of $k$ ($k\gg q\log q$).

\proclaim{Proposition 3.2.1}  There exists a universal constant $q_0$ 
such that, if $g\ge q_0$, we have
$$\leqalignno{P(\{ f(A,\dots ,A,x)\ge k\} ) &\le \left({e \over (e-1)
q\log q}\right)^k\left({1 \over P(A)}\right)^{q\log q}\,. &(3.2.4)\cr}$$
\endproclaim

\demo{Proof}  We take $\alpha =\log q$, and we show that for $q$ 
large enough, we have
$$a(q,\alpha )\ge a:= 1+\left( 1-{1 \over e}\right) q\log q\,.$$
For large $q$, we have $a\ge q$, so that $a^{1/\alpha}\ge e$, so that
$$a-1\le q\alpha\left( 1-{1 \over a^{1/\alpha}}\right)$$
and thus $a+q\alpha a^{-1/\alpha}\le 1+q\alpha$.\rbx
\enddemo

It is interesting to note that Proposition 3.2.1 is rather sharp.  
Consider the case where $\Omega =\{ 0,1\}$, and where $\mu$ gives 
weight $p$ to $1$ ($p\le 1/2$).  Assume for simplicity that $r=pN$ 
is an integer.  Consider the set $A=\{ x\in\Omega^N\,;\,
\sum\limits_{i\le N}x_i\le r\}$.  Then $P(A)$ is of order $1/2$.  
Considering $s=rq+k$, we clearly have that $\sum\limits_{i\le N} x_i
=s$ implies $f(A,\dots ,A,x)\ge k$.  Thus $P(\{ f(A,\dots ,A,x)\ge 
k\})\ge p^s(1-p)^{N-s}{N \choose s}$.

When $s\le N/2$, we have ${N\choose s}\ge (N/2s)^s$, so that
$$\eqalign{P(\{ f(A,\dots ,A,x\}\ge k\}) &\ge \left({pN \over 2s}
\right)^s e^{-2pN}\ge \left({r \over 2es}\right)^s\cr
&\ge \left({1 \over 2e(q+k/r)}\right)^{rq+k}\,.\cr}$$
If we take $k\ge q\log q$, fixed, and then $r$ of order $k/q\log q$, 
we get a lower bound of order $(1/K q\log q)^k$.

{\bf 3.3.~~Penalties.}

The result of this section is the one single major theorem of Part I 
that has not been motivated by direct applications.  Rather, it has 
been motivated by a desire of symmetry with Sections 2.7 and 4.4.

We consider a ``penalty function'' $h(\omega ,\omega^1,\dots ,
\omega^q)$ on $\Omega^{q+1}$.  We assume $h\ge 0$ and
$$\leqalignno{\omega\in\{\omega^1,\dots ,\omega^q\} &\Rightarrow h(
\omega ,\omega^1,\dots ,\omega^q) =0\,. &(3.3.1)\cr}$$

For subsets $A_1,\dots ,A_q$ of $\Omega^N$, we consider
$$\leqalignno{f_h(A_1,\dots ,A_q,x) &=\inf\left\{\sum_{i\le N}h(x_i,
y^1_i,\dots ,y^q_i)\,; y^1\in A_1,\dots ,y^q\in A_q\right\}\,. 
&(3.3.2)\cr}$$
The case considered in Section 3.1 is where $h(\omega ,\omega^1,\dots 
,\omega^q)=1$ unless $\omega\in\{\omega^1,\dots ,\omega^q\}$, in which 
case it is zero.

Given subsets $B_1,\dots ,B_q$ of $\Omega$, we set
$$\leqalignno{h(\omega ,B_1,\dots ,B_q)=\inf\{ h(\omega ,\omega^1,
\dots ,\omega^q)\,;\,\omega^1\in B_1,\dots ,\omega^q\in B_q\}\,. 
&&(3.3.3)\cr}$$

To control how large $h$ is, we will consider a nonincreasing function 
$\gamma$ from $]0,1]$ to $\Bbb R^+$, and assume that
$$\leqalignno{\forall\omega\in\Omega\,,\,\forall B_1,\dots ,B_q
\subset\Omega\,,\,h(\omega ,B_1,\dots ,B_q) &\le \sum_{i\le q}\gamma 
(\mu (B_i))\,. &(3.3.4)\cr}$$

A typical case where this condition is satisfied is when
$$h(\omega ,\omega^1,\dots ,\omega^q)=\sum_{i\le q} h_i(\omega^i)$$
for functions $h_i$ that satisfy the tail condition $\mu (\{ h_i\ge
\gamma (t)\})\le t$ and when $\gamma$ is left continuous.  Indeed, if 
$t<\mu (B_i)$, then $B_i$ contains a point $y_i$ with $h_i(y_i)<
\gamma (t)$.

We consider a convex function $\theta\colon ]0,1]\to\Bbb R^+$, and we 
make the mild technical assumption that the inverse function $\xi$ 
satisfies
$$\leqalignno{\vert\xi '(x+1)\vert &\ge {1 \over 3}\vert\xi '(x)\vert
\,. &(3.3.5)\cr}$$
(We put ${1 \over 3}$ rather than ${1 \over 2}$ simply to allow the 
case $\xi (x)=e^{-x}$.)

\proclaim{Theorem 3.3.1}  There exists a universal constant $K$ such 
that for $q\ge K$, under (3.3.1), (3.3.4), (3.3.5), if, for each $s
\le 1$, we have
$$\leqalignno{\int^1_s\gamma^{-1}(\theta(s)-\theta (t))d\lambda (t) 
&\le {\log (q/K) \over q\vert\theta '(s)\vert} &(3.3.6)\cr}$$
then, for each subsets $A_1,\dots ,A_q$ of $\Omega^N$, we have
$$\leqalignno{\int e^{f_h(A_1,\dots ,A_q,x)}dP(x) &\le \exp\left(
\sum_{i\le q}\theta (P(A_i))\right)\,.&(3.3.7)\cr}$$
\endproclaim

To understand (3.3.6) better, we observe that the term $\theta '(s)$ 
arises simply because $\theta (s)-\theta (t)$ resembles $(s-t)\theta 
'(s)$ for $t$ close to $s$.  Actually, since $\theta (s)-\theta (t)
\le (s-t)\theta '(s)$, change of variable and Lebesgue theorem show 
that (3.3.6) implies that $\int^\infty_0\gamma^{-1}(u)du\le q^{-1}
\log (q/K)$.  In the case where $\gamma$ is constant, one can take 
$h(\omega ,\omega^1,\dots ,\omega^q)=q\gamma$ whenever $\omega\not\in
\{\omega^1,\dots ,\omega^q\}$ (and otherwise $h=0$).  Then the 
integral in (3.3.6) has to be interpreted as $\vert\{ t\colon s\le 
t\,;\,\theta (t)\ge\theta (s)-\gamma\}\vert$.  When $\theta (x)=-\log 
x$, this is $s(e^\gamma -1)$, and (3.3.6) holds whenever $\gamma\le 
q^{-1}\log (q/K)$.  We then almost recover Theorem 3.1.1.

To prove Theorem 3.3.1, it suffices, by the induction method, to prove 
the following.

\proclaim{Proposition 3.3.2}  There exists a universal $K$ such that, 
under conditions (3.3.1), (3.3.4), (3.3.5), (3.3.6) if, we consider 
functions $(u_i)_{i\le q}$ on $\Omega$, $0\le u_i\le 1$, and define
$$\leqalignno{v(\omega ) &=\inf_{\omega^1,\dots ,\omega^q}\sum_{i\le 
q}\theta (u_i(\omega^i))+h(\omega ,\omega^1,\dots ,\omega^q) 
&(3.3.8)\cr}$$
then we have
$$\leqalignno{\int_\Omega e^vd\mu &\le \exp \left(\sum_{i\le q}\theta
\left(\int_\Omega u_id\mu \right)\right)\,. &(3.3.9)\cr}$$
\endproclaim

\demo{Proof}  For clarity, we will replace (3.3.6) by
$$\leqalignno{\int^1_s \gamma^{-1}(\theta (s)-\theta (t))d\lambda (t) 
&\le {\tau \over \vert\theta '(s)\vert} &(3.3.10)\cr}$$
and we will determine in due time a good choice for $\tau$.  We 
already assume $\tau\le 1$.  The two main parts of the proof are the 
research of upper bounds for $\int_\Omega e^vd\mu$, and of lower 
bounds for $\sum\limits_{i\le q}\theta (\int_\Omega u_id\mu )$.

{\bf Step 1}  For $i\le q$, we set $S_i=\inf\limits_{\omega}\theta 
(u_i(\omega ))=\theta (\sup\limits_\omega u_i(\omega ))$.  By (3.3.8) 
and (3.3.1), taking $\omega^i=\omega$, we see that if we set $S=
\sum\limits_{i\le q}S_i$, we have
$$\leqalignno{v(\omega ) \le \theta (u_i(\omega )) +\sum_{j\not= i}
S_j &= \theta (u_i(\omega )) +S- S_i\,. &(3.3.11)\cr}$$

{\bf Step 2}  We make the convention that $\gamma (0)=\infty$.  For 
$i\le q$, we define $s_i$ by
$$\theta (s_i)=\inf_{t\ge 0}\{\gamma (\mu (\{ u_i\ge t\} ))+\theta 
(t)\}\,.$$
Thus we have $\theta (s_i)\ge S_i$ and for $t>s_i$  we have
$$\leqalignno{\mu (\{ u_i\ge t\} ) &\le \gamma^{-1}(\theta (s_i)-
\theta (t))\,. &(3.3.12)\cr}$$

{\bf Step 3}  We show that for any subset $C$ of $\Omega$, we have
$$\leqalignno{\int_C e^vd\mu &\le \mu (C)\exp\left(\sum_{i\le q}\theta 
(s_i)\right)\,. &(3.3.13)\cr}$$
By definition of $s_i$, given $\varepsilon >0$, we can find $t_i$ such 
that $\gamma (\mu (B_i))+\theta (t_i)\le \theta (s_i)+\varepsilon$, 
where $B_i=\{ u_i\ge t_i\}$.  Since $\theta (u_i(\omega^i )) \le\theta 
(t_i)$ for $\omega^i\in B_i$, we have, by (3.3.8)
$$v(\omega )\le\sum_{i\le q} \theta (t_i)+h(\omega ,B_1,\dots ,B_q)$$
so that (3.3.13) follows by (3.3.4), since $\varepsilon$ is arbitrary.

{\bf Step 4}  Consider now a number $m$.  We set
$$\leqalignno{z &=\int \min ((v-m)^+,1)d\mu~~\text{and}~~C=\{ v\ge m+1
\}\,. &(3.3.14)\cr}$$
Thus, in particular $\mu (C)\le z$.

Using the inequality $e^x\le 1+2x^+$ for $x\le 1$, we get, using 
(3.3.13)
$$\int_\Omega e^{v-m}d\mu\le\int_{\Omega\backslash C}+\int_C \le 1+
2z+\mu (C)\exp\left(\sum_{i\le q}\theta (s_i)-m\right)$$
so that
$$\leqalignno{\int_\Omega e^{v-m}d\mu &\le 1+z\left( 2+\exp\left(
\sum_{i\le q}\theta (s_i)-m\right)\right). &(3.3.15)\cr}$$

{\bf Step 5} We now turn to lower bounds for $\sum\limits_{i\le q}
\theta (\int_\Omega u_id\mu )$.  For each $i\le q$, consider a number 
$m_i$, and set
$$\leqalignno{w_i(\omega ) &= v(\omega )-S+S_i -\theta (m_i) 
&(3.3.16)\cr
W_i &= \int_\Omega \min (w^+_i,1)d\mu\,. &(3.3.17)\cr}$$

We show that
$$\leqalignno{\int_{\{ w_i\ge 0\}} (u_i-m_i)d\mu &\le {W_i \over 3
\theta '(m_i)}\qquad \left( = - {W_i \over 3\vert\theta '(m_i)\vert}
\right)\,. &(3.3.18)\cr}$$

By (3.3.11), we have
$$\leqalignno{\theta (u_i(\omega )) &\ge v(\omega )-S+S_i = w_i(
\omega )+\theta (m_i)\,. &(3.3.19)\cr}$$

Now, by (3.3.5), for $y\ge x$, we have
$$\xi (y) \le \xi (x)+{1 \over 3}\xi '(x)\min (1, y-x)\,.$$

Taking $x=\theta (m_i)$, $y=x+w_i(\omega )$, combining with (3.3.19), 
and recalling that $\xi '(\theta (m_i))=\theta '(m_i)^{-1}$ yields, 
when $w_i (\omega )\ge 0$ that
$$u_i(\omega )\le m_i+{1 \over 3\theta '(m_i)}\min (1,w_i(\omega ))$$
from which (3.3.18) follows by integration.

{\bf Step 6}  We take $m_i=s_i$.  It follows from (3.3.12), (3.3.10) 
that
$$\int_{\{ u_i\ge s_i\}}(u_i-s_i) d\mu\le{\tau \over \vert\theta '
(s_i)\vert}\,.$$
Combining with (3.3.18), observing that $w_i(\omega )>0$ implies 
$u_i(\omega )<m_i$ by (3.3.19), and using convexity of $\theta$ yield
$$\leqalignno{\theta (\int_\Omega u_id\mu ) &\ge\theta (s_i)-\tau +
{W_i \over 3}\,. &(3.3.20)\cr}$$

We choose the number $m$ of Step 4 as the smallest for which
$$\card\{ i\le q\,;\, S-S_i+\theta (s_i)\le m\}\ge {q \over 2}\,.$$

We observe that if $S-S_i+\theta (s_i)\le m$, then $W_i\ge z$, where 
$W_i$ is given by (3.3.17) and $z$ by (3.3.14).  Thus (3.3.20) shows 
that if we set $R=\sum\limits_{i\le q}\theta (\int_\Omega u_id\mu )$, 
we have
$$\leqalignno{\sum_{i\le q}\theta (s_i) &\le R+q\tau -{q \over 6} z
\,. &(3.3.21)\cr}$$
Combining with (3.3.15) gives
$$\eqalign{\int_\Omega e^{v-m} d\mu &\le 1+z(2+e^{R+q\tau -{q \over 6}
z-m})\cr
&\le 3+ze^{q\tau-{q \over 6} z} e^{R-m}\cr}$$
Calculus show that $\sup\limits_{z} ze^{-qz/6}=6/qe$.  Thus if we 
assume
$$\leqalignno{e^{q\tau} &\le {qe \over 12} &(3.3.22)\cr}$$
we have
$$\int e^{v-m}d\mu\le 3+{1 \over 2} e^{R-m}\,.$$
For $R-m\ge 2$, this is $\le e^{R-m}$, so the proof is finished.

{\bf Step 7}  Thus, we only have to consider the case $R\le m+2$.  By 
definition of $m$, the set
$$I=\{ i\le q\,;\,m\le S-S_i+\theta (s_i)\}$$
has cardinality $\ge q/2$.  For $i$ in $I$, we have
$$R\le m+2\le 2+S+\theta (s_i)-S_i$$
and summation over $i\in I$ yields
$$\leqalignno{R-S &\le 2+{1 \over \card I}\sum_{i\in I}(\theta (s_i)
-S_i)\le 2+{2 \over q}\sum_{i\le q}(\theta (s_i)-S_i) &(3.3.23)\cr}$$
since $\theta (s_i)-S_i\ge 0$ for all $i\le q$.  On the other hand, 
(3.3.21) implies that
$$\sum_{i\le q}(\theta (s_i)-S_i)\le R-S+q\tau$$
and combining with (3.3.23) yield (for $q\ge 3$, $\tau\le 1$) that
$$\leqalignno{\sum_{i\le q}\theta (s_i)-S &\le \left( 1-{2 \over q}
\right)^{-1}(2+q\tau )\le K+q\tau\,. &(3.3.24)\cr}$$

{\bf Step 8}  We assume that $q\ge 3$, $\tau\le 1$, so that (3.3.24) 
holds, and we finish the proof.  In Step 5, we take $m_i=\sup u_i(
\omega )$, so that $\theta (m_i)=S_i$, and $w_i=v-S$ does not depend 
on $i$.  From (3.3.18) and convexity, we get
$$\theta (\int_\Omega u_i d\mu )\ge S_i+{W \over 3}$$
where $W=W_i=\int_\Omega \min (1,w^+)d\mu$.

We now have by summation that
$$\leqalignno{{qW \over 3} &\le R-S\,. &(3.3.25)\cr}$$

In Step 4, we take $m=S$, so that $z=W$.  From (3.3.15), (3.3.24) we 
get
$$\leqalignno{\int_\Omega e^{v-S}d\mu &\le 1+3W\exp (K+q\tau ) 
&(3.3.26)\cr
&\le \exp (3W\exp (K+q\tau ))\,.\cr}$$
According to (3.3.25), this is less than $\exp (R-S)$ provided $\exp 
(K+q\tau )\le q/9$, i.e. $\tau\le q^{-1}\log (q/K)$.  Moreover, this 
requirement implies (3.3.22).

The proof is now complete.\rbx
\enddemo

{\bf 3.4.~~Interpolation}

One can express Proposition 2.1.1 as the fact that, if $P(A)>1/2$, 
then for most of the elements $x$ of $\Omega^N$, all but of order 
$\sqrt {N}$ coordinates can be copied by an element of $A$.  On the 
other hand, Theorem 3.1.1 asserts that for most of the elements $x$ 
of $\Omega^N$, all but a bounded number of coordinates of $x$ can be 
copied by one of two elements of $A$.  A rather natural question is 
whether both phenomenon can be achieved simultaneously (using the
{\it same} elements of $A$).  In this section, we will show that this 
is indeed the case.

This fact seems to be a special case of a rather general phenomenon 
that can be informally formulated as follows:  Suppose we have defined 
two notion of the idea ``the points $x$ and $y$ are within $\ll
\text{distance}\gg t$''; we call these I and II respectively.  Assume 
that there is good concentration of measure when the fattening $A_t$ 
of $A$ is defined as the collection of points $x$ that are within 
distance $t$ of $A$, when the meaning of this is defined with respect 
to notion I (resp. II).  Then, in all the cases we have considered, 
it remains true that we have good concentration of measure when $A_t$ 
is now defined as the collection of points $x$ for which there 
exists a point $y$ which is within distance $t$ of $x$ with respect 
of the two notions {\it simultaneously}.  Two specific examples are 
presented, one in this section, the other in Section 4.5.  In both 
sections, we present an inequality, that quantitatively contains two 
rather separate inequalities presented before.  Considerably more 
difficult (if at all possible) would be the task to find a formulation 
that would allow to recover sharp forms of these two inequalities.  
This direction of finding inequalities that ``merge'' several other 
inequalities is very natural.  It remains at an embryonic stage.  The 
reason is partly the intrinsic difficulty; partly the lack of 
concrete applications that would help to formulate precise needs.

We now go back to question of finding an inequality encompassing at 
the same time the essence of Proposition 2.1.1 and Theorem 3.1.1.  
For simplicity, we consider only the case $q=2$ in Theorem 3.1.1.  
For two subsets $A_1$, $A_2$ of $\Omega^N$, $x\in\Omega^N$, $a,t>0$, 
we set
$$\leqalignno{f(A_1,A_2,a,t,x) &= \inf\{ f(y^1,y^2,a,t,x)\,;\, y^1\in 
A_1\,,\,y^2\in A_2\} &(3.4.1)\cr}$$
where
$$\eqalign{f(y^1,y^2,a,t,x) &= a\card\{ i\le N\,;\,x_i\not= y^1_i\,;
\,x_i\not= y^2_i\}\cr
&\quad +t\card\{ i\le N\,;\,x_i\not= y^1_i~\text{or}~x_i\not= y^2_i
\}\,.\cr}$$

\proclaim{Theorem 3.4.1}  For each $a<\log 2$, there exists $t_0>0$ 
such that
$$\leqalignno{t<t_0\Rightarrow\int e^{f(A_1,A_2,a,t,x)}dP(x) &\le 
{e^{4Nt^2} \over P(A_1)P(A_2)}\,. &(3.4.2)\cr}$$
\endproclaim

In particular, by Chebyshev inequality, this implies that for $u\le 
8Nt^2_0$, we have
$$P\left(\left\{ f\left( A_1,A_2,a,\sqrt{{u \over 8N}},x\right)\ge u
\right\}\right)\le {e^{-u/2}\over P(A_1)P(A_2)}\,.$$
When $f\left( A_1,A_2,a,\sqrt{{u \over 8N}},x\right)\le u$, by 
definition, we can find $y^1\in A_1$, $y^2\in A_2$ such that
$$a\,\card\{ i\le N\,;\, x_i\not\in\{ y^1_i,y^2_i\}\}+\sqrt{{u \over 
8N}} \card\{i\le N\,;\,x_i\not= y^1_i~\text{or}~x_i\not= y^2_i\}
\le u$$
so that in particular
$$\eqalign{\card\{ i\le N\,;\,x_i\not\in\{ y^1_i,y^2_i\}\} &\le {u 
\over a}\cr
\card\{ i\le N\,;\, x_i\not= y^1_i~\text{or}~x_i\not= y^2_i\} &\le 
\sqrt{8Nu}\,.\cr}$$

We would like to point out that the factor $e^{4Nt^2}$ in (3.4.2) is 
not optimal.  This factor can be improved, in particular, with greater 
effort on the calculus computations of the proof we will present.  
Further improvement would be possible as in Section 1.2, but we have 
not pursued that direction since it is not clear at the present time 
what would be an optimal quantitative form of the phenomenon described 
by Theorem 3.4.1.

The key to Theorem 3.4.1 is the following.

\proclaim{Proposition 3.4.2}  Given $b<\log 2$, there exists $t_0>0$ 
such that, if $t<t_0$, for any two functions $g_1,g_2\le 1$ on 
$\Omega$, we have
$$\leqalignno{\int_\Omega\min\left( e^b,{e^t \over g_1(\omega )} , 
{e^t \over g_2(\omega )} , {1 \over g_1(\omega )g_2(\omega )}\right) 
d\mu (\omega ) &\le {e^{4t^2} \over \int g_1d\mu\int g_2d\mu}\,. 
&(3.4.3)\cr}$$
\endproclaim

\demo{Proof}  The relatively simple method we present does not yield 
the optimal dependence in $t$ in the right hand side of (3.4.3), but 
it avoids lengthy unpleasant computations.  Arguing as in the  proof 
of Lemma 3.3.2, we see that
$$\int hd\mu\int g_1d\mu\int g_2d\mu \le \exp\int (h+g_1+g_2-3)d
\mu\,.$$
Thus, if we set
$$h(g_1,g_2) = \min\left( e^b,{e^t \over g_1} , {e^t \over g_2} , {1 
\over g_1g_2}\right)$$
it suffices to show that for $t$ small enough, and all numbers $g_1,
g_2\le 1$, we have
$$\leqalignno{h(g_1,g_2)+g_1+g_2 &\le 3+4t^2\,. &(3.4.4)\cr}$$
\enddemo

Certainly, we can assume $g_1\ge g_2$ and $2t\le b$.

\remark{Case 1} $g_2\le g_1\le e^{t-a}$.  In that case
$$h(g_1,g_2)+g_1+g_2-3\le e^b+2e^{t-b}-3\,.$$
Since $e^b<2$, we have $e^b+2e^{-b}-3<0$, so that we can find $t_0$ 
such that $e^b+2e^{t-b}-3\le 0$ if $t\le t_0$.
\endremark

\remark{Case 2} $g_2\le e^{t-b}\le g_1$.  In that case
$$\eqalign{h(g_1,g_2)+g_1+g_2-3 &\le {e^t \over g_1}+g_1+g_2-3\cr
&\le e^b+2e^{t-b} -3\cr}$$
since the function $x+e^t/x$ decreases for $x\le 1$, and we conclude 
as above.
\endremark

\remark{Case 3}  $e^{t-b}\le g_2\le e^{-t}$.  In that case, using 
again that the function $x+e^t/x$ decreases for $x\le 1$, we have, 
since $g_1\ge g_2$,
$$\eqalign{h(g_1,g_2)+g_1+g_2-3 &\le {e^t \over g_1}+g_1+g_2-3\cr
&\le {e^t \over g_2} +2g_2-3\le e^{2t}+2e^{-t}-3\cr}$$
since the function $2x+e^t/x$ is convex, and thus on the interval 
$[e^{t-b},e^{-t}]$ is bounded by the maximum of its values at the 
endpoints.  Also, we note that $e^{2t}+2e^{-t} -3\le 4t^2$ 
\endremark

\remark{Case 4}  $g_2\ge e^{-t}$.  Then
$$h(g_1,g_2)+g_1+g_2-3\le {1 \over g_1g_2} +g_1+g_2-3\le e^{2t}+2
e^{-t} -3$$
since, when $c>1$, the function $x+c/x$ decreases for $x\le 1$.  We 
then conclude as above.\rbx
\endremark

We will let the reader complete the proof of Theorem 3.4.1 using the 
induction method and Proposition 3.4.2.  The basic observation is 
that, if $B_i$ denotes the projection of $A_i$ on $\Omega^N$, we 
have for $x\in\Omega^N$, $\omega\in\Omega$,
$$\eqalign{f(A_1,A_2,a,t,(x,\omega )) &\le a+t+f(B_1,B_2,a,t,x)\cr
f(A_1,A_2,a,t,(x,\omega )) &\le t+f(B_1,A_2(\omega ),a,t,x)\cr
f(A_1,A_2,a,t,(x,\omega )) &\le t+f(A_1(\omega ),B_2,a,t,x)\cr
f(A_1,A_2,a,t,(x,\omega )) &\le f(A_1(\omega ), A_2(\omega ),a,t,
\omega )\,.\cr}$$
For the induction hypothesis, one then fixes $a<b<\log 2$, and take 
$t_0$ small enough that $a+t_0\le b$.
\vfil\eject

\noindent{\bf 4.~~Convex Hull}

{\bf 4.1.~~The basic result}

The main idea of this section is the introduction of a rather 
different way of measuring how far a point $x$ is from a subset $A$ 
of $\Omega^N$.  We introduce the set
$$U_A(x)=\{ (s_i)_{i\le N}\in\{ 0,1\}^N\,;\,\exists y\in A\,,\,s_i=0
\Rightarrow x_i=y_i\}\,.$$
We denote by $V_A(x)$ the convex hull of $U_A(x)$, when $U_A(x)$ is 
seen as a subset of $\Bbb R^N$.  Thus $V_A(x)$ contains zero if and 
only if $x$ belongs to $A$.  We denote by $f_c(A,x)$ the $\ell^2$ 
distance from zero to $V_A(x)$ (the letter $c$ refers to ``convexity'').  
The corresponding notion of ``enlargement'' of $A$ is as follows:
$$\leqalignno{A^c_t &= \{ x\in\Omega^N\,;\, f_c(A,x)\le t\}\,. 
&(4.1.1)\cr}$$
These notations will be kept throughout the paper.

\proclaim{Theorem 4.1.1} For every subset $A$ of $\Omega^N$, we have
$$\leqalignno{\int \exp {1 \over 4} f^2_c(A,x)dP(x) &\le {1 \over 
P(A)} &(4.1.2)\cr}$$

In particular
$$\leqalignno{P(A^c_t) &\ge 1-{1 \over P(A)}e^{-t^2/4}\,. 
&(4.1.3)\cr}$$
\endproclaim

In order to understand better (4.1.1) it is worthwhile to note the 
following simple result.

\proclaim{Lemma 4.1.2}  The following are equivalent
$$\leqalignno{&x\in A^c_t &(4.1.4)\cr
&\forall (\alpha_i)_{i\le N}\,,\,\exists y\in A\,,\,\sum_{i\le N}
\{\alpha_i\,;\,x_i\not= y_i\}\le t\sqrt{\sum\limits_{i\le N}
\alpha_i^2}\,. &(4.1.5)\cr}$$
\endproclaim

\demo{Proof}  The linear functional $\overline{\alpha}\colon x\to
\sum\limits_{i\le N}\alpha_ix_i$ on $\Bbb R^N$, provided with the 
Euclidean norm, has a norm $\Vert\overline{\alpha}\Vert =
\sqrt{\sum\limits_{i\le N}\alpha^2_i}$.  Since $V_A(x)$ contains a 
point of norm $\le f_c(A,x)$, the infimum of $\overline{\alpha}$ on 
$V_A(x)$ is $\le f_c(A,x)\Vert\overline{\alpha}\Vert$; but since 
$V_A(x)$ is the convex hull of $U_A(x)$, the infimum of 
$\overline{\alpha}$ on $U_A(x)$ is the same as the infimum on 
$V_A(x)$.  Thus (4.1.4) implies (4.1.5).  The converse (that is not 
needed in the paper) follows from the Hahn-Banach theorem.\rbx
\enddemo

It is very instructive to compare (4.1.3) with (2.1.3).  If one takes 
$t=k/\sqrt{N}$, $\alpha_i=1$, one sees that (4.1.3) implies
$$P(f(A,x)\ge k)\le {1 \over P(A)} e^{-k^2/4N}\,.$$
The only difference with (2.1.3) is the worst numerical coefficient 
in the exponential.  But the strength of (4.1.3) is, of course, that 
{\it all} choices of $\alpha_i$ are possible.  This makes Theorem 
4.1.1 a principle of considerable power, as will be demonstrated at 
length in Part II.  It does, however, take some effort to fully 
understand the potential of Theorem 4.1.1.  To illustrate one use of 
Theorem 4.1.1, let us consider the case where $\Omega =\{ 0,1\}$, and 
where the probability $\mu$ gives mass $p$ to $1$ (and mass $1-p$ 
to zero), where $p\le 1/2$.  Consider a subset $A$ of $\{ 0,1\}^N$, 
and assume that $A$ is hereditary, i.e., that if $y=(y_i)_{i\le N}\in 
A$, and if $(z_i)_{i\le N}$ is such that $z_i\le y_i$ for all $i$, 
then $z\in A$.  Consider $x\in\{ 0,1\}^N$, and $J=\{ i\le N\,;\, 
x_j=1\}$.  Set $m(x)=\card J$.  Define $\alpha_i=1$ if $i\in J$, 
$\alpha_i=0$ otherwise.  Then Lemma 4.1.2 shows that we can find $y
\in A$ such that
$$\card\{ i\in J\,;\, x_i\not= y_i\le f_c(A,x)\sqrt{m(x)}\}\,.$$
Since $A$ is hereditary, we have $f(A,x)\le f_c(A,x)\sqrt{m(x)}$.

Thus we have, for all $m'$
$$\eqalign{P(\{ f(A,\cdot )\ge t\}) &\le P(\{f_c(A,\cdot )\ge {t 
\over \sqrt{m'}}\} ) +P(m(y)>m')\cr
&\le {1 \over P(A)}\exp \left(-{t^2 \over 4m'}\right) +P(m(y)>m')
\,.\cr}$$
Since the last term becomes very small for $m'>pN$, we recover the 
correct order $1/Np$ of the coefficient of $t^2$ in (2.3.5).

The key to Theorem 4.1.1 is the following simple lemma.

\proclaim{Lemma 4.1.3}  Consider $0\le r\le 1$.  Then
$$\leqalignno{\inf_{0\le\lambda\le 1} r^{-\lambda}\exp {(1-\lambda 
)^2 \over 4} &\le 2-r\,. &(4.1.5)\cr}$$
\endproclaim

\demo{Proof}  Taking $\lambda =1+2\log r$ if $r\ge e^{-1/2}$, and 
$\lambda =0$ otherwise, and taking logarithms, it suffices to show 
that
$$f(r)=\log (2-r)+\log r+(\log r)^2\ge 0\,.$$
Now $f(1)=0$, so it suffices to show that $f'(r)\le 0$.  Since $f'(1)
=0$, it suffices to show that $(rf'(r))'\ge 0$, or, equivalently, by 
calculation that $(2-r)^{-2}-r^{-1}\le 0$.  But $(2-r)^{-2}\le 1\le 
r^{-1}$.\rbx
\enddemo

We now prove Theorem 4.1.1, by induction upon $N$.  We leave to the 
reader the easy case $N=1$.  For the induction step from $N$ to $N+1$, 
consider a subset $A$ of $\Omega^{N+1}$ and its projection $B$ on 
$\Omega^N$.  For $\omega\in\Omega$, we set as usual
$$A(\omega )=\{ x\in\Omega^N\,;\,(x,\omega )\in A\}\,.$$
Consider $x\in \Omega^N$, $\omega\in\Omega$, $z=(x,\omega )$.  The 
basic observation is that
$$\eqalign{s\in U_{A(\omega )}(x) &\Rightarrow (s,0)\in U_A(z)\cr
t\in U_B(x) &\Rightarrow (t,1)\in U_A(z)\,.\cr}$$
Thus, for $s\in V_{A(\omega )}(x)$, $t\in V_B(x)$, $0\le\lambda\le 
1$, we have $(\lambda s+(1-\lambda )t,1-\lambda )\in V_A(z)$.  The 
convexity of the function $u\to u^2$ shows that
$$\leqalignno{f^2_c(A,z) &\le (1-\lambda )^2+\lambda f^2_c(A(\omega 
),x)+(1-\lambda )f^2_c(B,x)\,. &(4.1.7)\cr}$$
The main trick of the proof is to resist the temptation to optimize 
now over $\lambda$.  By Holder's inequality and induction hypothesis, 
we have
$$\eqalign{\int\exp{1 \over 4}f^2_c &(A,(x,\omega ))dP(x)\cr
&\le \exp{1 \over 4}(1-\lambda )^2(\int_{\Omega^N}\exp 
{1 \over 4}f^2_c(A(\omega ),x)dP(x))^\lambda (\int_{\Omega^N}\exp {1 
\over 4} f^2_c(B,x)dP(x))^{1-\lambda}\cr
&\le \exp{1 \over 4}(1-\lambda )^2\left({1 \over P(A(\omega ))}
\right)^\lambda\left({1 \over P(B)}\right)^{1-\lambda}\cr
&={1 \over P(B)}\exp {1 \over 4}(1-\lambda )^2\left({P(A(\omega )) 
\over P(B)}\right)^{-\lambda}\,.\cr}$$
This inequality holds for all $0\le\lambda\le 1$.  Using (4.1.6) with 
$r=P(A(\omega ))/P(B)\le 1$, we get
$$\int_{\Omega^N}\exp{1 \over 4} f^2_c(A,(x,\omega ))dP(x)\le {1 
\over P(B)}\left( 2-{P(A(\omega )) \over P(B)}\right)\,.$$
Integrating with respect to $\omega$ and using Fubini theorem yields
$$\eqalignno{\int\exp{1 \over 4} f^2_c(A,\cdot )d(P\otimes\mu ) &\le 
{1 \over P(B)}\left( 2-{P\otimes \mu (A) \over P(B)}\right)\cr
&\le {1 \over P\otimes\mu (A)}\,, &\bx\cr}$$
since $x(2-x)\le 1$ for all $x$ real.
\vfil\eject

{\bf 4.2.~~Sharpening}

We now try to improve (4.1.2) by allowing a right hand side $P(A
)^{-\alpha}$ for some $\alpha\ge 0$.  In that case, it will be 
advantageous to measure the ``distance'' of $s$ to $V_A(x)$ by the 
function
$$f_\alpha (A,x) = \inf\{s\in V_A(x)\,;\,\sum_{i\le N}\xi (\alpha ,
s_i)\}$$
where
$$\leqalignno{\xi (\alpha ,u) &= \alpha (1-u)\log (1-u)-(\alpha +1-
\alpha u)\log\left({1+\alpha -\alpha u \over 1+\alpha}\right)\,. 
&(4.2.1)\cr}$$

The reader should observe right away that $f_\alpha (A,x)$ corresponds 
(with the notation of Section 4.1) to $f^2_c(A,x)$ rather than to 
$f_c(A,x)$.  This will be the case for all the extensions of Theorem 
4.1.1 we will consider.

As pointed out, Lemma 4.1.3 is the key to Theorem 4.1.1.  It is a 
somewhat magic fact that when one tries to improve upon Lemma 4.1.3, 
the best possible function that can be used instead of the function 
$(1-\lambda )^2/4$ can be computed exactly, leading to the formula 
(4.2.1).

\proclaim{Lemma 4.2.1}  Consider $0<r<1$.  Then
$$\leqalignno{\inf_{0\le\lambda\le 1}r^{-\lambda\alpha}\exp\xi (
\alpha ,1-\lambda ) &= 1+\alpha -\alpha r\,. &(4.2.2)\cr}$$
\endproclaim

\demo{Proof}  We will not give the shortest possible proof (that 
consists in checking by computation that for $\lambda =r(\alpha +1-
\alpha r)^{-1}$, we have $r^{-\lambda\alpha}\exp\xi (\alpha ,1-\lambda 
)=1+\alpha -\alpha r$).  Rather, we will explain how (4.2.2) was 
discovered.  We fix $\alpha$, and we set $f(x)=\alpha^{-1}\xi (\alpha 
,x)$.  The best choice for $\lambda$ is such that $\alpha\log r+
\alpha f'(1-\lambda )=0$, i.e., $r=\exp (-f'(1-\lambda ))$.  So we 
would like to have the identity, for $0\le\lambda\le 1$,
$$\exp (\alpha f(1-\lambda )+\alpha\lambda f'(1-\lambda )) = 1+\alpha 
-\alpha\exp (-f'(1-\lambda ))\,.$$
Setting $\nu =1-\lambda$, and taking logarithms, we want
$$\alpha f(\nu) +\alpha (1-\nu )f'(\nu )=\log (1+\alpha -\alpha\exp 
(-f'(\nu )))\,.$$
Differentiating in $\nu$ and setting $g(\nu )=\exp (-f'(\nu ))$, we 
get
$$\alpha (1-\nu )f''(\nu )= {\alpha f''(\nu )g(\nu ) \over 1+\alpha 
-\alpha g(\nu )}$$
so that $g(\nu )={(\alpha +1)(1-\nu ) \over \alpha +1 - \alpha\nu}$.  
Taking logarithms and integrating yields (4.2.1).\rbx
\enddemo

\proclaim{Lemma 4.2.2}  The function $\xi (\alpha ,\cdot )$ is 
increasing and convex on $[0,1]$ and $\xi (\alpha ,u)\ge {\alpha 
\over 2(\alpha +1)} u^2$.
\endproclaim

\demo{Proof}  Computation shows that $\xi (\alpha ,0)={d\xi \over du} 
(\alpha ,0)=0$, and
$${d^2\xi \over du^2} ={\alpha \over (\alpha +1-\alpha u)(1-u)}\ge 
{\alpha \over \alpha +1}$$
since $u\ge 0$.\rbx
\enddemo

\proclaim{Lemma 4.2.3}  For $\alpha ,a>0$, we have
$$\leqalignno{1+\alpha -\alpha a &\le a^{-\alpha} &(4.2.3)\cr
a+(1-a)\exp\xi (\alpha ,1) &\le a^{-\alpha}\,. &(4.2.4)\cr}$$
\endproclaim

\demo{Proof}  To prove (4.2.3), we observe that the graph of the 
convex function $x^{-\alpha}$ is above its tangent at the point $x=1$.  
To prove (4.2.4), we observe that $\xi (\alpha ,1)=\log (1+\alpha )$, 
so that the left hand side is
$$a+(1-a)(1+\alpha )=1+\alpha -\alpha a$$
and the result follows from (4.2.3).\rbx
\enddemo

\proclaim{Theorem 4.2.4}  For a subset $A$ of $\Omega^N$, we have
$$\leqalignno{\int^\ast \exp f_\alpha (A,x)dP(x) &\le {1 \over P
(A)^\alpha}\,. &(4.2.5)\cr}$$
\endproclaim

\demo{Proof}  It is an obvious adaptation of the proof of Theorem 
4.1.1.  The case $N=1$ follows from (4.2.4), and (4.2.3) is used as 
a substitute for the last inequality of (4.1.8).\rbx
\enddemo

If we use Lemma 4.2.2, we see that (4.1.3) can be generalized into
$$\leqalignno{P(A^c_t) &\ge 1-{1 \over P(A)^\alpha} \exp\left( -
{\alpha t^2 \over 2(\alpha +1)}\right)\,. &(4.2.6)
\cr}$$

Optimization over $\alpha$ as in Corollary 2.2.3 yields:

\proclaim{Corollary 4.2.5}  For each subset $A$ of $\Omega^N$,
$$\leqalignno{t &\ge \sqrt{2\log{1 \over P(A)}} \Rightarrow P(A^c_t)
\ge 1-\exp\left( -{1 \over 2}\left( t-\sqrt{2\log {1 \over P(A)}}
\right)^2\right)\,. &(4.2.7)\cr}$$
\endproclaim

It is an interesting question whether the term $\sqrt{2\log 1/P(A)}$ 
can be removed in (4.2.7).  We will, however, see in Section 4.3 that 
the coefficient $1/2$ cannot be improved.  It must be pointed out that 
Theorem 4.2.4 bring considerably more than a simple improvement of the
coefficient of $t^2$ in (4.1.3).  The reason is that $\xi (\alpha ,1)
=\log (\alpha +1)$ becomes very large when $\alpha$ is large.  In that 
case, (4.2.5) recovers certain features of (3.1.2), and in some ways, 
improves simultaneously upon Theorem 3.1.1 and Proposition 2.1.1.  To 
see this, consider $q\ge 1$.  We fix $A\subset \Omega^N$, and for $x
\in\Omega^N$, we consider
$$k(x)=\inf\{ k\,;\,\exists s\in V_A(x)\,;\,\card\{i\le N\,;\, s_i\ge 
1-{1 \over q}\}\le k\}\,.$$
Then, certainly, we have
$$k(x)\xi \left( q,1-{1 \over q}\right) \le f_q(A,x)\,.$$
Now,
$$\eqalign{\xi \left( q,1-{1 \over q}\right) &=\log{1 \over q} -2\log 
{2 \over 1+q}\cr
&=\log{(1+q)^2 \over 4q} \ge \log {q \over 4}\cr}$$
so that
$$\leqalignno{k(x)\log {q \over 4} &\le f_q(A,x)\,. &(4.2.8)\cr}$$

On the other hand, by (4.2.5), we have
$$P(f_q(A,x)\ge t)\le {e^{-t} \over P(A)^q}$$
so that, by (4.2.8)
$$\leqalignno{P(k(x)\ge k) &\le {e^{-k\log {q \over 4}} \over P(A)^q} 
= \left({4 \over q}\right)^k {1 \over P(A)^q}\,. &(4.2.9)\cr}$$

The relationship with (3.1.2) is as follows.

If $k(x)\le k$, we can find a family $(y^j)_{j\le m}$ of points of 
$A$, and coefficients $(\alpha_j)_{j\le m}$, $0\le \alpha_j\le 1$, 
$\sum\limits_{j\le m}\alpha_j =1$, such that
$$\leqalignno{\card\left\{ i\le N\,;\,\sum_{j\le m}\alpha_j1_{\{x_i
\not= y^j_i\}}\ge 1-{1 \over q}\right\} &\le k\,. &(4.2.10)\cr}$$
On the other hand, if $f(A,\dots ,A,x)\le k$, we can find $y^1,\dots ,
y^q$ in $A$ such that
$$\leqalignno{\card\left\{ i\le N\,;\,\sum_{j\le q}{1 \over q}1_{\{ 
x_i\not= y^j_i\}}> 1-{1 \over q}\right\} &\le k\,. &(4.2.11)\cr}$$

Certainly (4.2.11) is more precise than (4.2.10); however, for some 
important applications (see [T3]) (4.2.10) is just as powerful as 
(4.2.11).

{\bf 4.3.~~Two point space}

In this section, we consider the case where $\Omega =\{ 0,1\}$ and 
where $\mu$ gives weights $1-p$ to zero and $p$ to $1$.  The miracle 
of Lemma 4.2.1 does not seem to happen again, so we will only 
consider statements of the type
$$\leqalignno{\int\exp f_u(A,x)dP(x) &\le {1 \over P(A)^\alpha} 
&(4.3.1)\cr}$$
where, for a couple $u=(u_0,u_1)$ of positive numbers, we set
$$f_u(A,x) =\inf\left\{u_0\sum\{ s^2_i;x_i=0\}+u_1\sum\{s^2_i;x_i=1\} 
\colon s\in V_A(x)\right\}\,.$$
In other words, we take into account the fact that the points $0$ and 
$1$ do not play the same role.

If one analyzes the arguments of Sections 2.3, 4.1, 4.2, one sees 
that the best value the induction method allows to take for $u_0$ is 
the largest number $s$ such that, whenever $a<b$, we have
$$(1-p)\inf_{0\le\lambda\le 1}{1 \over a^{\alpha (1-\lambda )}}~{1 
\over b^{\alpha\lambda}} e^{\lambda^2s}+{p \over b^\alpha}\le {1 
\over ((1-p)a+pb)^\alpha}$$
or equivalently
$$\leqalignno{(1-p){1 \over x^\alpha}\inf_{0\le\lambda\le 1}
x^{\alpha\lambda}e^{\lambda^2s} &\le {1 \over ((1-p)x+p)^\alpha}-p 
&(4.3.2)\cr}$$
for all $0\le x\le 1$.

(The best possible value of $u_1$ is obtained in a similar way, 
changing $p$ in $1-p$, and will not be considered.)

The infimum in (4.3.2) is obtained for
$$\lambda =\max \left( 0,-{\alpha\log x \over 2s}\right)\,.$$
The left hand side of (4.3.2) is constant for $x\le\exp (-2s/\alpha 
)$; thereby (4.3.2) holds provided, for $x\ge\exp (-2s/\alpha )$, we 
have
$$\leqalignno{(1-p)e^{-{\alpha^2 \over 4s}(\log x)^2} &\le {1 \over (
(1-p)x+p)^\alpha}-p &(4.3.3)\cr}$$

Determining the best value of $s$ for which this holds is an 
unpleasant task, so we will content ourselves with finding good 
values of $s$.  Taking logarithms and differentiating, one sees
that (4.3.3) will hold provided, for $x>0$, we have
$$\leqalignno{-{\alpha \over 2s} \log x &\ge 1-{(1-p)x \over A-p
A^{\alpha +1}} = {p(1-A^{\alpha +1}) \over A-pA^{\alpha +1}} 
&(4.3.4)\cr}$$
where we have set $A=(1-p)x+p$.

It suffices that for $x\ge 0$ we have
$$\leqalignno{-{\alpha \over 2s}\log x &\ge {p \over 1-p}~{1-
A^{\alpha +1} \over A}\,. &(4.3.5)\cr}$$

We first consider the case $p ={1 \over 2}$, and we show that in this 
case we can take $s={\alpha \over \alpha +1}$.  Since
$$\leqalignno{1-A^{\alpha +1} &\le (\alpha +1)(1-A)=(\alpha +1)(1-p)
(1-x) &(4.3.6)\cr}$$
it suffices to see that
$$0<x\le 1 \Rightarrow \log x\ge {2(1-x) \over 1+x}\,.$$
But the function $f(x)=\log x -{2(1-x) \over 1+x}$ satisfies $f(1)=0$,
$f'(x)=-(1-x)^2/(1+x)^2\le 0$.  Using the notation $f_c(A,x)$ of 
Section 4.1.1, we then have proved the following.

\proclaim{Theorem 4.3.1}  When $\Omega =\{ 0,1\}$ and $\mu$ is 
uniform, for each $\alpha\ge 1$ and each subset $A$ of $\Omega^N$, we 
have
$$\leqalignno{\int \exp \left({\alpha \over \alpha+1} f^2_c(A,x)
\right) dP(x) &\le {1 \over P (A)^\alpha}\,. &(4.3.7)\cr}$$
\endproclaim

Compared with (4.2.6), we have gained a factor $2$ in the exponent in 
the special case of the two point space.

\proclaim{Corollary 4.3.2}  When $\Omega =\{ 0,1\}$, and $\mu$ is 
uniform, for each $\alpha\ge 1$, and each subset of $\Omega^N$, we 
have
$$\leqalignno{t &\ge \sqrt{\log{1 \over P(A)}} \Rightarrow P(A^c_t)
\ge 1 - \exp\left( -\left(t-\sqrt{\log{1 \over P(A)}}\right)^2\right)
\,. &(4.3.8)\cr}$$
\endproclaim

\demo{Proof}  From (4.3.7) and Chebyshev inequality, we get
$$P(A^c_t) \ge 1-{1 \over P(A)^\alpha}\exp\left( -{\alpha \over 
\alpha +1} t^2\right)$$
and we optimize over $\alpha$ as in the proof of Corollary 2.2.3.

It is a natural question whether (4.3.8) can be improved into
$$\leqalignno{P(A^c_t) &\ge 1-K\exp (-t^2)\,. &(4.3.9)\cr}$$

It should, however, be pointed out that the coefficient of $t^2$ is 
optimal.  We will now show this, and at the same times, the optimality 
of the coefficient $1/2$ in (4.2.7).  Provide $\Omega =\{ 0,1\}$ 
with the probability $\mu$ that gives mass $p$ to $1$.  Set
$$A=\left\{ (x_i)_{i\le N}\,;\,\sum_{i\le N}x_i\le pN\right\}\,.$$
(Thus, for $N$ large, $P(A)$ is about $1/2$.)  Consider $y\in\{ 0,
1\}^N$, such that $\card J=m$, where $J=\{i\le N\,;\,y_i=1\}$.  Assume 
$m>pN$.  Then any element $x$ of $A$ differs of $y$ in at 
least $m-pN$ of the coordinates indexed by $J$.  Using Lemma 4.1.2 
for $\alpha_i=1/\sqrt{m}$ when $i\in J$, $\alpha_i=0$ otherwise, we 
see that
$$\leqalignno{f_c(A,y) &\ge {(m-pN) \over \sqrt{m}} ={m-pN \over 
\sqrt{N}}\sqrt{{N \over m}}\,. &(4.3.10)\cr}$$

If we think of $m=m(y)$ as a r.v., the central limit theorem shows 
that, as $n\to\infty$, $(m-pN)/\sqrt{N}$ is asymptotically normal, 
with standard deviation $\sqrt{p(1-p)}$.  On the other hand, 
$\sqrt{N/m}$ converges to $\sqrt{1/p}$ in probability.  Thus
$$\eqalign{\lim_{N\to\infty} P(f_c(A,\cdot )\ge t) &\ge {1 \over
\sqrt{2\pi}}\int^\infty_{t/\sqrt{(1-p)}}\exp \left(-{u \over 2}
\right)^2 du\cr
&\ge {1 \over Kt}\exp\left( - {t^2 \over 2(1-p)}\right)\,.\cr}$$
If $p=1/2$, the coefficient of $t^2$ is $-1$; and if we let $p$ 
arbitrary, we cannot do better than the coefficient $-1/2$ of (4.2.7).

We now go back to our main line of discussion, and we consider the 
case $p\le 1/2$; we will show that in this case we can take
$$\leqalignno{s &= \min\left({\alpha \over K}\log {1 \over p}\,,\,
{\alpha \over 4(\alpha +1)p}\right)\,. &(4.3.11)\cr}$$
In particular, for $\alpha$ large, this is of order $1/p$, rather 
than order $\log (1/p)$.  This remarkable fact is closely connected 
to Theorem 4.4.1 below.  To prove (4.3.11), we prove (4.3.5), 
depending on the value of $x$.

{\bf Case 1.}  $x\ge 1/2$.  Then $-\log x\ge 1-x$, $1-A^{\alpha +1}
\le (\alpha +1)(1-p)(1-x)$, so that it suffices that
$$-{\alpha \over 2s}\ge {p(1+\alpha ) \over A}\,.$$
Now $A\ge {1 \over 2}$, so that it suffices that $s\le {\alpha \over 
4p(1+\alpha )}$.

{\bf Case 2.}  $x\le \sqrt{p}$.  Then $-\log x\ge -{1 \over 2}\log 
1/p$, so that it suffices that
$${\alpha \over 4s}\log \left({1 \over p}\right) \ge { p \over (1-p)
A}\,.$$
Since $A\ge p$, it suffices that
$$s\le {\alpha (1-p) \over 4}\log {1 \over p}\,.$$

{\bf Case 3.}  $\sqrt{p}\le x\le {1 \over 2}$.  It then suffices, 
since $-\log x\ge \log 2$, $A\ge (1-p)\sqrt{p}$ that
$$s\le {(1-p)^2 \over \sqrt{p}}~{\alpha \over 2}\log 2$$
which holds when $s\le{\alpha\over K}\log {1 \over p}$.

{\bf 4.4.~~Penalties}

We consider now a function $h$ on $\Omega\times\Omega$, such that $h
\ge 0$, and $h(\omega ,\omega )=0$ for $\omega\in\Omega$.  For a 
subset $A$ of $\Omega^N$, and $x\in\Omega^N$, we set
$$\leqalignno{U_A(x) &= \{ (s_i)\in\Bbb R^N_+\,;\,\exists y\in A\,;\,
\forall i\le N\,,\, s_i\ge h(x_i,y_i)\}\,. &(4.4.1)\cr}$$
We denote by $V_A(x)$ the convex hull of $U_A(x)$.  The situation of 
Section 4.1 corresponds to the case where $h(\omega ,\omega ')=1$ if 
$\omega\not= \omega '$.

In order to measure the ``distance'' of zero to $V_A (x)$, we consider 
a convex function $\psi$ on $\Bbb R$, with $\psi (0)=0$.  We will 
assume
$$x\le 1 \Rightarrow \psi (x)\le x^2\,;\qquad\qquad x\ge 1\Rightarrow
\psi (x)\ge x\,.$$

We set
$$f_{h,\psi}(A,x)=\inf\left\{\sum_{i\le N}\psi (s_i)\,;\,s=(s_i)_{i
\le N}\in V_A(x)\right\}\,.$$
(Thus, the situation of Section 4.1 corresponds to the case $\psi (s)
=s^2$.)  The material of this section is connected to that of Section 
2.6, and the notations of Section 2.6 are in force in the present 
section.  Thus $\theta$ denotes a convex function from $]0,1]$ to 
$\Bbb R^+$, with $\theta (1)=0$, $\lim\limits_{x\to 0}\theta (x)=
\infty$, and $\xi$ denotes the inverse function.  We assume that 
(2.6.1) holds, and assume moreover that for a certain number $\gamma 
>0$, we have
$$\leqalignno{b\ge 0 &\Rightarrow \vert\xi '(b+1)\vert\ge\gamma\vert
\xi '(b)\vert &(4.4.2)\cr}$$
$$\leqalignno{\vert\theta '(1)\vert &\ge\gamma\,;\qquad w(1/2)\ge
\gamma\,. &(4.4.3)\cr}$$

We recall the function $\Xi$ of (2.6.2), as well as condition $H(\xi 
,w)$ of (2.6.3).

\proclaim{Theorem 4.4.1}  Consider a nonincreasing function $w$ on 
$]0,1]$, $w\le\theta$.  Assume that $\int^1_0 w^2d\lambda\le 1$, and 
that condition $H(\xi ,w)$ holds.  Assume that for each subset 
$B$ of $\Omega$, we have
$$\leqalignno{0<\mu (B) &\le {1 \over 2}\Rightarrow\int_\Omega\exp
\psi (h(x,B))d\mu (x)\le \exp w(\mu (B)) &(4.4.4)\cr}$$
$$\leqalignno{\mu (B) &\ge {1 \over 2}\,,~~t\ge 1\Rightarrow\mu (\{ x
\,;\,\psi (h(x,B))\ge t\})\le e^{-t}(1-\mu (B))\,. &(4.4.5)\cr}$$
Then, for each subset $A$ of $\Omega^N$, we have
$$\leqalignno{\int_{\Omega^N}\exp {1 \over K} f_{h,\psi} (A,x)dP(x) 
&\le \exp \theta (P(A)) &(4.4.6)\cr}$$
where $K$ depends on $\gamma$ only.
\endproclaim

We should observe first that only the values of $w(x)$ for $x\le 1/2$ 
matter.

In order to compare Theorem 4.4.1 with Theorems 2.6.5 and 2.7.1, we 
first have to keep in mind that it is the function $\psi\circ h$ here 
that plays the role of $h$ in these theorems.  The conclusion 
of Theorem 4.4.1 is stronger than that of Theorem 2.6.5 (the way 
Theorem 4.1.1 improves on Proposition 2.1.1) but weaker than the 
conclusion of Theorem 2.7.1 (since one takes convex hulls).  
Condition (4.4.5) strongly resembles (2.7.2).  Condition (4.4.4) 
coincides with Condition (2.6.12) when $\mu (B)\le 1/2$.  A simple 
calculation using (4.4.5) shows that for $\mu (B)\ge 1/2$, condition 
(4.4.5) is of a somewhat stronger nature than (2.6.12).

An interesting case where it is worth to spell out (4.4.4) and 
(4.4.5) is when $h(x,y)=h(y)$ depends on $y$ only.  Denoting by $m$ 
a median of $h$, (4.4.5) will hold if $\psi (m)<1$.  And, as seen 
after Theorem 2.6.5, (4.4.4) holds provided $w(\mu (\{ h\ge t\} ))
\ge\psi (t)$ (a tail condition of $h$).

To prove Theorem 4.4.1 when $N=1$, we observe that, since $w\le
\theta$, (4.4.6) follows from (4.4.4) when $\mu (B)\le 1/2$.  When 
$\mu (B)\ge 1/2$, a simple computation using (4.4.5) shows that given 
$\gamma$, if $K$ is large enough, then
$$\int_\Omega\exp {1 \over K}\psi (h(x,B))d\mu (x)\le 1+\gamma (1-\mu 
(B)) \le \theta (\mu (B))\le \exp \theta (\mu (B))$$
since $\theta '(1)\ge\gamma$.

For the induction step, comparison with the proof of Theorem 4.1.1 
shows that it suffices to prove the following (used for $g=\xi (f)$).

\proclaim{Proposition 4.4.2}  There exists a constant $L$, depending 
on $\gamma$ only, with the following property.  Under the conditions 
of Theorem 4.4.1, consider a function $f\ge 0$ on $\Omega$.  Set
$$\leqalignno{\widehat{f}(x) &= \inf_{y\in\Omega\,,\,0\le\lambda\le 1}
(\lambda f(x)+(1-\lambda )f(y)+{1 \over L}\psi ((1-\lambda )h(x,y)))
\,. &(4.4.7)\cr}$$
Then we have
$$\leqalignno{\int e^{\widehat f}d\mu &\le e^{\theta (\int\xi (f)d
\mu )}\,. &(4.4.8)\cr}$$
\endproclaim

Understandably, with the level of generality considered here, the 
proof cannot be very short.  The reason why we have opted for great 
generality is that Theorem 4.4.1 is a principle of considerable 
power (as will be demonstrated in Chapter 8) and that thereby it 
seems worthwhile to prove extensions of it under weak hypothesis on 
the function $h$.  The proof will incorporate in particular ideas 
from Theorems 4.1.1, 2.6.5, 2.7.1.

A simple idea is that we will need to control $\theta (\int\xi (f)d
\mu )$ from below.  This means controlling the lower tail of $f$.  Set 
$B_s=\{f\le s\}$, and denote by $m$ a median of $f$, so that $\mu (B_m
)\ge 1/2$.  We set
$$\leqalignno{b &= \inf_{s\le m}\left\{ s+{1 \over L}w(\mu (B_s))
\right\}\,. &(4.4.9)\cr}$$

The first step of the proof will be to show that $\mu (B_s)$ is not 
too big, i.e. that $b$ is not too small.

\proclaim{Proposition 4.4.3}  To prove Proposition 4.4.2, if $L> 4/
\gamma$, we can assume
$$\leqalignno{m &\le b+{4 \over L\gamma}\,. &(4.4.10)\cr}$$
\endproclaim

\demo{Proof}  We assume $m>b$, for otherwise there is nothing to 
prove.  Using (4.4.7) with $\lambda =0$, we see that for each $s$ we 
have $\widehat{f}(x)\le s+L^{-1}\psi\circ h(x,B_s)$.  Using (4.4.4) 
together with H\"older's inequality, it follows that
$$\int_\Omega e^{\widehat f}d\mu\le \exp \left(s+{1 \over L}w(\mu 
(B_s))\right)$$
so that $\int_\Omega e^{\widehat f}d\mu\le e^b$ by taking the infimum 
over $s\le m$.  On the other hand, (4.4.9) implies
$$\leqalignno{s\le m &\Rightarrow b-s\le{1 \over L}w(\mu (B_s)) 
&(4.4.11)\cr}$$
i.e.
$$\leqalignno{\left|\left\{ {1 \over L} w\ge b-s\right\}\right| &\ge 
\mu (B_s)\,. &(4.4.12)\cr}$$

We can hence appeal to Lemma 2.6.4 with $C=\{ f<b\}$ and $t=1/L$ to 
see that
$$\int_C\xi (f)d\mu\le\mu (C)\xi (b)+\xi '(b)\int_C(f-b)d\mu +{1 
\over L^2}\vert\xi '(b)\vert\,.$$
But, by (4.4.12) we have
$$\int_C\vert f-b\vert d\mu\le {1 \over L}\int wd\mu\le {1 \over L}
(\int w^2d\mu )^{1/2}\le {1 \over L}$$
and thus
$$\leqalignno{\int_C\xi (f)d\mu &\le\mu (C)\xi (b)+{2 \over L}
\vert\xi '(b)\vert\,. &(4.4.13)\cr}$$

On the other hand, when $f(\omega )>b$, we have
$$\xi (f(\omega ))\le\xi (b)-(\xi (b)-\xi (m))1_{\{ f\ge m\}}(\omega 
)$$
and, by integration, since $\mu (\{ f\ge m\} )\ge 1/2$, we have (since 
$m>b$)
$$\int_{\Omega\backslash C}\xi (f)d\mu\le (1-\mu (C))\xi (b)-{1 \over 
2} (\xi (b)-\xi (m))\,.$$
Combining with (4.4.13) we get
$$\int_\Omega\xi (f)d\mu\le\xi (b)+{2 \over L}\vert\xi '(b)\vert -{1 
\over 2}(\xi (b)-\xi (m))\,.$$

Since we have shown that $\int_\Omega e^{\widehat f}d\mu\le e^b$, 
there is nothing to prove unless $\int_\Omega\xi (f)d\mu\ge\xi (b)$ 
(for otherwise $\theta (\int\xi (f)d\mu )\ge b$).  Thus we can 
assume
$$\leqalignno{{1 \over 2}(\xi (b)-\xi (m)) &\le {2 \over L}\vert\xi 
'(b)\vert\,. &(4.4.14)\cr}$$
Now, since $m>b$, from (4.4.2) follows that $\xi (m)\le\xi (b)-\gamma
\vert\xi '(b)\vert\min ((m-b),1)$.  Comparing with (4.4.14), we see 
that $\min (m-b,1)\le 4/(L\gamma )$, so that if $L> 4/\gamma$, 
we must have $m-b\le 4/L\gamma$.\rbx
\enddemo

We consider the smallest number $\alpha$ for which
$$\forall s\le m\,,\qquad m-s\le\alpha w(\mu (B_s))$$
so that
$$\leqalignno{\forall s\le m\,,\qquad \vert\{\alpha w\ge m-s\}\vert 
&\ge \mu (B_s)\,. &(4.4.15)\cr}$$

It is rather important to note that
$$\leqalignno{\alpha &\le {8 \over L\gamma^2}\,. &(4.4.16)\cr}$$
Indeed, if $m-s\le 8/L\gamma$, then 
$${m-s \over w(\mu (B_s))} \le {m-s \over w(1/2)}\le {8 \over L
\gamma^2}\,.$$
On the other hand, if $m-s\ge 8/L\gamma$, then, by (4.4.10), we have 
$m-s\le 2(b-s)$, so that
$${m-s \over w(\mu (B_s))} \le 2{b-s \over w(\mu (B_s))} \le {2 \over 
L}\,.$$

We consider a second parameter $M\le L$.  Throughout the rest of this 
section, we will have to put conditions on $L$, $M$, $L/M$.  For 
simplicity we make the convention that the expression ``If $L$ 
is large enough''... means ``there exists a constant $K(\gamma )$, 
depending on $\gamma$ only, such that, if $L\ge K(\gamma )$...'' and 
similarly for $M$, $L/M$.

We set $m'=m- 16/L\gamma^2$.  We consider the function
$$f'=\min \left( f,m+{1 \over M}\right)$$
and the function $g$ defined as
$$\eqalign{g(\omega ) &= f(\omega )\qquad~~\qquad\text{if}~~f(\omega 
)\le m'\cr
g(\omega ) &= \max\left(m',\min\left(\widehat{f}(\omega ),m+{1 \over
M}\right)\right)\qquad\text{if}~~
f(\omega )>m'\,.\cr}$$
Since $\widehat{f}\le f$, it is simple to see that $g\le f'$.  It is 
also simple to see that
$$\leqalignno{g(\omega )\not= f'(\omega ) &\Rightarrow g(\omega )\,,
\,f'(\omega )\in\left[ m',m+{1 \over M}\right]\,. &(4.4.17)\cr}$$
Indeed, the right hand side does not occur only when $f(\omega )< m'$, 
and then $f(\omega )=f'(\omega )=g(\omega )$.  We set
$$C = \{ f\ge m\}\,;\qquad D=\left\{ f\ge m+{1 \over M}\right\}\,.$$

\proclaim{Lemma 4.4.4}  We have
$$\leqalignno{\theta\left(\int_\Omega\xi (f)d\mu\right) &\ge 
\int_\Omega f'd\mu -\alpha^2 - \int_C(f'-m)^2d\mu\,. &(4.4.18)\cr}$$
\endproclaim

\demo{Proof}  Since $f'\le f$, we have $\xi (f')\ge\xi (f)$ and
$$\theta\left(\int_\Omega\xi (f)d\mu\right) \ge \theta\left(
\int_\Omega\xi (f')d\mu \right)\,.$$

We now appeal to Lemma 2.6.4 with $t=\alpha$.  We have
$$\int_\Omega\xi (f')\le\xi (m)+\xi '(m)\int_\Omega (f'-m)d\mu +
\alpha^2\vert\xi '(m)\vert +\xi ''(m)\int_C(f'-m)^2 d\mu\,.$$
By convexity of $\theta$ and since $\xi ''(m)\le\vert\xi '(m)\vert$ 
this implies
$$\eqalignno{\theta\left(\int_\Omega\xi (f)\right) &\ge m+\int_\Omega 
(f-m)d\mu -\alpha^2 -\int_C (f'-m)^2 d\mu\,. &\bx\cr}$$
\enddemo

\proclaim{Lemma 4.4.5}  If $L$ and $M$ are large enough, we have
$$\leqalignno{&~~\quad~~\int_\Omega e^{\widehat f}d\mu &(4.4.19)\cr
&~~\qquad~~ \le \exp\left({1 \over 2}\int_\Omega (g+f')d\mu +2
\alpha^2+2\int_C (f'-m)^2d\mu
+\int_\Omega(e^{\widehat{f}-m}-e^{1/M})^+d\mu\right)\,.\cr}$$
\endproclaim

\demo{Proof}  First, we observe that
$$\leqalignno{\int_\Omega e^{\widehat{f}-m}d\mu &\le
\int_\Omega\exp\left(\min\left(\widehat{f}-m,{1 
\over M}\right)\right)d\mu +\int_\Omega (e^{\widehat{f}-m}-e^{1/M})^+
d\mu\,. &(4.4.20)\cr}$$

We observe that $\min\left(\widehat{f}-m,{1 \over M}\right)\le g-m$.  
Since $e^x\le 1+x+x^2$ for $x\le 1/M\le 1$, we have
$$\leqalignno{\int_\Omega\exp &\left(\min\left(\widehat{f}-m,{1 
\over M}\right)\right)d\mu &(4.4.21)\cr
&\le\int_\Omega e^{g-m}d\mu \le 1+\int_\Omega (g-m)d\mu +\int_\Omega
(g-m)^2d\mu\,.\cr}$$
Now, by (4.4.17), and provided $L$, $M$ are large enough,
$$\leqalignno{(g-m)^2 &\le 2(f'-m)^2+2(g-f')^2\le 2(f'-m)^2+{1 \over 
2}(f'-g)\,. &(4.4.22)\cr}$$
We recall also that
$$\int_{\Omega\backslash C}(f'-m)^2d\mu\le\int_{\Omega\backslash
C}(f-m)^2d\mu\le\alpha^2\,.$$
The result follows by combining these inequalities, and using that 
$1+x\le e^x$.\rbx
\enddemo

It follows from Lemmas 4.4.4 and 4.4.5 that to prove Proposition 
4.4.2, it suffices to prove the following when $M$, $L/M$ are large 
enough.
$$\leqalignno{\int_\Omega (f'-g)d\mu &\ge 6\alpha^2+6\int_C(f'-m)^2
d\mu +2\int_\Omega (e^{\widehat{f}-m}-e^{1/M})^+d\mu\,. &(4.4.23)\cr}$$

This follows from the next three lemmas.

\proclaim{Lemma 4.4.6}  We have
$$\leqalignno{\int_\Omega (e^{\widehat{f}-m}-e^{1/M})^+d\mu &\le {K 
\over L}\mu (D)\le {KM^2 \over L}\int_C (f'-m)^2d\mu\,. &(4.4.24)\cr}$$
\endproclaim

\proclaim{Lemma 4.4.7} If $L/M$ is large enough, we have
$$\int_\Omega (f'-g)d\mu\ge {M \over K}\int_C(f'-m)^2d\mu\,.$$
\endproclaim

\proclaim{Lemma 4.4.8}  If $L/M$ is large enough, we have
$$\int_\Omega (f'-g)d\mu\ge {L\gamma^4\alpha^2 \over K}\,.$$
\endproclaim

\demo{Proof of Lemma 4.4.6}  The definition of $\widehat{f}$ (with 
$\lambda =1$) shows that $\widehat{f}(\omega )\le m+1/M+L^{-1}\psi (h
(\omega,\Omega\backslash D))$.  Thus by (4.4.5) we have
$$\mu\left(\left\{\widehat{f}\ge m+{1 \over M} +{k \over L}\right\}
\right) \le e^{-k}\mu (D)$$
and thus
$$\int_\Omega (e^{\widehat{f}-m}-e^{1/M})^+d\mu \le\sum_{k\ge 1} 
e^{{1 \over M}}(e^{{k \over L}}-1)e^{-k+1}\mu (D)$$
from which the first inequality of (4.4.24) follows by elementary 
estimates.  (The second inequality of (4.4.24) is obvious.)\rbx
\enddemo

\demo{Proof of Lemma 4.4.7}  {\bf Step 1.}  For $k\ge 0$, we define
$$a_k=\sup\left\{ t\,;\,\mu (\{ f'\ge t\} )\ge {1 \over 2e^k}\right\}
\,.$$
Thus $m\le a_k\le a_{k+1}\le m+1/M$.  We consider a set $Z_k\subset
\{ f\le a_k\}$ such that
$$\mu (Z_k) = 1-{1 \over 2e^k}\,.$$
We set $Z'_k=\{\omega\,;\,h(\omega ,Z_k)\le 2\}$.  Since $\psi (x)\ge 
x$ for $x\ge 1$, we have
$$Z'_k\supset\{\omega\,;\,\psi (h(\omega ,Z_k))\le 2\}$$
so that by (4.4.5) we have $\mu (Z'_k)\ge 1-1/2e^{k+2}$.  We set, 
for $k\ge 0$
$$W_k=Z'_k\cap (Z_{k+2}\backslash Z_{k+1})\,.$$
We observe that the sets $(W_k)_{k\ge 0}$ are disjoint, and that
$$\leqalignno{\mu (W_k) &\ge {1 \over 2e^k}\left({1 \over e}-{2 \over 
e^2}\right) \ge {1 \over 2e^{k+3}}\,. &(4.4.25)\cr}$$

{\bf Step 2.}  We show that
$$\leqalignno{\int_{W_k}(f'-g)d\mu &\ge {M \over K} (a_{k+1}-a_k)^2
\mu (W_k\backslash D)\,. &(4.4.26)\cr}$$

Consider $\omega\in W_k\backslash D$.  Then $f'(\omega )=f(\omega )$, 
so that given $\lambda\in [0,1]$, $\omega '\in\Omega$
$$\leqalignno{f'(\omega )-\widehat{f}(\omega ) &= f(\omega )-
\widehat{f}(\omega ) &(4.4.27)\cr
&\ge (1-\lambda )(f(\omega )-f(\omega '))-{1 \over L}\psi ((1-
\lambda )h(\omega ,\omega '))\,.\cr}$$
We can find $\omega '\in Z_k$ such that $h(\omega ,\omega ')\le 3$.  
Then $f(\omega )-f(\omega ')\ge a_{k+1}-a_k$.  We can take $0\le
\lambda\le 1$ such that $1-\lambda =M(a_{k+1}-a_k)/3$.  Then (4.4.27) 
yields, since $\psi (x)\le x^2$ for $x\le 1$, that
$$f'(\omega )-\widehat{f}(\omega )\ge {M \over 3}(a_{k+1}-a_k)^2-{9
M^2 \over L}(a_{k+1}-a_k)^2\,.$$
Thus, if $L/M$ is large enough,
$$f'(\omega )-\widehat{f}(\omega )\ge {M \over 4}(a_{k+1}-a_k)^2\,.$$
Thus
$$\widehat{f}(\omega )\le f'(\omega )-{M \over 4}(a_{k+1}-a_k)^2\,.$$
Since $a_{k+1}-a_k\le 1/M$, and $f'(\omega )\ge a_{k+1}$, the 
right-hand side is $\ge a_k\ge m$, so that
$$g(\omega )\le f'(\omega )-{M \over 4}(a_{k+1}-a_k)^2$$
and thus
$$f'(\omega )-g(\omega )\ge {M\over 4}(a_{k+1}-a_k)^2$$
from which (4.4.26) follows by integration.

{\bf Step 3.}  Denote by $k_0$ the largest integer such that $1/4
e^{k_0+3}\ge \mu (D)$.  Thus $\mu (W_k)\ge 2\mu (D)$ for $k\le k_0$, 
and by (4.4.26) and summation, we get, since $\mu (W_k\backslash 
D)\ge\mu (W_k)/2$:
$$\leqalignno{\int (f'-g)d\mu &\ge {M \over K}\sum_{k\le k_0}(a_{k+1}
-a_k)^2e^{-k}\,. &(4.4.28)\cr}$$
By the argument of Lemma 2.7.8, we have
$$\sum_{k\le k_0} (a_{k+1}-a_k)^2 e^{-k}\ge {1 \over K}\int_C (\min
(f',m+a_{k_0+1})-m)^2d\mu\,.$$
Thus the proof is completed if $a_{k_0+1}\ge 1/2M$.

{\bf Step 4.}  Assuming now $a_{k_0+1}\le 1/2M$, we show that
$$\leqalignno{\int (f'-g)d\mu &\ge {1 \over KM} e^{-k_0}\,. 
&(4.4.29)\cr}$$
Since
$$e^{-k_0}\ge {M^2 \over K}\sum_{k>k_0}(a_{k+1}-a_k)^2e^{-k}\,,$$
combining with (4.4.28), we get
$$\int (f'-g)d\mu \ge {M \over K}\sum_{k\ge 0}(a_{k+1}-a_k)^2e^{-k}
\ge {M \over K}\int_C(f-m)^2 d\mu$$
by (the argument of) Lemma 2.7.8, completing the proof of Lemma 4.4.7.

To prove (4.4.29), we observe that by definition of $k_0$ we have 
$e^{-k_0-6}\le \mu (D)$.  Consider the set
$$Z = \{\psi (h(\cdot ,Z_{k_0+1}))\le 6\}\,.$$
Then, by (4.4.5), we have $\mu (Z)\ge 1-e^{-k_0-7}$, so that $\mu 
(Z\cap D)\ge e^{-k_0}/K$.  Now, if $\omega\in D$, we have $f'(\omega 
)\ge m+1/M$ while if $\omega\not\in D$, we have
$$\eqalign{\widehat{f}(\omega ) &\le a_{k_0+1}+m+{1 \over L}\psi (h(
\omega ,Z_{k_0+1}))\cr
&\le m+{1 \over 2M}+{6 \over L}\,.\cr}$$
Thus $g(\omega )\le m+{3 \over 4M}$ if $L/M$ is large enough.  Hence, 
$f-g\ge 1/4M$ on $Z\cap D$.\rbx
\enddemo

\demo{Proof of Lemma 4.4.8}  {\bf Step 1.}  We show that we can 
assume $\mu (D)\le 1/8$.  Indeed otherwise by Lemma 4.4.7 we have 
$\int (f'-g)d\mu\ge 1/KM$ and, since $\alpha\le 8/\gamma^2L$, this 
is $\ge L\alpha^2$ when $L/M$ is large enough.

{\bf Step 2.}  By definition of $\alpha$, there exists $s<m$ with 
$m-s >\alpha w(\mu (B_s))/2$.  By (4.4.5) and Chebyshev inequality, 
the set
$$H=\{\psi (h(\cdot ,B_s))\le 2+w(\mu (B_s))\}$$
has measure $\ge 3/4$.  Thus if we set $G=H\cap (C\backslash D)$, we 
have $\mu (G)\ge 1/8$.

{\bf Step 3.}  Set
$$\beta ={m-s \over 3+w(\mu (B_s))}\,.$$
Since $w(\mu (B_s))\ge w(1/2)\ge \gamma$, and $m-s\ge\alpha w(\mu 
(B_s))/2$, we have
$${\gamma\alpha \over K}\le\beta\le\alpha\,.$$
Since $\mu (G)\ge 1/8$, it suffices to show that
$$\leqalignno{\forall\omega\in G\,,\quad f'(\omega )-g(\omega ) 
&\ge {L\gamma^2 \over 8}\beta^2\,. &(4.4.30)\cr}$$

{\bf Step 4.}  We prove (4.4.30).  Consider $\omega\in G$.  Then 
$f'(\omega )=f(\omega )\ge m$.  Consider $\omega '\in B_s$ with 
$h(\omega ,\omega ')\le 3+w(\mu (B_s))$.  Then
$$\leqalignno{f(\omega )-\widehat{f}(\omega ) &\ge\sup_{0\le\lambda
\le 1} ((1-\lambda )(m-s)-{1 \over L}\psi ((1-\lambda )(3+w(\mu 
(B_s))))\,. &(4.4.31)\cr}$$

We choose $0\le \lambda\le 1$ such that
$$1-\lambda ={L\gamma^2 \over 4}~{\beta \over 2+w(\mu (B_s))}\,.$$
This is possible since $\beta\le\alpha\le 8/L\gamma^2$.  Then 
(4.4.31) yields, since $\psi (x)\le x^2$ for $x\le 1$, that
$$f(\omega )-\widehat{f}(\omega )\ge {L\gamma^2\beta^2 \over 8}\,.$$
Thus $\widehat{f}(\omega )\le f(\omega )-L\gamma^2\beta^2/8$.  Since 
the right hand side is $\ge m'$, we have $g(\omega )\le f(\omega )-
L\gamma^2\beta^2/8$.  The proof is complete.
\enddemo

{\bf 4.5.~~Interpolation}

The result of this section will interpolate between (a weak form of) 
Theorem 3.1.1, for $q=2$, and (a weak form of) Theorem 4.1.1.  
Consider three points $x=(x_i)_{i\le N}$, $y^1=(y^1_i)_{i\le N}$, 
$y^2 =(y^2_i)_{i\le N}$ of $\Omega^N$.  Set
$$r_i(x,y^1,y^2)=(1_{\{x_i\not= y^1_i\}},1_{\{x_i\not=
y^2_i\}},1_{\{x_i\not\in\{y^1_i,y^2_i\}\}})
\,.$$
Thus $r_i(x,y^1,y^2)\in\{ 0,1\}^3$.  Set
$$r(x,y^1,y^2)=(r_i(x,y^1,y^2))_{i\le N}\in (\{ 0,1\}^3)^N\,.$$
Given two subsets $A_1$, $A_2$ of $\Omega^N$, let
$$U_{A_1,A_2}(x) = \{ r(x,y^1,y^2)\,;\, y^1\in A_1\,,\, y^2\in A_2\}
\,,$$
and consider the convex hull $V_{A_1,A_2}(x)$ of $U_{A_1,A_2}(x)$, 
when $U_{A_1,A_2}(x)$ is seen as a subset of $(\Bbb R^3)^N$.

Throughout this section, we set $b=1/6$, $a=\log (3-2e^{-b})$.  We 
make the convention to write a point $r\in (\Bbb R^3)^N$ as $(r_{1,i},
r_{2,i},r_{3,i})_{i\le N}$.  We set
$$f(A_1,A_2,x)=\inf\left\{\sum_{i\le N}ar^2_{1,i}+ar^2_{2,i}+br_{3,i}
\,;\,r\in V_{A_1,A_2}(x)\right\}\,.$$

\proclaim{Theorem 4.5.1}  We have
$$\int_{\Omega^N}\exp f(A_1,A_2,x)dP(x)\le {1 \over P(A_1)P(A_2)}\,.$$
\endproclaim

To understand better this statement, set $u=f(x,A_1,A_2)$.  Consider 
$r\in V_{A_1,A_2}(x)$ such that
$$\sum_{i\le N}ar^2_{1,i}+ar^2_{2,i}+br_{3,i}\le u\,.$$
Consider numbers $(c_{1,i})_{i\le N}$, $(c_{2,i})_{i\le N}$.  Then, 
for $j=1,2$
$$\eqalign{\sum_{i\le N} c_{j,i}r_{j,i} &\le \left(\sum_{i\le N} 
c^2_{j,i}\right)^{1/2}\left( \sum_{i\le N} r^2_{j,i}\right)^{1/2}\cr
&\le \left({u \over a}\right)^{1/2}\left(\sum_{i\le N}c^2_{j,i}
\right)^{1/2}\,.\cr}$$
Thus
$$\sum_{i\le N}(c_{1,i}r_{1,i}+c_{2,i}r_{2,i}+br_{3,i})\le t=:u+
\left({u \over a}\right)^{1/2}\left(\sqrt{\sum\limits_{i\le N}
c^2_{1,i}}+\sqrt{\sum\limits_{i\le N}c^2_{2,i}}\right)\,.$$

If we recall that $V_{A_1,A_2}(x)$ is the convex hull of $U_{A_1,A_2}
(x)$, this implies that we can find $y^1\in A_1$, $y^2\in A_2$ such 
that
$$\sum\{ c_{1,i}\,;\,x_i\not= y^1_i\} +\sum\{c_{2,i}\,;\,x_i\not= 
y^2_i\}+b\card\{ i\,;\,x_i\not\in \{y^1_i,y^2_i\}\}\le t\,.$$

The proof of Theorem 4.5.1 goes by induction over $N$.  The case 
$N=1$ is left to the reader.  For the induction from $N$ to $N+1$, 
one observes, with the usual notations, that, when $a_{0,0},a_{1,0},
a_{0,1},a_{1,1}\ge 0$, are of sum one, then
$$\eqalign{f(A_1,A_2,(x,\omega )) &\le a_{0,0}f(A_1(\omega ),A_2(
\omega ),x)+a_{1,0}f(B_1,A_2(\omega ),x)\cr
&\quad +a_{0,1} f(A_1(\omega ),B_2 ,x)+a_{1,1}f(B_1,B_2,x)\cr
&\quad +b(a_{1,0}+a_{1,1})^2+b(a_{0,1}+a_{1,1})^2+aa_{1,1}\,.\cr}$$
Thereby, to perform the induction it suffices to show that, when 
$g_1$, $g_2$ are two functions on $\Omega$, $g_1,g_2\le 1$, then
$$\leqalignno{\int\inf\exp &(aa_{1,1}+b(a_{0,1}+a_{1,1})^2+b(a_{1,0}
+a_{1,1})^2){1 \over (g_1g_2
)^{a_{0,0}}}~{1\over g^{a_{1,0}}_2}~{1 \over g^{a_{0,1}}_1}d\mu 
&(4.5.1)\cr
&\quad \le {1 \over \int g_1d\mu\int g_2 d\mu}\cr}$$
where the infimum is taken over all the allowed choices of $a_{0,0}$, 
$a_{0,1}$, $a_{1,0}$, $a_{1,1}$.

\proclaim{Lemma 4.5.2}  We have
$$\leqalignno{\inf\exp &(aa_{1,1}+b(a_{1,0}+a_{1,1})^2+b(a_{0,1}+
a_{1,1})^2){1 \over (g_1g_2)^{a_{0,0}}}~
{1 \over g^{a_{1,0}}_2}~{1 \over g^{a_{0,1}}_1}&(4.5.2)\cr
&\quad \le \sqrt{(3-2g_1)(3-2g_2)}\,.\cr}$$
\endproclaim

We first use (4.5.2) to prove (4.5.1).  By (4.5.2) and Cauchy-Schwarz, 
the left-hand side of (4.5.1) is bounded by
$$\sqrt{\int (3-2g_1)d\mu\int (3-2g_2)d\mu} = \sqrt{(3-2\int g_1d\mu)
(3-2\int g_2d\mu )}\,.$$
Thus it suffices to observe that for $0\le x\le 1$, we have $3-2x
\le x^{-2}$, which expresses the fact that the convex function 
$x^{-2}$ is above its tangent at $x=1$.

\demo{Proof of Lemma 4.5.2}  We will actually restrict the infimum to 
the cases $a_{1,1}=1$ or $a_{1,1}=0$.  We will prove
$$\leqalignno{\min\left( e^a ,\inf_{0\le a_1+a_2\le 1}e^{b(a^2_1+
a^2_2)}{1 \over g^{1-a_1}_1}~{1 \over g^{1-a_2}_2}\right) &\le 
\sqrt{(3-2g_1)(3-2g_2)}\,. &(4.5.3)\cr}$$
We distinguish cases.

{\bf Case 1.}  $g_1g_2\le e^{-2b}$.

It suffices to see that
$$e^a\le \sqrt{(3-2g_1)(3-2g_2)}\,.$$
The right hand side has minimum at $g_1=g_2=e^{-b}$, and our value of 
$a$ has been chosen so that inequality holds in that case.

{\bf Case 2.}  $g_1g_2\ge e^{-2b}$.

For $j=1,2$, we take $a_j=-{\log g_j \over 2b}$.  The purpose of the 
condition $g_1g_2\ge e^{-2b}$ is to ensure $a_1+a_2\le 1$.  It 
suffices to show that
$${1 \over g_1}e^{-{(\log g_1)^2 \over 4b}}\le \sqrt{3-2g_1}\,.$$
We will show that, for $0\le x\le 1$, we have
$$e^{-(\log x)^2/2b}\le x^2(3-2x)$$
or, equivalently that
$$\varphi (x)={(\log x)^2 \over 2b} +2\log x+\log (3-2x)\ge 0\,.$$
Since $\varphi '(0)=0$, $\varphi (0)=0$ it suffices to show that 
$(x\varphi '(x))'\ge 0$, i.e.
$${1 \over b} -{6x \over (3-2x)^2}\ge 0\,.$$
But, since ${1 \over b}=6$, it suffices to show that $x\le (3-2x)^2$, 
which is true since $x\le 1$, $(3-2x)^2\ge 1$.\rbx
\enddemo

\vfil\eject

\noindent{\bf 5.~~The Symmetric Group}

We denote by $S_N$ the group of permutations of $\{ 1,\dots N\}$.  
Our interest in the symmetric group stems from the fact that it is 
closely related to a product.  To see this, let us denote by 
$t_{i,j}$ the transposition of $i$ and $j$.  Then, it is easily seen 
that every $\sigma\in S_N$ can be written in a unique way as
$$\leqalignno{\sigma &= t_{N,i(N)}\circ t_{N-1,i(N-1)}\cdots \circ 
t_{2,i(2)} &(5.1)\cr}$$
where, for $j\le N$, we have $i(j)\le j$.  This decomposition allows 
to transfer some of the results of Chapter 2 to $S_N$.  The purpose 
of the present chapter is to prove a version of Theorem 4.1.1 for 
$S_N$.  The reason for which this is not such an easy task is that 
the decomposition (5.1) is highly noncommutative.

For a subset $A$ of $S_N$, and $\sigma\in S_N$, we set
$$U_A(\sigma )=\{ s\in\{ 0,1\}^N\,;\,\exists\tau\in A\,;\,\forall\ell
\le N\,,\,s_\ell =0\Rightarrow \tau (\ell )=\sigma (\ell )\}$$
and we consider the convex hull $V_A(\sigma )$ of $U_A(\sigma )$ in 
$[0,1]^N$.  We set
$$f(A,\sigma )=\inf\left\{\sum_{\ell\le N}s^2_\ell\,;\, s=(s_\ell )
\in V_A(\sigma )\right\}\,.$$
We denote by $P_N$ the canonical (= homogenous) probability on $S_N$.

\proclaim{Theorem 5.1}  For every subset $A$ of $S_N$ we have
$$\leqalignno{\int_{S_N}\exp{1 \over 16} f(A,\sigma )dP_N(\sigma ) 
&\le {1 \over P_N(A)}\,. &(5.2)\cr}$$
\endproclaim

In a natural way, $S_N$ can be considered as a subset of $\{ 1,\dots 
,N\}^N$, by the map $\sigma\to (\sigma (i))_{i\le N}$.  If $S_N$ were 
equal to all of $\{ 1,\dots ,N\}^N$, (5.2) would be a consequence of 
Theorem 4.1.1; but $S_N$ is only a very small subset of $\{ 1,\dots
,N\}^N$.

The challenge of Theorem 5.1 is that it is apparently not possible to 
prove (5.2) by induction over $N$.  Rather, we will use a stronger 
induction hypothesis.  Given $p\le N$, we set
$$f(A,\sigma ,p) =\inf\left\{ s^2_p+\sum_{\ell\le N}s^2_\ell\,;\, s
\in V_A (\sigma )\right\}\,.$$
Theorem 5.1 is obviously a consequence of the following.

\proclaim{Proposition 5.2}  For each subset $A$ of $S_N$ and each $p
\le N$, we have
$$\leqalignno{\int_{S_N}\exp{1 \over 16} f(A,\sigma ,p)dP_N(\sigma ) 
&\le {1 \over P_N(A)} &(5.3)_N \cr
\int_{S_N}\exp {1 \over 16}f(A,\sigma ,\sigma^{-1}(p))dP_N(\sigma ) 
&\le {1 \over P_N(A)}\,. &(5.4)_N\cr}$$
\endproclaim

We let to the reader to prove Proposition 5.2 when $N=1$.  We now 
assume that Proposition 5.2 has been proved for $N$ and we prove it 
for $N+1$.  A noticeable feature of this proof is that the proof 
of (5.3)$_{N+1}$ (resp. (5.4)$_{N+1}$) will require the use of 
(5.4)$_N$ (resp. (5.3)$_N$).  Before the proof starts, we need to 
introduce some notation.  Given $p,m\le N$, $p\not= m$, we set
$$\leqalignno{f(A,\sigma ,p,m) &= \inf\left\{ s^2_p+\sum_{\ell\le 
N}s^2_\ell\,;\,s\in V_A(\sigma )\,,\,s_m=0\right\}\,. &(5.5)\cr}$$

Given $i,j\le N$, we set
$$\leqalignno{g(A,\sigma ,i,j) &= \inf\left\{\sum_{\ell\not= i,j}
s^2_\ell\,;\, s\in V_A(\sigma )\right\}\,. &(5.6)\cr}$$
We start the proof of (5.4)$_{N+1}$.  Certainly there is no loss of 
generality to assume that $p=N+1$.

\proclaim{Lemma 5.3}  Consider $i,j\le N+1$, $i\not= j$, $\sigma\in 
S_{N+1}$, $0\le\lambda\le 1$.  Then
$$\leqalignno{f(A,\sigma ,i) &\le 4(1-\lambda )^2+(1-\lambda )g(A,
\sigma ,i,j)+\lambda f(A,\sigma ,j,i)\,. &(5.7)\cr}$$
\endproclaim

\demo{Proof}  Consider $s\in V_A(\sigma )$, $t\in V_A(\sigma )$, with 
$t_i=0$.  By convexity of $V_A (\sigma )$, we have
$$u=(1-\lambda )s+\lambda t\in V_A(\sigma )\,.$$
Thus
$$f(A,\sigma, i)\le\sum_{\ell\not= i}u^2_\ell +2u^2_i\,.$$
Since $s_i\le 1$ and since $t_i=0$, we have
$$f(A,\sigma ,i)\le\sum_{\ell\not= i,j}u^2_\ell +2(1-\lambda )^2+((
1-\lambda )s_j+\lambda t_j)^2\,.$$
Since $s_j\le 1$, we have
$$((1-\lambda )s_j+\lambda t_j)^2\le 2(1-\lambda )^2 s^2_j+2
\lambda^2 t^2_j\le 2(1-\lambda )^2+2 \lambda t^2_j\,.$$
Since the function $x\to x^2$ is convex, we have
$$u^2_\ell\le (1-\lambda )s^2_\ell +\lambda t^2_\ell\,.$$
Thus we have
$$f(A,\sigma ,i)\le (1-\lambda )\sum_{\ell\not= i,j}s^2_\ell +\lambda 
\left( 2t^2_j+\sum_{\ell\le N} t^2_\ell \right) +4(1-\lambda )^2\,.$$
The result follows by taking the infimum over $s,t$.\rbx
\enddemo

Following the idea of Theorem 4.1.1, (5.7) will be used together with 
Holder's inequality.  Some work is, however, needed to relate the 
resulting terms to the induction hypothesis.  For $i\le N+1$, we set
$$G_i=\{\sigma\in S_{N+1}\,;\,\sigma (i)=N+1\}\,.$$

For simplicity, we denote by $t_i=t_{N+1,i}$ the transposition of 
$N+1$ and $i$.  We consider the map $R\colon\rho\to\rho\circ t_i$.  
We observe that, if $\rho\in G_i$, then
$$R(\rho )(N+1)=\rho\circ t_i(N+1)=\rho (i)=N+1\,.$$
Thereby, we can consider $R$ as a map from $G_i$ to $S_N$.  We set 
$A_i=A\cap G_i$.

\proclaim{Lemma 5.4}  If $\sigma\in G_i$, we have
$$\leqalignno{f(A,\sigma ,j,i) &\le f(R(A_i),R(\sigma ),t_i(j))\,. 
&(5.8)\cr}$$
\endproclaim

\demo{Proof}  We let the reader consider the essentially obvious case 
where $i=N+1$, and we assume $i\not= N+1$.  Given a sequence $s\in
\{ 0,1\}^N$, we consider the sequence $\overline{s}=(\overline{s}_\ell 
)\in\{ 0,1\}^{N+1}$ defined by $\overline{s}_i=0$, $\overline{s}_{N+1}
=s_i$, $\overline{s}_\ell =s_\ell$ if $\ell\not= i,N+1$.  We note 
that $\overline{s}_\ell =s_{t_i(\ell )}$ for $\ell \not= i$.  Thus it 
suffices to prove that $\overline{s}\in U_A(\sigma )$ whenever $s\in 
U_{R(A_i)}(R(\sigma ))$.  Consider $s\in U_{R(A_i)}(R(\sigma ))$.  By 
definition, there exists $\tau\in R(A_i)$ such that, for $\ell\le N$
$$s_\ell =0\Rightarrow\tau (\ell )=R(\sigma )(\ell )\,.$$
Since $\tau\in R(A_i)$, we have $\tau =R(\rho )$ for a certain $\rho
\in A_i$.  Thus
$$\leqalignno{s_\ell &= 0\Rightarrow \rho (t_i(\ell ))=\sigma (t_i
(\ell ))\,. &(5.9)\cr}$$
We will show that, for $\ell\le N+1$
$$\overline{s}_\ell =0\Rightarrow\rho (\ell )=\sigma (\ell )\,.$$
This holds for $\ell =i$, since $\rho (i)=\sigma (i)=N+1$.  For $\ell
\not= i$, this follows from (5.9), since $\overline{s}_\ell =s_{t_i(
\ell )}$, and $t_i\circ t_i$ is the identity of $S_N$.\rbx
\enddemo

We denote by $Q_i$ the uniform probability on $G_i$.

\proclaim{Corollary 5.5}
$$\leqalignno{\int\exp {1 \over 16} f(A,\sigma ,j,i)dQ_i(\sigma )\le 
{1 \over Q_i(A_i)} &= {1 \over Q_i(A)}\,. &(5.10)\cr}$$
\endproclaim

\demo{Proof}  Using (5.8), the left-hand side of (5.10) is bounded by
$$\eqalign{\int\exp{1 \over 16} f(R(A_i),R(\sigma ),t_i(j))dQ_i(
\sigma ) &=\int\exp {1 \over 16} f(R (A_i),\rho ,t_i(j))dP_N(\rho )\cr
&\le {1 \over P_N(R(A_i))} = {1 \over Q_i(A)}\cr}$$
using (5.3)$_N$.\rbx
\enddemo

\proclaim{Lemma 5.6}  Assume $j\not= i$.  Then
$$\leqalignno{\int\exp{1 \over 16} g(A,\sigma ,i,j)dQ_i(\sigma ) &\le 
{1 \over Q_j(A)}\,. &(5.11)\cr}$$
\endproclaim

\demo{Proof}  The map $S\colon\rho\to\rho\circ t_{ij}$ is one to one 
from $G_j$ to $G_i$.  We will prove that setting $B=R(S(A_j))$, we 
have
$$\leqalignno{g(A,\sigma ,i,j) &\le f(B,R(\sigma )) &(5.12)\cr}$$
where we recall that $R$ is seen as a map from $G_i$ to $S_N$.  Since 
$P_N(B)=Q_j(A)$, (5.11) will follow from either (5.3)$_N$ or 
(5.4)$_N$ as in the proof of Corollary 5.5.

Given a sequence $s\in\{ 0,1\}^N$, we consider the sequence 
$\overline{s}\in\{ 0,1\}^{N+1}$ defined as follows.  We set 
$\overline{s}_i=\overline{s}_j=1$.  If $N+1\not= i,j$, we set
$\overline{s}_{N+1}= s_i$.  If $\ell\not\in\{ i,j,N+1\}$, we set 
$\overline{s}_\ell =s_\ell$.

We will show that when $s\in U_B(R(\sigma ))$, then $\overline{s}
\in U_A(\sigma )$.  By definition of $U_B(R(\sigma ))$, there exists 
$\tau\in B$ such that
$$s_\ell =0\Rightarrow\tau (\ell )=R(\sigma )(\ell )=\sigma\circ 
t_i(\ell )\,.$$
Since $\tau\in B$, we can write $\tau =\rho\circ t_{ij}\circ t_i$, 
where $\rho\in A_j$.  Thus
$$s_\ell =0\Rightarrow\rho\circ t_{ij}\circ t_i(\ell )=\sigma\circ 
t_i(\ell )\,.$$
We will show that for $\ell\le N+1$ we have
$$\overline{s}_\ell =0\Rightarrow\rho (\ell )=\sigma (\ell )\,.$$

The only nontrivial case is $\ell=N+1$, when $N+1\not= i,j$.  In that 
case, when $\overline{s}_{N+1}=0$, we have $s_i=0$, so that $\tau (i)
=R(\sigma )(i)=\sigma (N+1)$.  But
$$\tau (i)=\rho\circ t_{ij}\circ t_i(i)=\rho\circ t_{ij}(N+1)=\rho 
(N+1)\,.$$
since $N+1\not= i,j$.\rbx
\enddemo

We now complete the proof of (5.4)$_{N+1}$.  We select $j$ such that 
$Q_j(A)$ is maximum.  If $i\le N+1$, $i\not= j$, for $0\le\lambda
\le 1$, we have, using Lemmas 5.3, 5.4, Corollary 5.5 and Holder's 
inequality
$$\eqalign{\int\exp{1 \over 16} f(A,\sigma ,i) Q_i(\sigma ) &\le 
\exp\left[ {1 \over 4}(\lambda -1)^2 \right]{1 \over Q_i(A
)^\lambda}~{1 \over Q_j(A)^{1-\lambda}}\cr
&= {1 \over Q_j(A)}\left({Q_i(A) \over Q_j(A)}\right)^{-\lambda} 
\exp{1 \over 4}(1-\lambda )^2\,.\cr}$$
If we appeal to Lemma 4.1.3, we have
$$\leqalignno{\int \exp {1 \over 16} f(A,\sigma ,i)dQ_i(\sigma ) 
&\le {1 \over Q_j(A)}\left( 2-{Q_i(A) \over Q_j(A)}\right)\,. 
&(5.13)\cr}$$

It should be obvious from the induction hypothesis that (5.13) still 
hold for $i=j$.  Since $P_{N+1}=\sum_{i\le N+1}{1 \over N+1}Q_i$,
we have, from (5.13), and since $i=\sigma^{-1}(N+1)$ for $\sigma\in 
G_i$ that
$$\eqalignno{\int\exp{1 \over 16}f(A,\sigma ,\sigma^{-1}(N+1))dP_{N+1}
(\sigma ) &\le {1 \over Q_j(A)}\left( 2-{P_{N+1}(A) \over Q_j(A)}
\right)\cr
&\le {1 \over P_{N+1}(A)}\,. &\bx\cr}$$
Having proved (5.4)$_{N+1}$, we turn towards the proof of 
(5.3)$_{N+1}$.  We can assume again $p=N+1$.  The proof is not 
identical to that of (5.4)$_{N+1}$, but is completely parallel.

\proclaim{Lemma 5.7}  For $\sigma\in S_{N+1}$, $j\le N+1$, $j\not=
\sigma (N+1)$, $0\le\lambda\le 1$, we have
$$\leqalignno{f &(A,\sigma ,N+1)&(5.14)\cr
&\quad \le 4(1-\lambda )^2+(1-\lambda )g(A,\sigma ,N+1,\sigma^{-1}
(j))+\lambda f(A,\sigma ,\sigma^{-1}(j),N+1)\,.\cr}$$
\endproclaim

\demo{Proof}  This is (5.7) if one replaces $i$ by $N+1$, $j$ by 
$\sigma^{-1}(j)$.\rbx
\enddemo

We set
$$G'_i=\{\sigma\in S_{N+1}\,;\,\sigma (N+1)=i\}\,.$$
We fix $i$, and we consider the map $R'\colon\rho\to t_i\circ\rho$.  
Thus, for $\rho\in G'_i$, we have $R'(\rho )(N+1)=t_i(i)=N+1$, and 
we can view $R'$ as a map from $G'_i$ to $S_N$. We set $A'_i=A\cap 
G'_i$.

\proclaim{Lemma 5.8}  If $\sigma\in G'_i$, $i\not= j$, we have
$$\leqalignno{f(A,\sigma ,\sigma^{-1}(j),N+1) &\le f(R'(A'_i),
R'(\sigma ),R'(\sigma )^{-1}(t_i(j)))\,. &(5.15)\cr}$$
\endproclaim

\demo{Proof}  Given a sequence $s\in\{ 0,1\}^N$, we consider the 
sequence $\overline{s}=(\overline{s}_\ell )\in\{ 0,1\}^{N+1}$ defined 
by $\overline{s}_\ell =s_\ell$ if $\ell\not= N+1$, and 
$\overline{s}_{N+1}=0$.  Since $\sigma^{-1}(j)=R'(\sigma )^{-1}(t_i
(j))\not= N+1$, it suffices to prove that $\overline{s}\in U_A(\sigma 
)$ whenever $s\in U_{R'(A'_i)}(R'(\sigma ))$.  Thus, consider $s$ in 
this later set.  By definition, there exists $\tau\in R'(A'_i)$ such 
that
$$\forall\ell\le N\,,~~s_\ell =0\Rightarrow\tau (\ell )=R'(\sigma )(
\ell )\,.$$
Since $\tau\in R'(A'_i)$, we have $\tau =R'(\rho )$, $\rho\in A'_i$.  
Thus,
$$\forall\ell\le N\,,~~s_\ell =0\Rightarrow t_i\circ\rho (\ell )=t_i
\circ\sigma (\ell )\Rightarrow \rho (\ell )=\sigma (\ell )\,.$$
Since $\rho (N+1)=\sigma (N+1)=i$, we then have
$$\forall\ell\le N+1\,,~~\overline {s}_\ell =0\Rightarrow \rho (\ell 
)=\sigma (\ell )\,.$$
Thus $\overline{s}\in U_A(\sigma )$.\rbx
\enddemo

We denote by $Q'_i$ the homogeneous probability on $G'_i$.

\proclaim{Corollary 5.9}  If $j\not= i$,
$$\leqalignno{\int\exp{1 \over 16} f(A,\sigma ,\sigma^{-1}(j),N+1)d
Q'_i(\sigma ) &\le {1 \over Q'_i(A)}\,. &(5.16)\cr}$$
\endproclaim

\demo{Proof}  Using (5.15) and the fact that $R'$ transports $Q'_i$ 
to $P_N$, the left-hand side of (5.16) is bounded by
$$\int\exp{1 \over 16} f(R'(A_i),\rho ,\rho^{-1}(t_i(j)))dP_N(\rho ) 
\le {1 \over P_N(R'(A_i))} = {1 \over Q'_i(A)}$$
using (5.4)$_N$.\rbx
\enddemo

\proclaim{Lemma 5.10}  If $i\not= j$, we have
$$\leqalignno{\int\exp{1 \over 16} g(A,\sigma ,N+1,\sigma^{-1}(j)) d
Q'_i(\sigma ) &\le {1 \over Q'_j(A)}\,. &(5.17)\cr}$$
\endproclaim

\demo{Proof}  The map $S'\colon\rho\to t_{ij}\circ\rho$ is one to one 
from $G'_j$ to $G'_i$.  We will prove that, setting $B=R'\circ S'
(A_j)$, we have, for $\sigma$ in $G'_i$ that
$$\leqalignno{ g(A,\sigma ,N+1,\sigma^{-1}(j)) &\le f(B,R'(\sigma )) 
&(5.18)\cr}$$
where we recall that $R'$ is seen as a map from $G'_i$ to $S_N$.  
Since $P_N(B)=Q'_j(A)$, (5.17) will then follow from either (5.3)$_N$ 
or (5.4)$_N$.

Given a sequence $s\in\{ 0,1\}^N$, we consider the sequence 
$\overline{s}\in\{ 0,1\}^{N+1}$ defined as follows.  We set 
$\overline{s}_{N+1}=\overline{s}_{\sigma^{-1}(j)}=1$.  We set 
$\overline{s}_\ell =s_\ell$ if $\ell\not\in\{ N+1,\sigma^{-1}(j)\}$.  
To prove (5.18) it suffices to prove that if $s\in U_B(R'(\sigma ))$, 
then $\overline{s} \in U_A(\sigma )$.  Thus, consider $s\in U_B(R'(
\sigma ))$.  By definition, there exists $\tau\in B$ such that
$$\leqalignno{s_\ell = 0 &\Rightarrow \tau (\ell )=R'(\sigma )(\ell 
) = t_i\circ\sigma (\ell )\,.&(5.19)\cr}$$
Since $\tau\in B$, we can write $\tau =t_i\circ t_{ij}\circ\rho$, 
where $\rho\in A'_j$.  Thus, by (5.19)
$$s_\ell = 0 \Rightarrow t_{ij}\circ\rho (\ell ) = \sigma (\ell )
\Rightarrow \rho (\ell )=t_{ij}\circ\sigma (\ell )\,.$$
Now, for $\ell\not= N+1$, $\sigma^{-1}(j)$, we have $\sigma (\ell )
\not= i,j$; thus $t_{ij}\circ\sigma (\ell ) = \sigma (\ell )$.  Thus 
for these values of $\ell$ we have
$$\eqalignno{\overline{s}_\ell &= 0\Rightarrow s_\ell =0\Rightarrow
\rho (\ell ) = \sigma (\ell )\,. &\bx\cr}$$
\enddemo

The end of the proof of (5.3)$_{N+1}$ is similar to the end of the 
proof of (5.4)$_{N+1}$, and is left to the reader.
\vfil\eject

\noindent{\bf 6.~~Bin Packing}

Given a collection $x_1,\dots ,x_N$ of items, of sizes $\le 1$, the 
bin packing problem requires finding the minimum number $B_N(x_1,
\dots ,x_N)$ of unit size bins in which the items $x_1,\dots ,x_N$ 
can be packed, subject to the restriction that the sum of the sizes 
of items attributed to a given bin cannot exceed one.  (For 
simplicity, we will denote items and item sizes by the same
letters.)   The bin packing problem is a fundamental question of 
computer science, and, accordingly, has received considerable 
attention.  Much work has been done on stochastic models [C-L].  In 
the model we will consider, the items $X_1,\dots ,X_N$ are 
independently distributed according to a given distribution $\mu$.  
One of the natural questions that arises is the study of the 
fluctuations of the random variable $B_N(X_1,\dots ,X_N)$.  One early 
result, [R-T1], [McD1], using martingales, is that for all $t>0$, one 
has
$$\leqalignno{P(\vert B_N(X_1,\dots ,X_N)-EB_N(X_1,\dots ,X_N)\vert
\ge t) &\le 2\exp\left(-{2t^2 \over N}\right)\,. &(6.1)\cr}$$
However, especially when $EX_1$ is small, one expects that the 
behavior of $B_N(X_1,\dots ,X_N)$ resembles the behavior of 
$\sum\limits_{i\le N}X_i$.  Thereby one should expect that the exponent 
in the right-hand side of (6.1) should be of order $t^2/N\var (X_1)$, 
or, at least, less ambitiously, $t^2/NE(X^2_1)$.  This is apparently 
not so easy to prove, and despite several attempts, was established 
only recently using non-trivial bin-packing theory [R1], [R2], [R3].  
The purpose of the present section is to prove this result as an 
application of Theorem 4.1.1.  Several features of the proof will 
appear repeatedly in future applications.  One advantage of our 
approach is that it uses only trivial facts about bin packing, such 
as the following observation.

\proclaim{Lemma 6.1}  We have
$$B_N(x_1,\dots ,x_N)\le 2\sum_{i\le N} x_i+1\,.$$
\endproclaim

\demo{Proof}  It suffices to construct a packing in which at most 
one bin is less than half full.  Such a packing exists since bins 
that are less than half full can be merged.\rbx
\enddemo

We take $\Omega =[0,1]$.  For a subset $A$ of $\Omega^N$, and $x\in 
\Omega^N$, we recall the notation $f_c(A,x)$ introduced in Section 
4.1.  For $x=(x_1,\dots ,x_N)\in\Omega^N$, we write simply $B_N(x)$ 
rather than $B_N(x_1,\dots ,x_N)$.  For $x\in\Omega^N$, we set
$\Vert x\Vert_2=\left(\sum\limits_{i\le N}x^2_i\right)^{1/2}$.  
Finally, for $a>0$, we set
$$A(a)=\{ y\in\Omega^N\,;\,B_N(y)\le a\}\,.$$
The crucial observation is as follows.

\proclaim{Lemma 6.2}  For all $x\in\Omega^N$, we have
$$\leqalignno{B_N(x) &\le a+2\Vert x\Vert_2 f_c(A(a),x)+1\,. 
&(6.2)\cr}$$
\endproclaim

\demo{Proof}  As follows from Lemma 4.1.2, (taking $\alpha_i$ there 
equal to $x_i$) we can find $y\in A(a)$ such that, if $I$ denotes the 
set of indices $i\le N$ for which $x_i=y_i$, we have
$$\sum_{i\not\in I} x_i\le\Vert x\Vert_2 f_c(A(a),x)\,.$$
By Lemma 6.2 the items $(x_i)_{i\not\in I}$ can be packed using at 
most $2\Vert x\Vert_2 f_c(A(a),x)+1$ bins.  The items $(x_i)_{i\in I}$ 
are exactly the items $(y_i)_{i\in I}$, so they can certainly 
be packed using at most $a$ bins, since $y\in B(a)$.  The result 
follows.\rbx
\enddemo

We provide $[0,1]$ with the measure $\mu$, and we denote by $P$ the 
product probability on $\Omega^N$.  The term $\Vert x\Vert_2$ of 
(6.2) will be disposed of by the following simple observation.

\proclaim{Lemma 6.3}  We have
$$\leqalignno{P(\Vert x\Vert_2\ge 2\sqrt{N}(EX^2_1)^{1/2}) &\le \exp 
(-2NEX^2_1)\,. &(6.3)\cr}$$
\endproclaim

\demo{Proof}  Since $e^x\le 1+2x$ for $x\le 1$, we have
$$E\exp X^{2}_i\le 1+2EX^2_i\le \exp 2EX^2_1$$
so that
$$E\exp\left(\sum_{i\le N}X^2_i\right)\le\exp 2NEX^2_1$$
for which (6.3) follows by Chebyshev inequality.\rbx
\enddemo

We can now prove the basic inequality.

\proclaim{Proposition 6.4}  We have, for all $t>0$, all $a>0$ that
$$\leqalignno{P(B_N(x)\le a)P(B_N(x)\ge a+2t\sqrt{N}(EX^2_1)^{1/2}+1) 
&\le e^{-t^2/4}+e^{-2NEX^2_1}\,. &(6.4)\cr}$$
\endproclaim

\demo{Proof}  Indeed, by (6.2), if $B_N(x)\ge a+2t\sqrt{N}(EX^2_1
)^{1/2}+1$, we have either $f_c(A(a),x)\ge t$ or $\Vert x\Vert_2\ge 
2\sqrt{N}(EX^2_1)^{1/2}$.  The result then follows from (4.1.2) and 
(6.2).
\enddemo

\proclaim{Theorem 6.5}  Denote by $M$ a median of $B_N(x)$.  Then for 
all $u\le 4\sqrt{2}N EX^2_1$ we have
$$P(\vert B_N(X_1,\dots ,X_N)-M\vert\ge 1+u)\le 8\exp \left( -{u^2 
\over 16NEX^2_1}\right)\,.$$
\endproclaim

\demo{Proof}  First, we take $a=M$ to obtain from (6.4) setting
$u=2t\sqrt{N} (EX^2_1)^{1/2}$, since 
$P(B_N\le M)\ge 1/2$,
$$\eqalign{P(B_N\ge M+u+1) &\le 2(e^{-t^2/4}+e^{-2NEX^2_1})\cr
&\le 4e^{-t^2/4}\,.\cr}$$
The bound for $P(B_N\le M-u-1)$ follows similarly taking $a=M-u-1$.
\rbx
\enddemo

\remark{Remarks}  1)  One can also get bounds for larger values of 
$u$, by adapting Lemma 6.3.

2)  It is instructive to find an alternate proof of Theorem 6.5 
using Corollary 2.2.4 rather than Theorem 4.1.1.
\endremark

\vfil\eject

\noindent{\bf 7.~~Subsequences}

{\bf 7.1.~~The longest increasing subsequence}

Consider points $x_1,\dots ,x_N$ of $[0,1]$.  We denote by $L_N(x_1,
\dots ,x_N)$ the length of the longest increasing subsequence of $x_1,
\dots ,x_N$.  That is, the largest integer $p$ such that we can find 
$i_1<\cdots <i_p$ for which $x_{i_1}\le\cdots\le x_{i_p}$.  It is 
simple to see that when $X_1,\dots ,X_N$ are independent uniformly 
distributed over $[0,1]$ (or, actually, distributed according to any 
non atomic probability), the r.v. $L_N(X_1,\dots ,X_N)$ is distributed 
like the longest increasing subsequence of a random permutation 
$\sigma$ of $\{ 1,\dots ,N\}$ (where the symmetric group $S_N$ is of 
course provided with the uniform probability).  The concentration of 
$L_N(X_1,\dots ,X_N)$ around its mean has been studied in particular 
in [F] and [B-B].  Sharper results will be obtained here as a simple 
consequence of Theorem 4.1.1.  We consider $\Omega =[0,1]^N$.  For 
$x=(x_i)_{i\le N}$ in $\Omega$, we set $L_N(x)=L_N(x_1,\dots ,x_N)$.  
For $a>0$, we set
$$A(a)=\{ x\in\Omega\,; L_N(x)\le a\}\,.$$
The basic observation is as follows.

\proclaim{Lemma 7.1.1}  For all $x\in\Omega^N$, we have
$$\leqalignno{a &\ge L_N(x)-f_c(A(a),x)\sqrt{L_N(x)}\,. &(7.1.1)\cr}$$
In particular,
$$\leqalignno{L_N(x) &\ge a+v\Rightarrow f_c(A(a),x)\ge {v \over 
\sqrt{a+v}}\,. &(7.1.2)\cr}$$
\endproclaim

\demo{Proof}  For simplicity, we write $b=L_N(x)$.  By definition, 
we can find a subset $I$ of $\{1,\dots ,N\}$ of cardinality $b$ such 
that if $i,j\in I$, $i<j$, then $x_i<x_j$.  By Lemma 4.1.2 (taking 
$\alpha_i=1$ if $i\in I$ and $\alpha_i=0$ otherwise), there exists 
$y\in A(a)$ such that $\card J\le f_c(A(a),x)\sqrt{b}$, where $J=\{i
\in I\,;\, y_i\not= x_i\}$.  Thus $(x_i)_{i\in I\backslash J}$ is an 
increasing subsequence of $y$; since $y\in A(a)$, we have $\card(I
\backslash J)\le a$, which proves (7.1.1).

To prove (7.1.2), we observe that by (7.1.1) we have
$$f_c(A(x),x)\ge {L_N(x)-a \over \sqrt{L_N(x)}}$$
and that the function $u\to (u-a)/\sqrt{u}$ increases for $u\ge a$.
\rbx
\enddemo

We denote by $M$ ($=M_N$) a median of $L_N$.

\proclaim{Theorem 7.1.2}  For all $u>0$ we have
$$\leqalignno{P(L_N\ge M+u) &\le 2\exp -{u^2 \over 4(M+u)} &(7.1.3)\cr
P(L_N\le M-u) &\le 2\exp -{u^2 \over 4M}\,. &(7.1.4)\cr}$$
\endproclaim

\demo{Proof}  To prove (7.1.3), we combine (7.1.2) with $M=a$ and 
(4.1.2).  To prove (7.1.4), we use (7.1.2) with $a=M-u$, $v=u$ to see 
that
$$L_N(x)\ge M\Rightarrow f_c(A(M-u),x)\ge {u \over \sqrt{M}}$$
so that
$$\leqalignno{P\left( f_c(A(M-u),x) \ge {u \over \sqrt{M}}\right) 
&\ge {1 \over 2}\,. &(7.1.5)\cr}$$
On the other hand, by (4.4.2),
$$\leqalignno{P\left( f_c(A(M-u),x) \ge {u \over \sqrt {M}}\right) 
&\le {1 \over P(A(M-u))} e^{-{u^2 \over 4M}}\,. &(7.1.6)\cr}$$
Comparing (7.1.5), (7.1.6) gives the required bound on $P(A(M-u))$.
\rbx
\enddemo

It seems worthwhile to state an abstract version of Theorem 7.1.2.  
Let us say that a function $L_N\colon\Omega^N\to\Bbb N$ is a {\it 
configuration function} provided it has the following property.

\noindent{\bf (7.1.7)}  Given any $x=(x_i)_{i\le N}$ in $\Omega^N$, 
there exists a subset $J$ of $\{ 1,\dots ,N\}$ with $\card J=L_N(x)$ 
such that, for each $y$ in $\Omega^N$, we have $L_N(y)\ge 
\card\{ i\in J\,;\, y_i=x_i\}$.

The reason for this name is that, intuitively, $L_N$ counts the size 
of the largest ``configuration'' formed by the points $x_i$.

The proof of the following is identical to that of Theorem 7.1.2.

\proclaim{Theorem 7.1.3}  If $L_N$ is a configuration function, then 
(7.1.3) and (7.1.4) hold.
\endproclaim

{\bf 7.2.~~Longest common subsequence}

Consider two sequences $x=(x_1,\dots ,x_N)$, $y=(y_1,\dots ,y_N)$ of 
numbers.  We define the length $L_{N,M}(x;y)$ of the longest common 
subsequence of $x$, $y$ as the largest integer $p$ for which 
there exists $1\le i_1<\cdots <i_p\le N$ and $1\le j_1<\cdots <j_p\le 
N'$ such that $x_{i_\ell}=y_{j_\ell}$ for each $\ell\le p$.  One 
interpretation of this is when $x_1,\dots x_N$ are chosen among a 
(small) finite number of possibilities (the letters of an alphabet)
$L_{N,N'}(x;y)$ is then the length of the longest ``subword'' of the 
words $x$, $y$ (and $N+N'-L_{N,N'}(x;y)$ is the so-called ``edit 
distance'' of the two words).  These considerations arise in a number 
of situations, such as genetics, speech recognition, etc.  Consider 
now a r.v. $X$, and two independent sequences $(X_i)_{i\le N}$, 
$(Y_j)_{j\le N'}$ independently distributed like $X$.  We are 
interested in the random variable $L_{N,N'} = L_{N,N'}(X_1,\dots ,X_N
\,;\,Y_1,\dots ,Y_{N'})$.

\proclaim{Theorem 7.2.1}  Consider a median $M$ ($=M_{N,N'}$) of 
$L_{N,N'}$.  Then for all $u>0$, we have
$$\leqalignno{P(L_{N,N'}\ge M+u) &\le 2\exp\left( -{u^2 \over 32(M+
u)}\right) &(7.2.1)\cr
P(L_{N,N'}\le M-u) &\le 2\exp\left( -{u^2 \over 32M}\right)\,.
&(7.2.2)\cr}$$
\endproclaim

\remark{Comments}  It is known that $\lim\limits_{N\to\infty}E(L_{N,
N})/N$ exists.  However, this limit can be very small, in the case 
where $X$ takes many possible values.  In this case, we have $M
\ll N$, and (7.2.1), (7.2.2) give a better result than Azuma's 
inequality.
\endremark

\demo{Proof}  The proof is very similar to the proof of Theorem 
7.1.2.  Consider $\Omega =[0,1]$, and for $x\in\Omega^{N+N'}$, 
consider
$$L(x)=L_{N,N'}(x_1,\dots ,x_N\,;\,x_{N+1},\dots ,x_{N+N'})\,.$$
Consider the set
$$A(a)=\{ x\,;\,L(x)\le a\}\,.$$
The basic inequality is that
$$\leqalignno{a &\ge L(x)-2\sqrt{2}f_c(A(a),x)\sqrt{L(x)}\,. 
&(7.2.3)\cr}$$

To see this, we set $b=L(x)$; we can find indices
$$1\le i_1<\cdots <i_b\le N<i_{b+1}<\cdots <i_{2b}\le N+N'$$
such that $x_{i_k}=x_{i_{k+b}}$ for $1\le k\le b$.  Consider the 
set $I=\{ i_k\,;\,1\le k\le 2b\}$.  
By Lemma 4.4.2, we can find $y\in A(a)$ such that
$$\leqalignno{\card\{ i\in I\,;\,x_i\not= y_i\} &\le f_c(A(a),x)
\sqrt{2b}\,. &(7.2.4)\cr}$$
Consider then
$$J=\{ k\le b\,;\,x_{i_k}=y_{i_k}\,;\,x_{i_{k+b}}=y_{i_{k+b}}\}\,.$$
By (7.2.4) we see that
$$\card J\ge b-2\card\{ i\in I\,;\, x_i\not= y_i\}\ge b-2\sqrt{2b} 
f_c(A(a),x)\,.$$
On the other hand, $L(y)\ge \card J$ since, for $k\in J$, we have 
$y_{i_k}=y_{i_{k+b}}$.  Also, since $y\in A(a)$, we have $L(y)\le a$.  
Condition (7.2.3) follows.  The rest of the proof is identical to 
that of Theorem 7.2.2.\rbx
\enddemo

\remark{Remark}  It is also possible to find a more general version 
of Theorem 7.1.3 that would contain Theorem 7.1.2.
\endremark

\vfil\eject

\noindent{\bf 8.~~Infimum and Percolation}

Consider an independent sequence $(X_i)_{i\le N}$ of positive r.v.  
Consider a family ${\Cal F}$ of $N$-tuple $\alpha =(\alpha_i)_{i\le 
N}$ of positive numbers.  Our prime topic of interest in the 
present section is the random variable
$$\leqalignno{Z' &= Z'_{\Cal F}=\inf_{\alpha\in{\Cal F}}\sum_{i\le
N}\alpha_iX_i\,. &(8.1)\cr}$$
It does matter a lot that we take an infimum rather than a supremum.  
The function of the $'$ in $Z'$ is to indicate that we take such an 
infimum.  Rather that (8.1) one can also write
$$Z'=-\sup_{\alpha\in{\Cal F}}\sum_{i\le N}(-\alpha_i)X_i$$
but the numbers $-\alpha_i$ are negative.  In Section 13, we will 
have to study the r.v.
$$\leqalignno{Z &= \sup_{\alpha\in{\Cal F}}\sum_{i\le N}\alpha_i 
X_i &(8.2)\cr}$$
where $\alpha_i$ and $X_i$ can possibly have any signs.  In order to 
avoid repetition, we will study the variables $Z$ given by (8.2).

{\bf 8.1.~~The basic result}

Consider a family ${\Cal F}$ of $n$-tuples $\alpha =(\alpha_i)_{i\le 
N}$.  We make no assumption on the sign of $\alpha_i$.  We set 
$\sigma =\sup\limits_{\alpha\in{\Cal F}}\Vert\alpha\Vert_2$, where
$\Vert\alpha\Vert_2=\left(\sum\limits_{i\le N}\alpha^2_i
\right)^{1/2}$.  We consider independent r.v. $X_i$, and we assume 
that for each $i$ there is number $r_i$ such that $r_i\le X_i\le r_i
+1$.

\proclaim{Theorem 8.1.1}  Consider the r.v. $Z$ given by (8.2), and 
a median $M$ of $Z$.  Then, for all $u>0$, we have
$$\leqalignno{P(\vert Z-M\vert\ge u) &\le 4\exp\left( -{u^2 \over 
4\sigma^2}\right)\,. &(8.1.1)\cr}$$
\endproclaim

\demo{Proof}  This will again follow from Theorem 4.1.1.

{\bf Step 1.}  Set $\Omega =[0,1]$, and for $x=(x_i)_{i\le N}\in
\Omega^N$, set
$$Z(x)=\sup\limits_{\alpha\in{\Cal F}}\sum_{i\le N}\alpha_i(r_i+
x_i)\,.$$
Consider $a\in\Bbb R$, and $A(a)=\{ y\in\Omega^N\,;\, Z(y)\le a\}$.  
The basic observation is that
$$\leqalignno{\forall x\in\Omega^N\,,~~Z(x) &\le a+\sigma f_c(A(a),
x)\,. &(8.1.2)\cr}$$
To prove this, consider $\alpha\in{\Cal F}$.  By Lemma 4.1.2, we can 
find $y\in A(x)$ such that, if $I=\{ i\le N\,;\, y_i\not= x_i\}$, 
then we have
$$\leqalignno{\sum_{i\in I}\vert\alpha_i\vert &\le\Vert\alpha\Vert_2 
f_c(A(a),x) \le \sigma f_c (A(a),x)\,. &(8.1.3)\cr}$$
We then have
$$\left|\sum_{i\le N}\alpha_i(r_i+y_i)-\sum_{i\le N}\alpha_i(r_i+x_i)
\right| \le\sum_{i\in I}\vert
\alpha_i\vert~\vert y_i-x_i\vert \le \sum_{i\in I}\vert\alpha_i
\vert\,.$$
Thus, by (8.1.3)
$$\sum_{i\le N}\alpha_i(r_i+x_i)\le Z(y)+\sigma f_c(A(a),x)\le a+
\sigma f_c(A(a),x)\,,$$
and taking the supremum over $\alpha$ proves (8.1.2).

{\bf Step 2.}  We provide the $i^{\text{th}}$ factor $[0,1]$ with 
the law $\mu_i$ of $X_i-r_i$.  We denote by $P$ the product 
probability.  Thus by (8.1.2) and (4.1.2)
$$P(Z(x)\ge b)\le P\left( f_c(A(a),x)\ge {b-a \over \sigma}\right) 
\le{1 \over P(A(a))}\exp\left(-{(b-a)^2 \over 4\sigma^2}\right)\,,$$
i.e.,
$$P(Z(x)\ge b)P(Z(x)\le a)\le\exp\left( -{(b-a)^2 \over 4\sigma^2}
\right)$$
from which (8.1.1) follows as in Chapter 7, since law of $Z(x)$ under 
$P$ coincides with the law of $Z$.\rbx
\enddemo

{\bf 8.2.~~General moments}

In the present section we rely on the theory of Section 4.4.  We 
start with some preliminaries.  Consider a convex function $\psi$ on 
$\Bbb R^+$ that satisfies (4.4.2) and $\psi (0)=0$.  Consider a 
family ${\Cal F}$ of $N$-tuples as in Section 8.1.  For $u>0$, we 
define
$$\psi_{\Cal F}(u) =\inf \left\{\sum_{i\le N}\psi (s_i)\,;\,\exists
\alpha\in{\Cal F}\,,\, \sum_{i\le N}s_i\vert\alpha_i\vert\ge u\right
\}\,.$$

The simplest case is when $\psi (x)=x^2$.  In that case it is easily 
seen that $\psi_{\Cal F}(u)=u^2/\sigma^2$, where $\sigma^2=\sup\{\Vert
\alpha\Vert^2_2\,;\,\alpha\in{\Cal F}\}$.  The most interesting case 
is arguably the case where $\psi =\psi_0$ is given by
$$\psi_0(x)=x^2~~\text{if}~~x\le 1\,;\qquad \psi_0(x)=2x-1~~\text{if}~~
x\ge 1\,.$$
If we set
$$\tau = \sup\{\vert\alpha_i\vert\,;\, i\le N\,;\,\alpha\in{\Cal F}
\}\,,$$
we note that, for given $\alpha\in{\Cal F}$, for each $s=(s_i)_{i\le 
N}$, setting $J=\{ i\le N\,;\,s_i\le 1\}$ we have
$$\eqalign{\sum_{i\le N} s_i\vert\alpha_i\vert &=\sum_{i\in J}s_i
\vert\alpha_i\vert +\sum_{i\not\in J}s_i\vert\alpha_i\vert\cr
&\le \sigma \left(\sum_{i\in J}s^2_i\right)^{1/2}+\tau\sum_{i\not\in 
J}s_i\cr
&\le \sigma\left(\sum_{i\le N}\psi (s_i)\right)^{1/2}+\tau\sum_{i\le 
N}\psi (s_i)\,.\cr}$$
Thus, if $\sum_{i\le N}s_i\vert\alpha_i\vert\ge u$, then either 
$\sum_{i\le N}\psi (s_i)\ge u^2/4\sigma^2$, or else $\sum\limits_{i
\le N}\psi (s_i)\ge\tau /2$, and thus
$$\leqalignno{\Psi_{\Cal F}(u) &\ge \min\left({u^2 \over 4\sigma^2}
\,,\,{u \over 2\tau}\right)\,. &(8.2.1)\cr}$$

The basic observation is as follows.

\proclaim{Proposition 8.2.1}  Consider ${\Cal F}$, $\psi$ as above.  
Set $\Omega =\Bbb R$, and consider the function $Z(x)=
\sup\limits_{\alpha\in{\Cal F}}\sum\limits_{i\le N}\alpha_ix_i$.  
Consider $a\in\Bbb R$, and $A(a)=\{ y\,;\,Z(y)\le a\}$.  Then
$$\leqalignno{\forall x\in\Omega^N\,,\quad f_{h,\psi}(A(a),x) &\ge 
\Psi_{\Cal F}(Z(x)-a) &(8.2.2)\cr}$$
when the function $h$ is defined on $\Bbb R\times\Bbb R$ by
$$\leqalignno{h(\omega ,\omega ') &=\vert\omega -\omega '\vert\,. 
&(8.2.3)\cr}$$
Moreover, when $\alpha_i\le 0$ for each $i\le N$ and each $\alpha
\in{\Cal F}$, we can take
$$\leqalignno{h(\omega ,\omega ') &= (\omega '-\omega )^+\,. 
&(8.2.4)\cr}$$
\endproclaim

\demo{Proof}  By definition of $f_{h,\psi}$, given $\varepsilon >0$, 
we can find $s\in V_{A(a)}(x)$ such that
$$\sum_{i\le N}\psi (s_i)\le f_{h,\psi }(A(a),x)+\varepsilon\,.$$
Consider $\alpha =(\alpha_i)\in{\Cal F}$.  Then there exists $s'\in 
U_{A(a)}(x)$ such that $\sum\limits_{i\le N}\vert\alpha_i\vert s'_i
\le\sum\limits_{i\le N}\vert\alpha_i\vert s_i$.  This means that there 
is $y\in A(a)$ for which $\sum\limits_{i\in I}\vert\alpha_i\vert h(
x_i,y_i)\le\sum\limits_{i\le N}\vert\alpha_i\vert s_i$, where $I=\{ 
i\le N\,;\,x_i\not= y_i\}$.  Now
$$\leqalignno{\sum_{i\le N}\alpha_i x_i &= \sum_{i\le N}\alpha_iy_i+
\sum_{i\in I}\alpha_i(x_i-y_i)\,. &(8.2.5)\cr}$$
We have $\alpha_i(x_i-y_i)\le \vert\alpha_i\vert~\vert x_i-y_i\vert$.  
If $\alpha_i$ is $\le 0$ we have $\alpha_i(x_i-y_i)\le\vert\alpha_i
\vert (y_i-x_i)^+$.  Thus in all cases under consideration, we have
$$\sum_{i\in I}\alpha_i(x_i-y_i)\le\sum_{i\in I}\vert\alpha_i\vert 
h(x_i,y_i)\le\sum_{i\le N}\vert\alpha_i\vert s_i$$
so that, by (8.2.5)
$$\eqalign{\sum_{i\le N}\alpha_ix_i &\le \sum_{i\le N}\alpha_iy_i+
\sum_{i\le N}\vert\alpha_i\vert s_i\cr
&\le a+\sum_{i\le N}\vert\alpha_i\vert s_i\,.\cr}$$
Taking the $\sup$ over $\alpha$ yields
$$\sup_{\alpha\in{\Cal F}}\sum_{i\le N}\vert\alpha_i\vert s_i\ge Z(x) 
-a$$
and the result follows by definition of $\Psi_{\Cal F}$.\rbx
\enddemo

\proclaim{Corollary 8.2.2}  Consider a family ${\Cal F}$ of $N$-tuples 
$\alpha =(\alpha_i)_{i\le N}$.  Consider a sequence of independent 
r.v. $(X_i)_{i\le N}$ with common law $\mu$.  Assume that (4.4.6) 
holds (for a certain function $\theta$) when $P=\mu^{\otimes N}$, and 
where $h$ is the function determined in Proposition 8.2.1.  Then the 
r.v. $Z=\sup\limits_{\alpha\in{\Cal F}} \sum\limits_{i\le N}\alpha_i
X_i$ satisfies
$$\leqalignno{u\ge 0 &\Rightarrow P(Z\ge M+u)\le\exp \left(\theta
\left({1 \over 2}\right)-{1 \over K}\Psi_{\Cal F}(u)\right) &(8.2.6)\cr
u\ge 0 &\Rightarrow P(Z\le M-u) \le\xi\left({1 \over K}\Psi_{\Cal F}
(u)-\log 2\right) &(8.2.7)\cr}$$
where $M$ is a median of $Z$.
\endproclaim

\demo{Proof}  Using (4.4.6) and Chebyshev inequality, we have 
$$P(f_{h,\psi}(A(a),x)\ge t)\le\exp\left(\theta ( P(A(x))) -{t \over 
K}\right)$$
where $A(a)$ is the set of Proposition 8.2.1; thus, by (8.2.2) 
setting $t=\Psi_{\Cal F}(b-a)$, for $Z(x)\ge b$ we have $f_{h,\psi 
}(A(a),x)\ge t$, so that
$$P(Z\ge b)\le\exp\left(\theta (P(Z\le a)) - {1 \over K}\Psi_{\Cal F}
(b-a)\right)\,.$$

Taking $a=M$, $b=M+u$ imply (8.2.6).  Taking $b=M$, $a=M-u$ imply
$${1 \over 2} \le \exp\left(\theta (P(Z\le M-u))-{1 \over K}
\Psi_{\Cal F}(u)\right)$$
from which (8.2.7) follows.\rbx
\enddemo

We now go back to our main line of study, that of the r.v. $Z'=
\sup\limits_{\alpha\in{\Cal F}}~\sum\limits_{i\le N}(-\alpha_i)X_i$.  
In order to apply Corollary 8.2.2, we need (4.4.6) for the 
penalty function $h(x,y)=(y-x)^+$.  Since $X_i$ is positive, its law 
$\mu$ is supported by $\Bbb R^+$.  Thereby, only the properties of 
$h$ on $\Bbb R^+\times\Bbb R^+$ matter; but then $(y-x)^+\le y$.  
Thus, to have (4.4.6) it suffices that the function $h(x,y)=y$ 
satisfies the conditions of Theorem 4.4.1.  The case where the 
function $h(x,y)$ depends on $y$ only has been discussed after Theorem 
4.4.1.  Thus, we have proved the following.

\proclaim{Theorem 8.2.3}  Consider a family ${\Cal F}$ on $N$-tuples 
of positive numbers, and independent identically distributed 
nonnegative r.v. variables $(X_i)_{i\le N}$.  Consider functions 
$\theta$, $\xi$, $w$ as in Theorem 4.4.1.  Assume that (2.6.1), 
(4.4.2), (4.4.3) hold.  Assume that condition $H(\xi ,w)$ holds.  
Assume that the median $m$ of $X_1$ is $\le 1$, and that for $t\ge 
m$, we have
$$\leqalignno{w(P(X_1\ge t)) &\ge \psi (t)\,. &(8.2.8)\cr}$$

Then if $M$ is a median of $Z'=\inf\limits_{\Cal F}~\sum\limits_{i\le 
N}\alpha_iX_i$, the following holds (where the constant $K$ depends 
only on the parameter of $\gamma$ of Theorem 4.4.1)
$$\leqalignno{u\ge 0 &\Rightarrow P(Z'\le M-u)\le \exp\left(\theta
\left({1 \over 2}\right) - {1 \over K}\Psi_{\Cal F}(u)\right) 
&(8.2.9)\cr
u\ge 0 &\Rightarrow P(Z'\ge M+u)\le \xi\left({1 \over K}\Psi_{\Cal F}
(u)-\log 2\right)\,. &(8.2.10)\cr}$$
\endproclaim

\remark{Comment}  A striking feature of this result is the different 
forms of (8.2.9) and (8.2.10).  This phenomenon is well known in the 
case where ${\Cal F}$ consists of single point $\alpha$.  In 
that case, ${\Cal F}$ is a sum of positive independent r.v. $Y_i$.  
The lower tails of $Z$ have a tendency to be ``subgaussian'' ([H]) 
while the upper tails of $Z$ certainly depend much on the upper 
tails of the variables $Y_i$.
\endremark

\proclaim{Corollary 8.2.4}  There exists a universal constant $U$ 
with the following property.  Assume that $\psi$ satisfies (4.4.2).  
Assume that
$$\leqalignno{\forall t\ge 1\,,\qquad P(X_1\ge t) &\le \exp (-2\psi 
(t))\,. &(8.2.11)\cr}$$
Then we have
$$\leqalignno{u\ge 0 &\Rightarrow P(\vert Z-M\vert\ge u) \le 3\exp 
\left( -{1 \over K}\Psi_{\Cal F}(u)\right)\,. &(8.2.12)\cr}$$
\endproclaim

\demo{Proof}  We take $\xi (x)=e^{-x}$, $\theta (x)=-\log x$.  
According to Proposition 2.6.1, condition $H(\xi ,w)$ holds if $\int 
e^w d\lambda\le 2$, so in particular, if $w(t)=-{1 \over 2}\log 
t$.  Also, by (4.4.2), $\psi (1)=1$, so that (8.2.11) implies that 
the median of $X_1\le 1$.  Thus Corollary 8.2.4 follows from Theorem 
8.2.3.
\enddemo

\proclaim{Corollary 8.2.5}  Assume that (2.6.5) holds for a certain 
number $L$.  Then, for some constant $K$ depending on $\xi$ only, if 
for all $t\ge 1$ we have 
$$\leqalignno{P(X_1\ge t) &\le {1 \over K}\vert\xi '(\psi (t))\vert 
&(8.2.13)\cr}$$
then (8.2.9), (8.2.10) hold (for a constant $K$  depending on $\xi$ 
only).
\endproclaim

\demo{Proof}  We simply have to find a function $w$ that satisfies 
(8.2.8) and such that condition $H(\xi ,w)$ holds.  It follows from 
Proposition 2.6.3 that if we take $R$ large enough ($R$ can 
actually be taken depending on $L$ and $\int^\infty_0 \xi\,d\lambda$ 
only) then the function $w$ such that
$$\forall b\ge c\qquad \vert\{ w\ge b\}\vert ={1 \over R}\vert\xi '
(b)\vert$$
satisfies condition $H(\xi ,w)$.  Now, if we take $K=R$ in (8.2.13), 
then
$$\vert\{ w\ge \psi (t)\}\vert ={1 \over R}\vert\xi '(\psi (t))\vert 
\ge P(X_1\le t)\,,$$
so that $w(P(X_1\le t))\ge \psi (t)$ since $w$ is non-decreasing.\rbx
\enddemo

We now explain why Corollary 8.2.5 is sharp.  Consider the case where 
${\Cal F}$ consists of the single element $\alpha =(\alpha_i)$, where 
$\alpha_i=1/\sqrt{N}$.  Consider $\psi$ such that $\psi (x)=x^2$ if 
$x\le 1$ and $\psi (x)=2x-1$ for $x\ge 1$.  Then, for $u=\sqrt{N}$, 
$\Psi_{\Cal F} (u)\ge N/4$ by (8.2.1).  Consider a r.v. $X_i\in\{ 0,
N\}$, with
$$P(X_i=N)=p=:{1 \over K}\vert\xi '(2N-1)\vert\,.$$
Under condition (2.6.5), we have $\lim\limits_{x\to\infty} x\xi '(x)
=0$, and it is not a restriction to assume $Np\le 1/2$.  Thus the 
median of $Z=N^{-1/2}\sum_{i\le N}X_i$ is zero.

Now
$$P(Z\ge u) = P\left({1 \over \sqrt{N}}\sum_{i\le N}X_i\ge u\right) 
\sim Np={N \over R}\xi '(2N-1)$$
and the bound $\xi\left({N \over K}\right)$ of (8.2.10) is indeed 
reasonably good, as $x\xi '(x)$ is of order $\xi (x)$ for many choices 
of $\xi$.

{\bf 8.3.~~First time passage in percolation}

Consider a graph $(V,E)$ where $V$ is the set of vertices, $E$ the 
set of edges.  Assume that we have a family $(X_e)_{e\in E}$ of 
positive r.v. distributed liked a given r.v. $X$.  ($X_e$ represents 
the passage time through edge $e$.)  Consider a family ${\Cal S}$ of 
sets of edges; and for $S\in {\Cal S}$, consider $X_S = 
\sum\limits_{e\in S}X_e$.  In the case where $S$ is a path, i.e., 
consists of the edges $e_{v_1v_2}, e_{v_2v_3}, \dots ,e_{v_{k-1},
v_k}$ linking vertices $v_1,\dots ,v_k$, $X_S$ represents the 
``passage time through $S$''.  Let us set $Z'_{\Cal S}=\inf\limits_{S
\in{\Cal S}}X_S$.  Let us set $r=\sup\limits_{S\in{\Cal S}}\card S$.  
Denote $M$ a median of $Z_{\Cal S}$.  The following is a consequence 
of (8.2.1) and Corollary (8.2.4).

\proclaim{Proposition 8.3}  There exists a universal constant $K$ such 
that if $E\exp K^{-1}X\le 2$, we have
$$\leqalignno{\forall u>0\,,\qquad P(\vert Z_{\Cal S}-M\vert\ge u) 
&\le 4\exp\left( - {1 \over K}\min 
\left({u^2 \over r} ,u\right)\right)\,. &(8.3.1)\cr}$$
\endproclaim

Consider the case where $V=\Bbb Z^2$, $E$ consists of the edges that 
link any two adjacent vertices.  Denote by ${\Cal S}$ the sets of 
self-avoiding paths linking the origin to the point $(0,n)$; and by 
${\Cal S}(C)$ the subset of ${\Cal S}$ consisting of paths of length 
$\le Cn$.  H. Kesten [K1] proved that if $P(X=0)<{1 \over 2}$, then, 
for some constant $C$ independent of $n$, we have
$$P(Z'_{\Cal S} = Z'_{{\Cal S}(C)})\ge 1-Ce^{-n/C}\,.$$
It then follows from (8.3.1) that for some constant $C'$ independent 
of $n$, we have
$$\leqalignno{u\le {n \over C'} &\Rightarrow P(\vert Z'_{\Cal S}-M
\vert\ge u) \le 5\exp\left( -{u^2 \over C'n}\right)\,. &(8.3.2)\cr}$$

This improves recent results of H. Kesten [K2], based on the use of 
martingales, who proves (8.3.2) with an exponent $u/C'\sqrt{n}$.  It 
should, however, be pointed out that the reason why martingales allow 
some success on this problem is because we consider only sums of the 
type $\sum\alpha_eX_e$ for very special families $\alpha =(\alpha_e)$.  
Martingales are apparently powerless to approach Corollary 8.2.5.

It is pointed out in the literature that (in the case $V=\Bbb Z^2$) 
(8.3.2) apparently does not give the correct rate.  In view of 
Corollary 8.2.5, the obvious approach to improve (8.3.2) would be to 
show that $Z_{\Cal S}$ is very close to $Z_{\Cal F}$, where the 
family ${\Cal F}$ of sequences $(\alpha_e)_{e\in E}$ satisfies 
$\sigma =\sup\limits_{\alpha\in{\Cal F}}\Vert\alpha\Vert_2\ll n$.  
There is an obvious candidate for ${\Cal F}$.  Indeed, consider the 
family ${\Cal F}'$, defined as follows:  ${\Cal F}'$, seen as a 
subset of $(\Bbb R^+)^E$, is the convex hull of the family of points 
$a_S$ given by $a_S(e)=1$ if $e\in S$ and $a_S(e)=0$ if $e\not\in S$, 
for all $S\in{\Cal S}$.  Then, obviously, $Z'_{\Cal S}=Z'_{{\Cal F}'}$.  
Then consider the family ${\Cal F}(\sigma )$ of sequences 
$(\alpha_e)_{e\in V}$ of ${\Cal F}'$ for which $\Vert\alpha\Vert_2\le
\sigma$.  Then $Z'_{\Cal S}\le Z'_{{\Cal F}(\sigma )}$.  Thus if one 
could show that for some $\sigma =o(n)$, and still have 
$Z'_{{\Cal F}(\sigma )}\le Z'_{\Cal S}+o(\sqrt{n})$, with probability 
$1-o(n^{-1})$, one would obtain that the likely fluctuations of 
$Z'_{\Cal S}$ from $M$ are $o(\sqrt{n})$.  Roughly speaking, this 
means that the shortest passage time from $(0,0)$ to $(0,n)$ is 
(within $o(\sqrt{n})$) obtained through a number of rather disjoint 
paths.  Proving such a statement is apparently a long range 
program in Percolation theory.
\vfil\eject

\noindent{\bf 9.~~Chromatic Number of Random Graphs}

The use of martingales has allowed several important progresses in 
the understanding of the chromatic number of random graphs.  Use of 
martingales does require ingenuity.  This chapter will demonstrate 
that Theorem 4.1.1 achieves somewhat better results than martingales 
in a completely straightforward 
manner.

For simplicity we call a graph $G$ with vertice set $V=\{1,\dots ,n\}$ 
a subset of $E_0=\{ (i,j)\,;\, i<j\}$.  If $(i,j)$ belongs to $G$, we 
say that $i,j$ are linked by an edge.

A subset $I$ of $V$ is called independent if no two points of $I$ are 
linked by an edge (the word independent there should not be confused 
with its probabilistic meaning).  The chromatic number 
$\chi (G,A)$ of a subset $A$ of $V$ is the smallest number of 
independent sets that can cover $A$; that is, the vertices of $A$ can 
be given $\chi (G,A)$ colors so that no two points with the same 
color are linked by an edge.  We set
$$\chi (G,m)=\inf\{\chi (G,A)\,;\,\card A=m\}\,.$$

Given $p$, $0<p<1$, the random graph $G=G(n,p)$ is defined by putting 
each possible edge $(i,j)$ in $G$ with probability $p$, independently 
of what is done for the other edges.

The chromatic number is remarkably concentrated, as the following 
shows.

\proclaim{Theorem 9.1}  Consider $k\in\Bbb N$ and $t>0$.  Then there 
exists an integer $a$ such that
$$\leqalignno{P &(\chi (G(n,p),m)\in [a-k,a])\cr
&\quad \ge 1-2e^{-t^2/8}-P(\sup\{\chi (G(n,p),F)\,;\,F\subset V
\,,\,\card F\le t\sqrt{m}\} >k)\,. &(9.1)\cr}$$
\endproclaim

\remark{Comments}  1)  The last term is always zero for $k> t
\sqrt{m}$.  But when $p=n^{-\alpha}$ ($\alpha >0$), it is still small 
for smaller values of $k$.  See [S-S], [A-S, p. 88].

2)  Another version of this Theorem could be proved, in the spirit of 
Theorem 7.1.3, concerning the concentration property of the number
$$\max\{\card F\,;\,\chi (G(n,p),F)\le m\}\,.$$

3)  With a bit of care, we can replace $m$ by $m-1$ in the right hand 
side of (9.1), and improve the coefficient $1/8$.
\endremark

\demo{Proof}  We set
$$b=P(\sup\{\chi (G(n,p),F)\,;\,\card F\le t\sqrt{m}\} >k)\,.$$
We then define $a$ as the largest integer for which
$$\leqalignno{P(\chi (G(n,p),m)\ge a) &\ge e^{-t^2/8}+b\,. &(9.2)\cr}$$
Thus
$$\leqalignno{P(\chi (G(n,p),m)>a) &< e^{-t^2/8}+b\,. &(9.3)\cr}$$

In order to apply Theorem 4.1.1, we must represent the underlying 
probability space as a product space.  The first idea that comes to 
mind would be to use $\{ 0,1\}^{E_0}$;  this is not a good choice.  
For $2\le j\le n$, set $\Omega_j=\{ 0,1\}^{j-1}$.  Set $\Omega '=
\prod\limits_{2\le j\le n} \Omega_j$.  We write $\omega\in\Omega '$ 
as $(\omega_j)_{j\le n}$, where $\omega_j =(\omega_{i,j})_{i
\le j-1}\in\Omega_j$.  To $\omega$ we associate the graph $G(\omega )$ 
such that, for $i<j$, $(i,j)\in G(\omega )$ if and only if 
$\omega_{i,j} =1$.  The only property of $G(n,p)$ we need is that it 
is distributed as $G(\omega )$ for a certain product measure $P$ on
$\prod\limits_{j\le n}\Omega_j$.

Define $A\subset\Omega '$ as the set of $\omega$ for which
$$\chi (G(\omega ),m)\ge a\,;\,\sup\{\chi (G(\omega ),F)\colon\card 
F\le t\sqrt{m}\}\le k\,.$$
Thus by (9.2) we have $P(A)\ge e^{-t^2/8}$.  
Combining Theorem 4.1.1 and Lemma 4.1.2, we see that $P(B)\ge 1-
e^{-t^2/8}$, where we have set
$$\leqalignno{B &=\{\omega\,;\,\forall (\alpha_j)_{2\le j\le n}\,,\,
\exists\omega '\in A\,;\,\sum\alpha_j1_{\{\omega_j\not=\omega_j'\}}
\le t\sqrt{\sum_j\alpha^2_j}\}\,. &(9.4)\cr}$$

To finish the proof, it suffices to show that
$$\leqalignno{\omega\in B &\Rightarrow \chi (G(\omega ),m)\ge a-k\,. 
&(9.5)\cr}$$

So, consider $\omega\in B$, and set $r=\chi (G(\omega ),m)$.  Consider 
a subset $F$ of $V$, of cardinal $m$, such that $\chi (G(\omega ),F)
=r$.  We use (9.3) with $\alpha_j=1$ if $j\in F$ and zero otherwise.  
Thus there is $\omega \in A$ such that if $J=\{j\in F\,;\,\omega_j
\not=\omega_j'\}$, then $\card J\le t\sqrt{m}$.  But obviously,
$$\chi (G(\omega '),F\backslash J) = \chi (G(\omega ),F\backslash I)
\le r$$
and thus
$$\eqalignno{a\le\chi (G(\omega '),F) &\le r+\chi (G(\omega '),I)\cr
&\le r+k\,. &\bx\cr}$$
\enddemo

In order to obtain an upper bound for $\chi_G$, the most obvious 
approach is the ``greedy'' one:  one chooses an independent set $W_1$ 
of maximal sizes, and remove its vertices and all edges adjacent.  
One is then left with a graph on fewer vertices, and one iterates the 
process until exhaustion.  To make this approach work one needs a 
competent bound on the probability that a random graph contains 
at least one independent set of size $r$.  Such bounds were first 
obtained by B. Bollobas [B], using martingales.  A recent powerful 
correlation inequality of Janson [J] is both simpler and more powerful 
than the martingale approach (compare [A-S] p. 87 and p. 148).  It is 
of some interest to note that Theorem 4.1.1 does as well as Janson's 
inequality.  We fix an integer $r$.  For $e=(i,j)\in E_0$, 
we denote by $N(G,e)$ the number of independent sets of size $r$ that 
contain $i,j$.

\proclaim{Proposition 9.2}  Consider a number $u$, and assume that
$$\leqalignno{P\left( u\sqrt{\sum_{e\in E_0}N(G(n,p),e)^2}\le\sum_{e
\in E_0}N(G(n,p),e)\right) &>{1 \over 2}\,. &(9.4)\cr}$$
Then
$$P(G(n,p)~\text{contains~no~independent~set~of~size~}r)\le 2\exp
\left( -{u^2 \over r^2(r-1)^2}\right)\,.$$
\endproclaim

\demo{Proof}  We set $\Omega =\{ 0,1\}$, provided with the probability 
that gives weight $p$ to $1$ (and $1-p$ to $0$).  Consider the product 
probability $P$ on $\Omega^{E_0}$.  For $x\in (x_e)_{e\in E_0}$ we 
define $G(x)$ by $e(i,j)\in G(x)$ if and only if $x_e =1$.  The graph 
$G(x)$ is distributed like $G(n,p)$.

Consider the set $A\subset\Omega^{E_0}$, given by
$$A = \{ y\,;\,G(y)~\text{contains~no~independent~set~of~size~}r\}\,.$$
Consider $t_0=2\sqrt{\log (2/P(A))}$.  If we combine (9.4), Theorem 
4.1.1 and Lemma 4.1.2, we see that there exists $x$ such that
$$\leqalignno{u\sqrt{\sum_{e\in E_0}N(G(x),e)^2} &\le \sum_{e\in E_0}
N(G(x),e) &(9.5)\cr}$$
with the property that
$$\forall (\alpha_e)_{e\in E_0}\,,\,\exists y\in A\,,\,\sum_{x_e
\not= y_e}\alpha_e\le t_0\sqrt{\sum_{e\in E_0}\alpha^2_e}\,.$$
In particular, there exists $y\in A$, such that if
$$C=\{ e\in E_0\,,\,x_e\not= y_e\}$$
we have
$$\leqalignno{\sum_{e\in C}N(G(x),e) &\le t_0\sqrt{\sum_{e\in E_0}N
(G(x),e)^2}\cr
&\le {t_0 \over u}\sum_{e\in E_0}N(G(x),e) &(9.6)\cr}$$
where the last inequality follows from (9.5).  The total number $N$ 
of independent sets of $G(x)$ of size $r$ is
$$\leqalignno{N &= \left({r(r-1) \over 2}\right)^{-1}\sum_{e\in E_0}
N(G(x),e)\,. &(9.7)\cr}$$

We must have
$$N\le\sum_{e\in C}N(G(x),e)$$
for otherwise there would be an independent set of size $r$ of $G(x)$ 
that would contain no edge of $C$, and thus would be an independent 
set of $G(y)$, which is impossible.  Combining with (9.6), (9.7), we 
get $t_0\ge {2u \over r(r-1)}$, so that
$$\eqalignno{P(A) &\le 2\exp -{u^2 \over r^2(r-1)^2}\,. &\bx\cr}$$
\enddemo

In order to take advantage of Proposition 9.2, one must find competent 
(= large) values of $u$ for which (9.4) holds.  For example, one can 
take $u=u_1/u_2$, where
$$\leqalignno{P\left(\sum_{e\in E_0}N((G,p),e)^2\le u^2_2\right) &> 
{3 \over 4} &(9.8)\cr
P\left(\sum_{e\in E_0} N((G,p),e)\ge u_1\right) &> {3 \over 4} 
&(9.9)\cr}$$
We then find values of $u_2$ (resp. $u_1$) using Chebyshev inequality 
(resp. the second moment method).  Not surprisingly that leads to 
unpleasant computations (as seems unavoidable in this topic).  
These are better not reproduced here, and left to the specialist that 
wants to evaluate the strength of Proposition 9.2.
\vfil\eject

\noindent{\bf 10.~~The Assignment Problem}

Consider a number $N$, and two disjoint sets $I$, $J$ of cardinal $N$.  
An assignment is a one to one map $\tau$ from $I$ to $J$.  Consider a 
matrix $a=(a_{i,j})_{i\in I,j\in J}$, such that $a_{i,j}$ represents 
the cost of assigning $j$ to $i$.  The cost of the assignment $\tau$ 
is $\sum\limits_{i\in I}a_{i,\tau (i)}$ and the problem is to find 
the assignment of minimal cost.

Assume now that the costs $a_{i,j}$ are taken equal to $X_{i,j}$, 
where the r.v. $(X_{i,j})_{i\in I,j\in J}$ are independent uniformly 
distributed over $[0,1]$.  Consider the r.v.
$$L_N=\inf\left\{\sum_{i\in I}X_{i,\tau (i)}\,;\,\tau~
\text{assignment}\right\}\,.$$
It is a remarkable fact [W] that $E(L_N)$ is bounded independently of 
$L_N$.  (Actually $E(L_N)\le 2$ [Ka].)

In this section we try to bound the fluctuations of $L_N$; the 
challenge is that the average value of $L_N$ is of the same order as 
the average value of the costs $X_{i,j}$, and that $N^2$ of these 
costs are involved.

We will first show that we can replace the costs $X_{i,j}$ by $Y_{i,j}
=\min (X_{i,j},v)$ for $v$ of order $N^{-1}(\log N)^2$; then we will 
appeal to Theorem 4.1.1.

A {\it digraph} $D$ will be a subset of $I\times J$.  (If $(i,j)\in 
D$, we think as $i,j$ being linked by an edge.)  The digraphs of use 
will mostly consist of those couples $(i,j)$ for which $X_{i,j}$ is 
small.  Consider a digraph $D$, and $S\subset I$.  Se set
$$D(S) = \{ j\in J\,;\,\exists i\in S\,,\, (i,j)\in D\}\,.$$
We will say that a digraph $D$ is $\alpha$-expanding ($\alpha\ge 2$) 
if the following occurs, for all subsets $S$ of $I$:
$$\leqalignno{\card S\le {N \over 2} &\Rightarrow \card D(S)\ge \min
\left(\alpha\card S,{N \over 2}\right) &(10.1)\cr
\card S\ge {N \over 2} &\Rightarrow \card D(S)\ge N-{1 \over \alpha} 
(N-\card S)\,. &(10.2)\cr}$$

Our first lemma mimics an argument of Steele and Karp [S-K].

\proclaim{Lemma 10.1}  Consider an $\alpha$-expanding digraph $D$ and 
an integer $m$ such that $\alpha^m \ge N/2$.  Consider a one to one 
map $\tau$ from $I$ to $J$.  Then, given any $i\in I$, we can find 
$n\le 2m$ and disjoint points $i_1=i,i_2,\dots ,i_{n+1}=i$ such that 
for $1\le\ell\le n$, we have $(i_\ell ,\tau (i_{\ell +1}))\in D$.
\endproclaim

\demo{Proof}  We fix $i\in I$.  Consider the set $S_p$ of points of 
$i_p\in I$ that have the property that, we can find $i_2,\dots ,i_p$ 
in $I$, for which $(i_\ell ,\tau (i_{\ell +1}))\in D$ for $1\le\ell 
< p$.  We observe that, obviously, $S_{p+1}\supset\tau^{-1}(D(S_p))$.  
Since we can assume without loss of generality that $\alpha^{m-1}\le 
N/2$, we see from (10.1) and induction that for $p\le m$, we have 
$\card S_p\ge \alpha^{p-1}$.  Then (10.1) shows that $\card S_{m+1}
\ge N/2$.  Then (10.2) shows that for $p\ge 1$, $N-\card S_{m+p+1}\le 
\alpha^{-p}N/2$.  Thus $N-\card S_{2m+1}\le \alpha^{-m}N/2<1$ 
which means $S_{2m+1}=I$.  Thus $i\in S_{2m+1}$.  Consider then the 
smallest $n$ for which $i\in S_{n+1}$; thus $n\le 2m$.  Then one can 
find $i_1=i$, $i_2,i_3,\dots ,i_{n+1}=i$ such that, for $1\le\ell\le 
n$ we have $(i_\ell ,\tau (i_{\ell +1}))\in D$.  The minimality of 
$n$ implies that the points $i_\ell$ are all disjoint.\rbx
\enddemo

Consider $u>0$, and consider the digraph $D_u$ given by $(i,j)\in D_u
\Leftrightarrow X_{i,j}\le 2uN^{-1}\log N$. 

\proclaim{Corollary 10.2}  Assume that the digraph $D_u$ is 
$\alpha$-expanding, and consider an integer $m$ such that $\alpha^m
\ge N/2$.  Then for an optimal assignment $\tau$ we have $X_{i,\tau 
(i)}\le  4muN^{-1}\log N$ for all $i\le N$.
\endproclaim

\demo{Proof}  Consider any $i\in I$, and consider $i=i_1,\dots ,
i_{n+1}=i$ as in Lemma 10.1, used for $D=D_u$.  Define $\sigma (i_\ell 
)=\tau (i_{\ell +1})$ for $1\le\ell\le n$, and $\sigma (i')=\tau (i')$ 
if $i'\not\in\{ i_1,\dots ,i_{n}\}$.  Since $\tau$ is optimal, we 
have
$$\sum_{i'\le N} X_{i',\tau (i')}\le\sum_{i'\le N}X_{i',\sigma (i')}$$
so that
$$\eqalignno{X_{i,\tau (i)} &\le \sum_{1\le\ell\le n} X_{i_\ell ,
\sigma (i_\ell )}\le 2nuN^{-1}\log N\,. &\bx \cr}$$
\enddemo

It remains to do computations.

\proclaim{Proposition 10.3}  For some constant $K$, and all $u>K$, 
$u\log N\le N$, the random digraph $D_u$ is $u\log N$-expanding with 
probability $\ge 1-N^{-{u \over K}}$.
\endproclaim

\demo{Proof}  We explain why (10.1) is satisfied with probability 
$\ge 1-N^{-u/K}$.  The case of (10.2) is similar and is left to the 
reader.  For simplicity, we set $\theta =uN^{-1}\log N$.
\enddemo

Consider a subset $S$ of $\{ 1,\dots ,N\}$ , and set $s=\card S$.  
For $j\in J$, we have
$$P(j\not\in D(S)) = \left( 1-{2u\log N \over N}\right)^s = (1-2
\theta )^s\le \exp (-2\theta s)$$
and thus
$$P(j\in D(S))\ge 1-\exp (-2\theta s)\,.$$
We observe that
$$0\le x\le 1\Rightarrow 1-e^{-x}\ge (1-e^{-1})x\,.$$
Thus, if we assume
$$\leqalignno{s\theta &\le {1 \over 2} &(10.3)\cr}$$
we have
$$\leqalignno{P(j\in D(S)) &\ge\gamma s\theta &(10.4)\cr}$$
where we have set $\gamma =2(1-e^{-1})>1$.

Consider $\gamma '=(1+\gamma )/2$.  We claim that, under (10.3) we 
have
$$\leqalignno{P(\card D(S)<\gamma 's\theta N) &\le \exp \left(-{s
\theta N \over K}\right)\,. &(10.5)
\cr}$$

This follows from (10.4) and the following general fact:

\proclaim{Lemma 10.4}  Consider independent events $(A_i)_{i\le N}$ 
with $P(A_i)=p$, and consider $\delta <1$.  Then, the probability 
than less that $\delta pN$ events occur is at most $\exp (-Np
/K(\delta ))$, where $K(\delta )$ depends on $\delta$ only.
\endproclaim

\demo{Proof}  Set $Y_i=1_{A_i}$, so that
$$E\exp (-\lambda Y_i)=1-p(1-e^{-\lambda})\le\exp (-p(1-e^{-\lambda})
)\,.$$
Thus
$$E\exp\left( -\lambda\sum_{i\le N}Y_i\right)\le\exp (-Np(1`-e^{-
\lambda}))\,.$$
By Chebyshev inequality we get
$$P\left(\sum_{i\le N}Y_i\le \delta pN\right)\le\exp Np(\lambda\delta 
-(1-e^{-\lambda}))$$
so the result follows by taking $\lambda$ small enough that $\lambda
\delta -(1-e^{-\lambda})<0$.\rbx
\enddemo

The number of subsets $S$ of $I$ of cardinal $s$ is at most $N^s$.  
For $u\ge K$, we have
$$N^s\exp\left( -{s\theta N \over K}\right)\le\exp\left( - {s\theta 
N \over K}\right)$$
and
$$\sum_{s\ge 1} \exp\left( - {s\theta N \over K}\right) \le N^{-u/K}
\,.$$
Thus, it follows that with probability $\ge 1-N^{-u/K}$, for all 
subsets $S$ of $I$ such that $s=\card S$ satisfies $\theta s \le 1/2$, 
we have $\card D(S)\ge\gamma '\theta Ns$.  Equivalently, we have
$$\leqalignno{u\log N\card S &\le {N \over 2}\Rightarrow\card D(S)\ge
\gamma 'u\log N\card S\,. &(10.6)\cr}$$
To complete the proof that (10.1) holds for $\alpha =u\log N$, it 
suffices to show that $\card D(S)\ge N/2$ whenever $\alpha\card S\ge 
N/2$.  This follows by applying (10.11) to a subset $S'$ of $S$ 
for which $\card S'$ satisfies $\alpha\card S'\le N/2$ and is as 
large as possible.\rbx

We can now prove the main result.

\proclaim{Theorem 10.5}  Denote by $M$ a median of $L_N$.  Then (for 
$N\ge 3$),
$$\leqalignno{t\le \sqrt{\log N} &\Rightarrow P\left(\vert L_N-M\vert
\ge {Kt(\log N)^2 \over \sqrt{N}\log\log N}\right) \le 2\exp(-t^2) 
&(10.7)\cr
t\ge\sqrt{\log N} &\Rightarrow P\left(\vert L_N-M\vert\ge {Kt^3\log N 
\over \sqrt{N}\log t^2}\right) \le 2\exp (-t^2)\,. &(10.8)\cr}$$
\endproclaim

\demo{Proof}  {\bf Step 1.}  Consider $u\le N/(2\log N)$, $\alpha = u
\log N$ and the smallest $m$ such that $\alpha^m\ge N/2$.  Set $v=4mu
N^{-1}\log N$, and $Y_{i,j}=\min (X_{i,j},v)$.  Consider the r.v. 
$L^u_N$ defined as $L_N$ but using the costs $Y_{i,j}$ rather that 
$X_{i,j}$.  It follows from Corollary 10.2 that $L^u_N=L_N$ whenever 
$D_u$ is $\alpha$-expanding, so that by Proposition 10.3
$$\leqalignno{P(L_N=L^u_N) &\ge 1-N^{-u/K}\,. &(10.9)\cr}$$

{\bf Step 2.}  When $N^{-u/K}\le 1/2$, it follows from (10.9) that 
$M$ is also a median of $L^u_N$.  It then follows from (8.1.1) (and 
scaling) that for all $w>0$,
$$P(\vert L^u_N-M\vert\ge w)\le 2\exp\left( -{w^2 \over 4Nv^2}\right)$$
and combining with (10.9) we get
$$P(\vert L_N-M\vert\ge w)\le 2\exp\left( -{w^2 \over 4Nv ^2}\right) 
+N^{-u/K}\,.$$

{\bf Step 3.}  We choose the parameters.  We take $w=3\sqrt{N} tv$.  
If $t^2\le\log N$, we take $u=K$; if $t^2\ge\log N$, we take $u=Kt^2/
\log N$.

Theorem 10.6 follows easily.\rbx
\enddemo

\remark{Remark}  A simple computation using Theorem 10.6 shows that 
the standard deviation of $L_N$ is not more than $K(\log N)^2/\sqrt{N}
\log\log N$.
\endremark

\vfil\eject

\noindent {\bf 11.~~Geometric Probability}

{\bf 11.1.~~Irregularities of the Poisson Point Process}

In this Chapter we will consider $N$ points $X_1,\dots ,X_N$ that are 
independent uniformly distributed in $[0,1]^d$, where, except on 
Section 13.5, $d=2$, and we will study certain functionals $L(X_1,
\dots ,X_N)$ of this configuration $X_1,\dots ,X_N$ (that is $L$ will 
depend only on $\{ X_1,\dots ,X_N\}$ rather than on the order in which 
the points are taken).

One would like to think that the sample $X_1,\dots ,X_N$ is rather 
uniform on $[0,1]^2$; say, that it meets every subsquare of side $K/
\sqrt{N}$.  This is not the case; there are empty squares of side of 
order $(N^{-1}\log N)^{1/2}$ (an empty square will informally be 
called a hole).  More importantly, in exceptional situations there are 
larger empty squares.  Several of the functionals we will study 
have the property that, if one delete or add a point to a finite set 
$F$, the amount by which $L(F)$ can vary depends whether $F$ has a 
``large'' hole close to $x$.  Thereby the first task is to study the 
size and number of holes.

It is not convenient to work with the sample $X_1,\dots X_N$.  The 
difficulty is that what happens say, in the left half of $[0,1]^2$ 
(for example, there is an excess of points here) affects what happens 
in the right half (there must then be a deficit of points there).  
Rather, one will work with a Poisson point process of constant 
intensity $\mu$.  This process generates a random subset $\Pi$ 
($=\Pi_\mu$) of $[0,1]^2$ with the following properties:

\exam{(11.1.1)}{If $A$ and $B$ are disjoint (Borel) subsets of $[0,1
]^2$, $\Pi\cap A$ and $\Pi\cap B$ are independent.}

\exam{(11.1.2)}{If $A$ is a (Borel) subset of $[0,1]^2$, the r.v. 
$\card (\Pi\cap A)$ is Poisson of parameter $\mu\vert A\vert$, where 
$\vert A\vert$ denotes the area of $A$.}

Let us recall that a r.v. $Y$ is Poisson of parameter $\lambda$ if
$P(Y=k)=e^{-\lambda}\lambda^k/k!$ for $k\ge 0$.  Thus
$$E(e^{uY})=\sum_{k\ge 0}e^{uk}e^{-\lambda}{\lambda^k \over k!} = 
\exp (\lambda (e^u-1))\,.$$

For the convenience of the reader, we recall some simple facts.

\proclaim{Lemma 11.1.1}  If a r.v. $Y$ satisfies
$$\leqalignno{E(e^{uY}) &\le \exp (\lambda (e^u-1)) &(11.1.3)\cr}$$
for $u\ge 0$, then for
$$\leqalignno{P(Y\ge t) &\le \exp\left( -t\log{t \over e\lambda}
\right)\,. &(11.1.4)\cr}$$
\endproclaim

\demo{Proof}  One can assume $t\ge\lambda$.  Write
$$P(Y\ge t)\le e^{-tu}E(e^{uY})\,,$$
use (11.1.3) and take $u=\log (t/\lambda )$.\rbx
\enddemo

\proclaim{Lemma 11.1.2}  If the r.v. $Y$ is Poisson of parameter 
$\lambda$, then
$$P\left(Y\le {\lambda \over 8}\right)\le \exp -{\lambda \over 2}\,.$$
\endproclaim

\demo{Proof}  Write, for all $u\ge 0$
$$P\left( Y\le {\lambda \over 8}\right)\le\exp{\lambda \over 8}
uEe^{-uY}=\exp\left({\lambda u \over 
8}+\lambda (e^{-u}-1)\right)$$
and take $u=2$.\rbx
\enddemo

For $k\ge 1$, we denote by ${\Cal C}_k$ the family of the $2^{2k}$ 
``dyadic squares'' of side $2^{-k}$.  So the vertices of these squares 
are of the type $(\ell_12^{-k},\ell_2 2^{-k})$, $0\le\ell_1,\ell_2\le 
2^{-k}$, $\ell_1,\ell_2\in\Bbb N$.

For $C\in{\Cal C}_k$, we set
$$\eqalign{Z_C = 1&~~~\text{if}~~\card (C\cap\Pi )\le\mu 2^{-2k-3}\cr
Z_C =0 &~~~\text{otherwise}\,.\cr}$$
From (11.1.2) and Lemma 11.1.2 follow that $\delta_k=P(Z_C=1)$ 
satisfies 
$$\leqalignno{\delta_k &\le \exp (-\mu 2^{-2k-1})\,. &(11.1.5)\cr}$$
Now, for $u>0$,
$$\leqalignno{Ee^{uZ_C} &= 1-\delta_k+\delta_k e^u &(11.1.6)\cr
&=1+\delta_k(e^u-1)\le\exp\delta_k(e^u-1)\,.\cr}$$
By (11.1.1) the variables $(Z_C)_{C\in{\Cal C}_k}$ are independent; 
so that, by (11.1.6)
$$Ee^{u\sum_{C\in{\Cal C}_k}Z_C}\le\exp 2^{2k}\delta_k(e^u-1)$$
and by Lemma 11.1.1 we have
$$\leqalignno{P\left(\sum_{C\in{\Cal C}_k}Z_C\ge v\right) &\le\exp
\left( -v\log {v \over e2^{2k}\delta_k}\right)\,. &(11.1.7)\cr}$$

Observe that $n_k=\sum\limits_{C\in{\Cal C}_k}Z_C$ is simply the 
number of squares of ${\Cal C}_k$ that contain no more than $1/8$ of 
the expected number of points of $\Pi$ they should contain.  Combining 
(11.1.5) and (11.1.7) we see that
$$\leqalignno{\forall k\,,~~~P(n_k\ge 2e^22^{2k}\exp (-\mu 2^{-2k-1})
) &\le \exp (-2e^2 2^{2k}\exp (-\mu 2^{-2k-1}))\,. &(11.1.8)\cr}$$

We now fix a number $t$, and we study how the number $n_k$ can be 
controlled if one rules out an exceptional set of probability $\le 
e^{-t^2}$.  We assume $t\ge 1$, $\mu\ge 4$.

We denote by $k_1$ the largest integer such that
$$\leqalignno{e^2 2^{2k_1}\exp (-\mu 2^{-2k_1-1}) &\le t^2\,. 
&(11.1.9)\cr}$$
Thus, $k_1\ge 0$ and for $k>k_1$ we have
$$e^2 2^{2k}\exp (-\mu 2^{-2k-1}) \ge t^2\,.$$

We now observe that if $a>1$, we have $\sum\limits_{\ell\ge 0}\exp 
(-2^{2\ell}a)\le 2\exp (-a)$, so that, combining with (11.1.8)
$$\leqalignno{P(\forall k>k_1\,,\, n_k\le 2e^2 2^{2k}\exp (-\mu 2^{-2
k-1})) &\ge 1-2e^{-2t^2}\,. &(11.1.10)\cr}$$

\proclaim{Lemma 11.1.3}  If $t\le\sqrt{\mu}/K$, we have
$${t^2 \over \mu\delta_{k_1-1}} \ge {1 \over \sqrt{\delta_{k_1-1}}}
\,.$$
\endproclaim

\demo{Proof}  It suffices to show that $\sqrt{\delta_{k_1-1}}\le t^2/
\mu$.  Now, by (11.1.5) and (11.1.9),
$$\sqrt{\delta_{k_1-1}}\le\exp (-\mu 2^{-2k_1}) \le\left({t^2 \over 
e^2 2^{2k_1}}\right)^2\,.$$
Thus it suffices to show that $2^{2k_1}\ge t\sqrt{\mu}/e^2$, i.e. 
$2^{2(k_1+1)}\ge 4t\sqrt{\mu}/e^2$.  The function $f(x)=e^2x\exp 
(-\mu /2x)$ is increasing for $x>0$.  Thereby, since $f(2^{2(k_1+1)})
\ge t^2$ by definition of $k_1$, it suffices to show that $f(4t
\sqrt{\mu}/e^2)<t^2$, which is equivalent to $\exp (-ae^2/2)< 1/16a$ 
for $a=\sqrt{\mu} /4t$.\rbx
\enddemo

We now apply (11.1.7), taking $k=k_1-1$ and $v=e2^{2k_1} t^2/\mu$.  
We observe that by Lemma 11.1.3 and (11.1.5) we have, for $t\le 
\sqrt{\mu}/K$,
$$\log{v \over e2^{2(k_1-1)}\delta_{k_1-1}}\ge\log{1 \over 
\sqrt{\delta_{k_1-1}}}\ge \mu 2^{-2k_1}$$
so that
$$\leqalignno{P\left( n_{k_1-1}\ge {e2^{2k_1}t^2 \over\mu}\right) 
&\le \exp (-2t^2)\,. &(11.1.11)\cr}$$

We now go back to the sample $X_1,\dots ,X_N$ and state our 
conclusions.

\proclaim{Proposition 11.1.4}  Consider $t\le \sqrt{N}/K$.  Denote by 
$k_0$ the largest integer for which $2^{2k_0}\le N$.  There exists an 
integer $k_1\le k_0$ such that
$$\leqalignno{{1 \over K}\log {N \over t^2}\le 2^{2(k_0-k_1)} &\le K 
\log{N \over t^2}  &(11.1.12)\cr}$$
and such that with probability $\ge 1-Ke^{-t^2}$ we have the 
following properties, where $m_k$ denotes the number of squares $C$ 
of ${\Cal C}_k$ such that
$$\leqalignno{\card (C\cap\{ X_1,\dots ,X_N\}) &\le N2^{-2k-6}\,. 
&(11.1.13)\cr}$$
For each $k_1\le k\le k_0$, we have
$$\leqalignno{m_k &\le K2^{2k}\exp (-N2^{-2k-6}) &(11.1.14)\cr}$$
and we have
$$\leqalignno{m_{k_1-1} &\le K{2^{2k_1}t^2 \over N}\,. 
&(11.1.15)\cr}$$
\endproclaim

\demo{Proof}  {\bf Step 1.}  Consider the process $\Pi = \Pi_\mu$, 
for $\mu =N/8$.  It follows from (11.1.4) that with probability $\ge 
1-\exp (-N/K)$, we have $\card\Pi\le N$.  It is obvious that, 
conditionally on the event $\{\card\Pi\le N\}$, the number $n_k$ of 
squares $C$ of ${\Cal C}_k$ for which $\card (C\cap\Pi )\le N2^{-2k
-6}=\mu /8$ stochastically dominates the number $m_k$.  Thus it 
suffices to prove (11.1.13) to (11.1.15) for $n_k$ rather than $m_k$, 
since, as we consider only $t\le \sqrt{N}/K$, the term $\exp (-N/K)$ 
is swallowed by the term $K\exp (-t^2)$.

{\bf Step 2.}  We define $k_1$ as in (11.1.9).  We observe that, 
since $1\le N2^{-2k_0}\le 4$ and $t\le\sqrt{N}/K$, we can assume 
$t^2\le e^22^{2k_0}\exp (-\mu 2^{-2k_0-1})$, so that $k_1\le k_0$.  
By (11.1.9) and definition of $k_1$ we have
$$\exp\mu 2^{-2(k_1+1)-1}\le {e^2 2^{2k_1+2} \over t^2}\le {N \over 
t^2}$$
so that $\mu 2^{-2k_1}\le K\log {N \over t^2}$, and thus $2^{2(k_0-
k_1)}\le K\log{N \over t^2}$.  By (11.1.9),
$$\exp (\mu 2^{-2k_1-1})\ge {e^2 2^{2k_1} \over t^2}\ge {2^{2k_0} 
\over t^22^{2(k_0-k_1)}}\ge {N \over t^2}\left( K\log {N \over t^2}
\right)^{-1}$$
and this finishes the proof of (11.1.12).

{\bf Step 3.}  By (11.1.10), with probability $\ge 1-2e^{-2t^2}$, for 
each $k>k_1$ we have
$$\leqalignno{n_k &\le 2e^2 2^{2k}\exp (-N2^{-2k-4})\,. 
&(11.1.16)\cr}$$
Now we observe that $n_{k_1}\le n_{k_1+1}$.  This is obvious, 
since, if $C\in {\Cal C}_{k_1}$, one of the $4$ squares of 
${\Cal C}_{k_1+1}$ contained in $C$ must contain at most $\card (\Pi
\cap C)/4$ points.  Thereby, by (11.1.16), for each $k\ge k_1$ we must 
have
$$n_k \le 8e^2 2^{2k}\exp (-N2^{-2k-6})\,.$$

Also, (11.1.11) shows that, with probability $\ge 1-e^{-2t^2}$, we 
have
$$\leqalignno{n_{k_1-1} &\le {K2^{2k_1}t^2 \over N}\,. 
&(11.1.17)\cr}$$

The events described above occur simultaneously with probability 
$\ge 1-3e^{-2t^2}$.\rbx
\enddemo

Having studied when and how the sample $X_1,\dots ,X_N$ can have a 
``deficit'' of points, we study how it can have excesses of points.  
While Proposition 11.1.4 is central to this chapter, the 
following result will be used only in Section 11.4.

\proclaim{Proposition 11.1.5}  Consider the integer $k_0$ of 
Proposition 11.1.4, and consider $k_2\le k_0$.  For $k_2\le k\le 
k_0$ consider a number $r_k$ such that $2^{2k}\ge r_k\ge 2^{2k}t^2/N$.  
Then, with probability $\ge 1-Ke^{-t^2}$, the following occurs

\exam{(11.1.18)}{Given $k_2\le k\le k_0$, and given a set $S\subset
{\Cal C}_k$ with $\card S\le r_k$, then}
$$\card\{ i\le N\,;\,X_i\in\cup\{ C\colon C\in S\}\}\le KN2^{-2k}
r_k+r_k\log {e2^{2k} \over r_k}\,.$$
\endproclaim

\demo{Proof}  For a subset $U$ of $[0,1]^2$, we have
$$\leqalignno{P(\card\{ i\le N\,;\,X_i\in U\}\ge u) &\le \exp 
\left( -u\log{u \over eN\vert U\vert}\right)\,. &(11.1.19)\cr}$$
This follows from (11.1.3) and (the argument of) (11.1.6).

For a subset $S$ of ${\Cal C}_k$, denote $U_S$ the union of the 
elements of $S$.  It suffices to consider the sets $S$ with $\card S
=r_k$.  For these we get from (11.1.19)
$$P(\card\{ i\le N\,;\, X_i\in U_S\}\ge u) \le\exp\left( -u\log {u 
\over eN2^{-2k}r_k}\right)\,.$$
There are at most ${2^{2k} \choose r_k}\le\exp (r_k\log (e2^{2k}
/r_k))$ choices for $S$.  We take 
$$u=r_k\log (e2^{2k}/r_k)+e^3N2^{-2k}r_k\,.$$
Thus we see that
$${2^{2k} \choose r_k} \exp\left( -u\log {u \over eN2^{-2k}r_k}
\right)\le {2^{2k}\choose r_k} \exp (-2u)\le\left({r_k \over e2^{2k}}
\right)^{r_k}\exp (-t^2)\,.$$
Since $\left({x \over e2^{2k}}\right)^x\le 2^{-2k}$ for $1\le x\le 
2^{2k}$, we see that (11.1.18) occurs with probability at least 
$1-Ke^{-t^2}$.\rbx
\enddemo

{\bf 11.2.~~The Traveling Salesman Problem}

The Traveling Salesman Problem (TSP) requires, given $N$ points 
$x_1,\dots ,x_N$ in the plane, to find the shortest tour through 
these points; in other words, to minimize
$$\Vert x_{\sigma (N)}-x_{\sigma (1)}\Vert +\sum^{N-1}_{i=1}\Vert 
x_{\sigma (i)}-x_{\sigma (i+1)}\Vert$$
over all permutations $\sigma\in S_N$.  The charm of the TSP is that 
it is the archetype of untractable question.  In this section, we 
denote by $L(F)$ the length of the shortest through $F$, and we study 
the r.v. $L_N=L(X_1,\dots ,X_N)$ where $X_1,\dots ,X_N$ are 
independent uniformly distributed over $[0,1]^2$.

While the TSP is usually very hard, somewhat surprisingly, it turns 
out that as far as the concentration of $L_N$ is concerned, it is the 
easiest problem we will consider.  The reason for this is its good 
regularity properties.  The only fact we will use about the TSP is as 
follows.

\proclaim{Lemma 11.2.1}  Consider $F\subset [0,1]^2$, $C\in
{\Cal C}_k$, $G\subset C$, and assume that there is a point of $F$ 
within distance $2^{-k+2}$ of $C$.  Then
$$\leqalignno{L(F) &\le L(F\cup G)\le L(F)+K2^{-k}\sqrt{\card G}\,. 
&(11.2.1)\cr}$$
\endproclaim

\demo{Proof}  An essential property of the TSP is its monotonicity:  
$L(F)\le L(F\cup\{ x\})$, as is seen by bypassing $x$ in a tour 
through $F\cup\{ x\}$.  This implies the left side inequality in 
(11.2.1).  To prove the right hand side inequality, one first uses 
the (well known, elementary) fact that there is a tour through $G$ of 
length $\le K2^{-k}\sqrt{\card G}$, and one connects this tour 
to a tour of $F$.

\proclaim{Theorem 11.2.2}  Assume that the functional $L$ satisfies 
the regularity condition of Lemma 11.2.1.  Then, if $X_1,\dots ,X_N$ 
are independent uniformly distributed over $[0,1]^2$, for each $t\ge 
0$ the r.v. $L_N=L(X_1,\dots ,X_N)$ satisfies $P(\vert L_N-M\vert\ge 
t)\le Ke^{-t^2/K}$, where $M$ is a median of $L_N$.
\endproclaim

Since the TSP is the simplest case we will consider, we will give the 
shortest possible proof, which is considerably simpler than the 
original proof.  The idea of this proof is, however, a bit tricky; 
a more straightforward, but somewhat longer proof will be given in 
Section 11.3.

The basic idea of the whole chapter is as follows:  consider $\Omega 
=[0,1]^2$, and the subset $A(a)$ of $\Omega^N$ that consists of the 
$N$-tuples $y_1,\cdots ,y_N$ for which $L(y_1,\dots ,y_N)\le a$.  
When $a=M$ is the median of $L$, Proposition 2.1.1 shows that, except 
for a set of probability $2e^{-t^2}$, given $X_1,\dots ,X_N$, we can 
find $(y_1,\dots ,y_N)\in A(a)$ such that $\card J\le Kt\sqrt{N}$, 
where $J=\{ i\le N\,;\,X_i\not= y_i\}$.  Thus we have a tour through 
$\{ X_i\,;\,i\not\in J\}$ of length $\le M$.  The points $X_i$, $i\in 
J$ should be in average at distance $\le K/\sqrt{N}$ of the 
set $\{ X_i\,;\,i\not\in J\}$; so each of them can be inserted in the 
tour by lengthening the tour of at most $K/\sqrt{N}$; for a total 
lengthening $\le Kt$.  This would prove that $P(L_N\ge M+Kt)\le 
e^{-t^2}$.  The problem with this argument is that the points $X_i$, 
$i\in J$ could be precisely chosen among those which are much further 
than $K/\sqrt{N}$ from their closest neighbor.  So we have to find a 
way to show that this does not happen, or at least that the effect of 
this phenomenon does not affect the final result.  The idea of this 
section is to give appropriate weights $\alpha (X_i)$ to each point 
$X_i$ (the more isolated the point is, the higher its weight) and 
then to use Theorem 4.1.1, to minimize the influence of points with 
large weights.

For $x\in [0,1]^2$, throughout this chapter, $C_k(x)$ denotes the 
square $C\in {\Cal C}_k$ containing $x$.  Throughout this section, we 
will set $F=\{ X_1,\dots ,X_N\}$,
$${\Cal H}_k=\{ C\in{\Cal C}_k\,;\,\card (F\cap C)\le N2^{-2k-6}\}
\,,$$
and $m_k=\card{\Cal H}_k$.

We fix $t\le\sqrt{N}/K$, and we recall the integers $k_0$, $k_1$ of 
Proposition 11.1.4.

For $x\in [0,1]^2$, we define
$$\alpha (x)=\sup\{ 2^{-k}\,;\,k_1\le k\le k_0\,;\,\card (\{ X_1,
\dots ,X_N\}\cap C_k(x)\})\le N2^{-2k-7}\}$$
when the set on the right is non-empty; and we set $\alpha (x)=
2^{-k_0}$ otherwise.

\proclaim{Proposition 11.2.3}  With probability $\ge 1-K\exp (-t^2)$, 
we have
$$\leqalignno{\sum_{i\le N}\alpha^2 (X_i) &\le K\,. &(11.2.3)\cr}$$
\endproclaim

\demo{Proof}  It should be obvious that
$$\eqalign{\sum_{i\le N}\alpha^2 (X_i) &\le K+\sum_{k_1\le k\le k_0}
2^{-2k}\card (F\cap\cup\{ C\,;\,C\in{\Cal H}_k\})\cr
&\le K+\sum_{k_1\le k\le k_0}2^{-2k}\times N2^{-2k-6}\card{\Cal H}_k
\,.\cr}$$
By Proposition 11.1.4, with probability $\ge 1-Ke^{-t^2}$, for all 
$k_1\le k\le k_0$, we have
$$m_k=\card{\Cal H}_k\le K2^{2k}\exp (-N2^{-2k-6})\,.$$
The result then follows from the elementary fact that 
$\sum\limits_{k\le k_0}2^{-2k}\exp (-N2^{-2k-6})\le K/N$.\rbx
\enddemo

\proclaim{Proposition 11.2.4}  In order to prove Theorem 11.2.2, it 
suffices to prove Proposition 11.2.5 below.
\endproclaim

\proclaim{Proposition 11.2.5}  Consider $X_1,\dots ,X_N$, and a 
subset $J$ of $\{ 1,\dots ,N\}$.  Assume that
$$\leqalignno{\sum_{i\not\in J}\alpha (X_i) &\le Kt &(11.2.4)\cr
\card{\Cal H}_{k_1-1} &\le {K2^{2k_1}t^2 \over N}\,. &(11.2.5)\cr}$$
Then
$$\leqalignno{L(X_1,\dots ,X_N) &\le L(\{ X_i\,;\,i\in J\} )+K't\,, 
&(11.2.6)\cr}$$
where $K'$ depends on the constants in (11.2.4) and (11.2.5) only.
\endproclaim

\demo{Proof of Proposition 11.2.4}  To prove Theorem 11.2.1, since 
$L_N\le K\sqrt{N}$, it suffices to consider the case $t\le\sqrt{N}/K$.  
We fix such a $t$, and we consider $a$ such that $P(L_N\le a)\ge 
e^{-t^2}$.  We will prove that
$$\leqalignno{P(L_N &\ge a+Kt)\le Ke^{-t^2} &(11.2.7)\cr}$$
and this clearly implies the result.  The condition $P(L_N\le a)\ge 
e^{-t^2}$ means $P(A(a))\ge e^{-t^2}$ (where $P$ denotes now the 
product measure on $\Omega^N$).  If we combine Lemma 4.1.2 and Theorem 
4.1.1, we see that with probability $\ge 1-e^{-t^2}$, the set $\{ X_1,
\dots ,X_N\}$ has the property that we can find $(y_1,\dots ,y_N)\in 
A(a)$ for which
$$\sum_{i\not\in J} \alpha_i(X_i)\le Kt\sqrt{\sum\limits_{i\le N}
\alpha (X_i)^2}$$
where $J=\{ i\le N\,;\, X_i=y_i\}$.  Now, by Proposition 11.2.3 and 
Proposition 11.1.4 we can moreover assume, with probability $\ge 1-K
e^{-t^2}$ that $\sum\limits_{i\le N}\alpha (X_i)^2\le K$ and that 
(11.2.5) holds.  By Proposition 11.2.5, we then have
$$\eqalignno{L(X_1,\dots ,X_N) &\le L(\{ y_i\,;\, i\in J\} )+Kt\le 
a+Kt\,. &\bx\cr}$$
\enddemo

\demo{Proof of Proposition 11.2.5}  We set $F'=\{ X_i\,;\, i\in J\}$, 
$G=\{ X_i\,;\, i\not\in J\}$.  We have to incorporate the points of 
$G$ into a tour through $F'$ without lengthening too much the 
tour.

{\bf Step 1.}  For $0\le k< k_0$, we denote by $U'_k$ the collection 
of those $C\in{\Cal C}_k$ that satisfy $C\cap F'=\emptyset$; we set 
$U'_{k_0}={\Cal C}_{k_0}$, and, for $0\le k\le k_0$, we denote by 
$U_k$ the collection of those $C\in U'_k$ that are not included in 
any $C'\in U_{k-1}$.  Thus, if $C\in U_k$, its distance to $F'$ is 
$\le 2^{-k+2}$.

By repeated applications of Lemma 11.2.1, we see that
$$L(X_1,\dots ,X_N)\le L(\{ X_i\,;\, i\in J\})+K\sum_{0\le k\le k_0}
\,\sum_{C\in U_k}2^{-k}\sqrt{\card (G\cap C)}\,.$$
Thereby, it suffices to show that this double sum is $\le Kt$.

{\bf Step 2.}  We consider three types of terms:

{\bf Type 1:}  $\card (G\cap C)\ge N2^{-2k-7}$.

In that case, since $\alpha (X_i)\ge N^{-1/2}$, we have
$$\leqalignno{2^{-k}\sqrt{\card (G\cap C)} &\le {K \over \sqrt{N}}
\card (G\cap C)\le K\sum\{\alpha (X_i)\,;\, X_i\in C\cap G\}\,. 
&(11.2.7)\cr}$$

{\bf Type 2:}  $k\ge k_1$, $\card (G\cap C)<N2^{-2k-7}$.

In that case, the definition of $\alpha (X_i)$ shows that $\alpha 
(X_i)\ge 2^{-k}$ for $X_i\in C$.  Thus
$$2^{-k}\sqrt{\card (G\cap C)}\le 2^{-k}\card (G\cap C)\le \sum
\{\alpha (X_i)\,;\, X_i\in C\cap G\}\,.$$

We observe that the total contribution of the terms of Types 1 and 
2 is $<Kt$ by (11.2.4), since the union of the sets $U_k$ are 
disjoint by construction.

{\bf Type 3:}  $k<k_1$, $\card (G\cap C)<N2^{-2k-7}$.

{\bf Step 3.}  We control the contribution of the terms of Type 3.  
We denote by $V_k$ the union of the sets $C\in U_k$ for which $\card 
(G\cap C)<N2^{-2k-7}$.  Denoting by $\vert V\vert$ the area 
of $V$, the key observation is that, under (11.2.5) we have
$$\leqalignno{\left|\bigcup\limits_{k<k_1} V_k\right| &\le {Kt^2 
\over N}\,. &(11.2.8)\cr}$$

The reason is simply that if $C\in {\Cal C}_k$ satisfies $\card 
(G\cap C)<N2^{-2k-7}$, when $C\in U_k$, $\card (G\cap C)=\card (F
\cap C)<N2^{-2k-7}$, so that, among the $2^{2(k_1-k-1)}$ squares $C'$ 
of ${\Cal C}_{k_1-1}$ that are contained in $C$, at least half of 
them must satisfy $\card (C'\cap F)<N2^{-2(k_1-1)-6}$, so belong to 
${\Cal H}_{k_1-1}$.  Thereby the area of $\bigcup\limits_{k<k_1}V_k$ 
can be at most twice the area of the union of ${\Cal H}_{k_1-1}$.

There are $2^{2k}\vert V_k\vert$ 
sets $C$ of ${\Cal C}_k$ included in $V_k$.  Thus, by Cauchy-Schwarz, 
we have
$$\eqalign{\sum_{C\in{\Cal C}_k ,C\subset V_k}2^{-k}\sqrt{\card (G
\cap C)} &\le 2^{-k}\sqrt{\card (G\cap V_k) 2^{2k}\vert V_k\vert}\cr
&= \sqrt{\card (G\cap V_k)\vert V_k\vert}\,.\cr}$$

Using Cauchy-Schwarz again, the sum of these terms over $k<k_1$ is at 
most\hfil\break $\sqrt{\vert V\vert\card (G\cup V)}\le \sqrt{N\vert 
V\vert}$ where $V=\bigcup\limits_{k<k_1}V_k$.  This is less than $Kt$ 
by (11.2.8).\rbx
\enddemo

\noindent{\bf 11.3.~~The Minimum Spanning Tree}

A spanning tree of a finite subset $F$ of $\Bbb R^2$ is a connected 
set that is a union of segments (called {\it edges}) each of which 
joins two points of $F$.  Its length is the sum of the lengths of 
these segments.  We denote by $L(F)$ the length of the shortest ($=$ 
minimum) spanning tree of $F$.  An interesting difference with the 
TSP is that it can happen that $L(F\cup \{ x\})<L(F)$.  This is 
e.g. the case if $F$ consists of the three vertices of an equilateral 
triangle and $x$ is its center.

The regularity property of $L$ that we will use is as follows.

\proclaim{Lemma 11.3.1}  Consider $C\in {\Cal C}_k$ ($k\ge 1$) and a 
subset $F$ of $[0,1]^2$.  Assume that each $C'\in{\Cal C}_{k-1}$ that 
is within distance $2^{-k+5}$ of $C$ meets $F$. Consider a subset $G$ 
of $C$.  Then
$$\leqalignno{\vert L(F\cup G)-L(F)\vert &\le K2^{-k}\sqrt{\card G}
\,. &(11.3.1)\cr}$$
\endproclaim

\demo{Proof}  {\bf Step 1.}  The inequality
$$L(F\cup G)\le L(F)+K2^{-k}(\card G)^{1/2}$$
is proved as in the case of the TSP.  The problem is the reverse 
inequality.

Consider a minimum spanning tree of $F\cup G$.  We remove all the 
edges adjacent to $G$.  This breaks the spanning tree in a number of 
pieces; and we have to add edges to connect it again.  We will prove 
$2$ facts.

{\it Fact 1.}  There is at most $6$ $\card G$ pieces;

{\it Fact 2.}  Each of the pieces contains a point within distance 
$K2^{-k}$ of $C$.

Once this is known, we simply take a point in each of these pieces 
within distance $K2^{-k}$ of $C$.  We build a tour of length $\le K
2^{-k}(\card G)^{1/2}$ through these points to reconnect the pieces.

{\bf Step 2.}  {\it Proof of Fact 1.}  Consider three points $x$, 
$a$, $b$ of $F\cup G$, such that the segments $[x,a]$, $[x,b]$ both 
belong to a minimum spanning tree of $F\cup G$.  Then we must have 
$\Vert a-b\Vert\ge\Vert x-a\Vert$ for otherwise we could remove the 
edge $[x,a]$ and replace it by $[a,b]$ to get a shorter spanning tree.  
Similarly, we have $\Vert a-b\Vert\ge\Vert x-b\Vert$.  Thus the 
angle between the lines $xa$, $xb$ is at least $\pi /3$.  Thereby the 
spanning tree must contain at most $6$ edges adjacent to each point.  
Thus removing $k$ points and the edges adjacent creates at most $6k$ 
connected components.

{\bf Step 3.}  {\it Proof of Fact 2.}  Consider a finite set $H$ of 
$[0,1]^2$.  Consider $a$, $b$ in $H$, and assume that $[a,b]$ belongs 
to a minimum spanning tree of $H$.  We show that the ``lens'' 
$$\leqalignno{L_{a,b} &=\{ x\,;\,\Vert a-x\Vert <\Vert a-b\Vert\,,\,
\Vert b-x\Vert <\Vert a-b\Vert\} &(11.3.2)\cr}$$
does not meet $H$.  Indeed if we remove $[a,b]$ from the minimum 
spanning tree, we split $H$ into the component $H_a$ containing $a$ 
and the component $H_b$ containing $b$.  If there existed $c\in L_{a,
b}\cap H_a$, we could remove the edge $[a,b]$ from the minimum 
spanning tree, and replace it by $[c,b]$ to get a shorter spanning 
tree.  Similarly, $L_{a,b}\cap H_b=\emptyset$.

We apply the above result to $H=F\cup G$.  An edge $[a,b]$ from a 
minimal spanning tree of $H$ is such that $L_{a,b}$ does not contain 
a square $C'$ in ${\Cal C}_{k-1}$ within distance $2^{-k+5}$ of $C$, 
because it is assumed that all such squares meet $F$, hence $H$.  
Thus, if $a\in C$, then, clearly, $\Vert b-a\Vert\le K2^{-k}$.\rbx
\enddemo

The main result of this section is as follows.

\proclaim{Theorem 11.3.2}  Assume that the functional $L$ satisfies 
the regularity condition of Lemma 11.3.1.  Then, if $X_1,\dots ,X_N$ 
are independent uniformly distributed over $[0,1]^2$, the r.v. $L_N=
L(X_1,\dots ,X_N)$ satisfies
$$\forall t\ge 0\,,~~P(\vert L_N-M\vert\ge t)\le Ke^{-t^2/K}$$
where $M$ is a median of $L_N$.
\endproclaim

One central idea of the approach will be to condition with respect to 
$X_1,\dots X_m$, where $m=[N/2]$.  The size of the holes of $\{ X_1,
\dots ,X_N\}$ are then controlled by the sizes of the holes of 
$\{ X_1,\dots ,X_m\}$, independently of $X_{m+1},\dots X_N$.  The 
main part of the proof of Theorem 11.3.2 is to obtain the following 
statement.  We set $\Omega =[0,1]^2$.

\proclaim{Proposition 11.3.3}  Consider an integer $n$ with $\left|{N 
\over 2}-n\right|\le 1$.  We write $\Omega_1=\Omega^n$, $\Omega_2=
\Omega^{N-n}$;  we denote by $P_1$, $P_2$ the product measures on 
$\Omega_1$, $\Omega_2$ respectively.  Given $0<t<\sqrt{N}/K$, there 
exists a subset $H_t$ of $\Omega_1$ such that $P_1(H_t)\le K_1e^{-
t^2}$, and that, whenever $(x_1,\dots ,x_n)\not\in H_t$, the r.v.
$$L'=L'(X_{n+1},\dots ,X_N)=L_N(x_1,\dots ,x_n,X_{n+1},\dots ,X_N)$$
defined on $\Omega_2$ has the following property

\exam{(11.3.3)}{If $P_2(L'\le a)\ge e^{-t^2}$, $P_2(L'\ge b)\ge e^{-
t^2}$, then $b-a\le Kt$.}
\endproclaim

First, we prove that Proposition 11.3.3 implies Theorem 11.3.2.  To 
prove that theorem, it suffices to prove the following statement:

If $P(L_N\le a)\ge 2e^{-t^2/2}$, $P(L_N\ge b)\ge 2e^{-t^2/2}$, then 
$b-a\le Kt$.  

Consider the set 
$A=\{L_N\le a\}$ in $\Omega^N$.  We will write
$\Omega^N=\Omega_1\times\Omega_2$ ($\Omega_1=\Omega^n$; 
$\Omega_2=\Omega^{N-n}$) and $P=P_1\otimes P_2$.  Thus, given
$\omega_1\in\Omega_1$, we define $L'$ 
on $\Omega_2$ by $L'
(\omega_2)=L_N(\omega_1,\omega_2)$.  For $\omega_1\in\Omega_1$, we 
write
$$A(\omega_1)=\{\omega_2\in\Omega_2\,;\,(\omega_1,\omega_2)\in A\}
\,.$$
Since $P(A)\ge 2e^{-t^2/2}$, the set
$$C_1=\{\omega_1\in\Omega_1\,;\,P_2(A(\omega_1))\ge e^{-t^2/2}\}$$
satisfies $P_1(C_1)\ge e^{-t^2/2}$.  Consider $C_2=C_1\backslash H_t$, 
so that $P_1(C_2)\ge (e^{-t^2/2}-K_1e^{-t^2})$.  When $\omega_1\in 
C_2$, we have $P_2(L'\le a)\ge e^{-t^2/2}$, so that by (11.3.3) 
we have $P_2(L'\le a+Kt)\ge 1-e^{-t^2}$.  By Fubini theorem, we get 
$$\leqalignno{P(W_1) &\ge (1-e^{-t^2})P(C_2\times\Omega_2 ) 
&(11.3.4)\cr}$$
where $W_1=\{ L_N\le a+Kt\}\cap (C_2\times\Omega_2)$.

We observe that (11.3.3) implies
$$P_2(L'\ge b)\ge e^{-t^2}\Rightarrow P_2(L'\ge b-Kt)\ge 1-e^{-t^2}
\,.$$
Thus, we can apply the same argument as above to show that
$$\leqalignno{P(W_2) &\ge (1-e^{-t^2})P(\Omega_1\times D_2) 
&(11.3.5)\cr}$$
where $W_2=\{L_N\ge b-Kt\}\cap (\Omega_1\times D_2)$, and $P_2(D_2)
\ge e^{-t^2/2}-K_1e^{-t^2}$.  For $t$ large enough,
$$P((C_2\times\Omega_2)\backslash W_1)+P((\Omega_1\times D_2)
\backslash W_2)<P(C_2\times D_2)$$
so that $W_1\cap W_2\not=\emptyset$.\rbx

We now start the proof of Proposition 11.3.3.  Consider $x_1,\dots 
,x_n\in\Omega =[0,1]^2$, and set $F'=\{x_1,\dots x_n\}$.  Denote by 
$m'_k$ the number of squares of ${\Cal C}_k$ that do not meet 
$F'$.  We consider the integers $k_1$, $k_0$ of Proposition 11.1.4 
(defined using $n$ rather than $N$).  We define $H_t$ as the set of 
$n$-tuples $(x_1,\dots ,x_n)$ for which

\exam{(11.3.6)}{For each $k$, $k_1\le k\le k_0$, we have}

$$\leqalignno{m'_k &\le K2^{2k}\exp (-n2^{-2k-6})\cr
m'_{k_1-1} &\le {K2^{2k_1}t^2 \over n}\,. &(11.3.7)\cr}$$

\noindent Thereby, $P_1(H_t)\ge 1-Ke^{-t^2}$ by Proposition 11.1.4.

We now fix $(x_1,\dots ,x_n)$ such that (11.3.6), (11.3.7) hold and 
we start the proof of (11.3.3).  For $x\in [0,1]^2$, we denote by 
$\ell (x)$ the smallest integer $\ell$ such that there is $C\in
{\Cal C}_\ell$, $C$ within distance $2^{-\ell +4}$ of $C_\ell (x)$, 
such that $F'\cap C=\emptyset$.  Thus, by definition, we observe

\exam{(11.3.8)}{If $\ell =\ell (x)$, any square $C'\in{\Cal C}_{\ell 
-1}$ that is within distance $2^{-\ell +5}$ of $C_\ell (x)$ meets 
$F'$.}

We also observe that if $y\in C_{\ell (x)}(x)$, then $\ell (y)=\ell 
(x)$, so that $V_\ell =\{ x\,;\,\ell (x)=\ell\}$ is a union of squares 
of ${\Cal C}_\ell$.

\proclaim{Lemma 11.3.4}  a)  We have, for each $k_1\le k\le k_0$ that
$$\leqalignno{\vert V_k\vert &\le K\exp (-n2^{-2k-6}) &(11.3.9)\cr}$$

$$\leqalignno{\left|\bigcup\limits_{\ell <k_1}V_\ell\right| &\le {K
t^2 \over n} \le {Kt^2 \over N}\,. &\text{b)}\cr}$$
\endproclaim

\demo{Proof}  Let us denote by $U'_\ell$ the union of the elements of 
${\Cal C}_\ell$ that do not meet $F'$, and set $U_\ell=U'_\ell
\backslash\bigcup\limits_{k<\ell}U'_k$.  It suffices to observe 
that if $x\in V_\ell$, then $C_{\ell (x)}$ is within distance $2^{-
\ell +4}$ of $U_\ell$, so that $\vert V_\ell\vert\le K\vert U_\ell
\vert$, and the result follows from (11.3.6), (11.3.7).\rbx
\enddemo

We consider the function $g(x)=2^{-\max (k_1,\ell (x))}$.  By 
(11.1.12), we have
$$\leqalignno{\Vert g\Vert_\infty &\le 2^{-k_1}\le {K \over\sqrt{n}}
\left(\log {n \over t^2} \right)^{1/2}\le {K \over t}\,. 
&(11.3.10)\cr}$$
By (11.3.9) and an obvious computation, we have
$$\leqalignno{\Vert g\Vert_2 &\le K/\sqrt{n}\,. &(11.3.11)\cr}$$

To prove (11.3.3), we have to prove that if $a$, $b$ are such that 
$P_2(L'\le a)\ge e^{-t^2}$, $P_2(L'\ge b)\ge e^{-t^2}$, then $b-a\le 
Kt$.  We now appeal to Corollary 2.4.5 with $u=Kt$, for the function 
$h(x,y)=g(x)+g(y)$.  From (11.3.10), (11.3.11), we see that we can 
find  $y_{n+1},\dots ,y_N$, $z_{n+1},\dots ,z_N$ such that
$$\leqalignno{L'(y_{n+1},\dots ,y_N) &\le a\,;~~L'(z_{n+1},\dots ,
z_N)\ge b\,, &(11.3.12)\cr}$$
and
$$\leqalignno{\sum_{i\in J}(g(y_i)+g(z_i)) &\le Kt &(11.3.13)\cr}$$
where $J=\{ n+1\le i\le N\,;\,y_i\not= z_i\}$.

Consider the set $F$ that consists of the points $x_1,\dots ,x_n$, 
as well as the points $y_i$, $i\not\in J$.  We will prove
$$\leqalignno{\vert L'(y_{n+1},\dots ,y_N)-L(F)\vert &\le Kt\,. 
&(11.3.14)\cr}$$

The same argument will show that
$$\vert L'(z_{n+1},\dots ,z_N) -L(F)\vert \le Kt$$
and this will finish the proof.

First we observe from (11.3.1) that if $F_1\supset F$, and if $\ell 
(x)\ge k_1$, then $\vert L(F_1\cup\{ x\})-L(F)\vert\le Kg(x)$.  
Thereby, it follows from (11.3.13) that we can add to $F$ all
the points $y_i$, $i\in J$, for which $\ell (y_i)\ge k_1$, without 
changing the value of $f$ by more than $Kt$.  Denote by $G$ the set 
of the other points $y_i$.  We observe that $G$ is contained in 
$\bigcup\limits_{\ell <k_1}V_\ell$.  Consider $C\in{\Cal C}_\ell$, 
$C\subset V_\ell$.  By (11.3.1), we have, for any set $F_1$ 
containing $F$, that
$$\vert L(F_1\cup (G\cap C)) - L(F_1)\vert\le K2^{-\ell}(\card G
\cap C)^{1/2}\,.$$
Thereby it suffices to show that
$$\sum_{\ell <k_1} 2^{-\ell}\sum_{C\subset V_\ell}(\card G\cap 
C)^{1/2}\le Kt\,.$$
But this is shown as in Step 3 of the proof of Proposition 
11.2.5.\rbx

\noindent{\bf 11.4.~~Gabriel Graph and Voronoi Polygons}

Given a subset $F$ of $[0,1]^2$ its Gabriel graph is the set of edges 
$[a,b]$ such that the closure $\overline{L}_{a,b}$ of the set 
$L_{a,b}$ of (11.3.2) meets $F$ only in $a$ and $b$.  When the set 
$F$ has the property that it does not contain points $x$, $y$, $z$ 
such that $\Vert x-y\Vert =\Vert x-z\Vert$, (a property that is 
satisfied with probability one for random sets) this is equivalent 
to saying that $F$ contains the edge $[a,b]$ if and only if $L_{a,b}$ 
does not meet $F$.  In that case, the Gabriel graph contains the 
minimum spanning tree, as is shown in the course of the proof 
of Lemma 11.3.1.  As in the case of the MST, at most $6$ edges are 
adjacent to each point of $F$.

We denote by $L(F)$ the length of the Gabriel graph.  An interesting 
feature of this functional is that in certain special configurations 
adding a single point creates a big decrease of $L(F)$.  A typical 
such configuration consists of the points $(0,k/n), (1,k/n)$, $0\le k
\le n$.  The Gabriel graph contains all the edges between $(0,k/n)$ 
and $(1,k/n)$.  All these edges will disappear when one adds the 
middle of the unit square to $F$.  The following lemma shows that the 
previous example is close to be the worst possible behavior.

\proclaim{Lemma 11.4.1}  Consider $C\in{\Cal C}_k$, $F$ a subset of 
$[0,1]^2$, and assume the following

\exam{{\rm (11.4.1)}}{Every element $C'$ of ${\Cal C}_{k-1}$ that is 
within distance of $2^{-k+3}$ of $C$ meets $F$.}

Then if $G\subset C$, we have
$$\leqalignno{\vert L(F)-L(F\cup G)\vert &\le K2^{-k}\card\{(F\cup 
G)\cap B(C,K2^{-k})\} &(11.4.2)\cr}$$
where $B(C,r)$ denotes the set of points within distance $r$ of $C$.
\endproclaim

\remark{Comment}  The difference with Lemma 11.3.1 is that the bound 
depends now upon $F\cup G$ rather than $G$ alone.
\endremark

\demo{Proof}  As already seen, a point is adjacent to at most $6$ 
edges, and, as in the case of the MST, edges adjacent to $G$ have a 
length $\le K2^{-k}$.  Thus
$$L(F\cup G)\le L(F)+K2^{-k}\card G\,.$$

To prove the reverse inequality, we observe that the edges $[a,b]$ 
that belong to the Gabriel graph of $F$ but not to the Gabriel graph 
of $F\cup G$ are exactly these for which $\overline{L}_{a,b}
\backslash\{ a,b\}$ meets $G$ but not $F$.  Then $\Vert a-b\Vert\le
K2^{-k}$, for otherwise there would exist $C'\in{\Cal C}_{k-1}$ 
within distance $2^{-k+3}$ of $C$ that would not meet $F$.  This 
implies, since $\overline{L}_{a,b}$ meets $G$, that $a,b\in B(C,K
2^{-k})$.  In the Gabriel graph of $F$, there are at most $6\cdot
\card (F\cap B(C,K2^{-k}))$ edges adjacent to points in $B(C,K
2^{-k})$, so at most that many edges can be removed.\rbx

Another natural example of functional that satisfies Lemma 11.4.1 is 
the total length of the Voronoi polygons.  If $F$ is a subset of $[0,
1]^2$, and $x\in F$, let us define the Voronoi polygon $V_x$ 
of $x$ as the set of all points $y$ of $[0,1]^2$ for which $d(x,y)=
d(y,F\backslash\{ x\})$.  (This name is a bit abusive since when $x$ 
is close to the boundary of $[0,1]^2$ this set is not a polygon).  
Denote by $L(F)$ the sum of the lengths of the Voronoi polygons of 
all points of $F$.  We sketch a proof that $L(F)$ satisfies the 
condition of Lemma 11.4.1.  First, we observe that if $y\in V_x$, 
there is no point of $F$ within distance less than $\Vert x-y\Vert$ 
of $y$.  Thus, if $x\in G$, the Voronoi polygon of $x$ (with respect 
to $F\cup G$) is under (11.4.1) entirely contained in $B(x,K
2^{-k})$, so is of length $\le K2^{-k}$.  Thus $L(F\cup G)\le L(F)+
K2^{-k}\card G$.  To prove the reverse inequality, consider a point 
$a$ belonging to the Voronoi polygon of $x\in F$, with respect 
to $F$, but not with respect to $F\cup G$.  Then there is no point 
of $F$ within distance less than $\Vert x-a\Vert$ of $a$, but there 
is at least a point of $G$.  Under (11.4.1) we have $a,x\in B(C,
K2^{-k})$; but the total length of the part of the Voronoi polygons 
of $F$ contained in $B(C,K2^{-k})$ is easily seen to be $\le K2^{-k}
\card (F\cap B(C,K2^{-k}))$.

\proclaim{Theorem 11.4.2}  Consider a functional that satisfies the 
condition of Lemma 11.4.1.  Set, as usual, $L=L_N=L(X_1,\dots ,X_N)$, 
and consider the median $M$ of $L_N$.  Then
$$\leqalignno{\forall t >0\,,~~P(\vert L-M\vert\ge t) &\le K\exp
\left(-{1 \over K}\min (t^2,(t \sqrt{N})^{2/3})\right)\,. 
&(11.4.3)\cr}$$
\endproclaim

In particular, the tails of $L_N$ are subgaussian for values of $t$ 
up to $N^{1/4}$.  We now sketch, in the case of the Gabriel graph, 
why, within logarithmic terms, the exponent in (11.4.3) is correct 
for $t\ge N^{1/4}$.  We give an informal argument, that could be 
made rigorous.  For simplicity, let us argue about $L(\Pi)$, where 
$\Pi$ is a Poisson point process of intensity $N$.  Consider $u\le
\sqrt{N}$, and let $a=u/\sqrt{N}\le 1$.  Denote by $k$ the cardinality 
of $\Pi\cap [0,a]^2$.  When $k$ is even, conditionally on $k$, with 
probability $\ge (1/Kk^4)^k$, the $k$ points of $\Pi\cap [0,a]^2$ are 
such that each of the discs of center $(\eta ,2\ell a/k)$, for $\eta
\in\{ 0,a\}$, $1\le\ell\le k/2$, and of radius $a/4k^2$ contains 
exactly one of these points.  Then the Gabriel graph of $\Pi$ contains 
the edge from the point in the disc of center $(0,2\ell a/k)$ to the 
point in the disc of center $(a,2\ell a/k)$, for a total length of 
order $ka$.  Now with overwhelming probability $k$ is of order $u^2$; 
so, with probability $\ge (1/K u^8)^{u^2}$ we get the exceptional 
configuration described above that creates an abnormal length of 
order $t=u^2a=u^3/\sqrt{N}$.  Now $u=(t\sqrt{N})^{1/3}$, and
$$\left({1 \over Ku^8}\right)^{u^2}\ge \exp\left( -{1 \over K}(t
\sqrt{N})^{2/3}\log t\sqrt{N} \right)\,.$$
So this later quantity is a lower bound on the probability that we 
get an abnormal length of order $t$ that will have $L$ exceed the 
median by $t$.

To prove Theorem 11.4.2, we observe that, since $\vert L_N\vert\le K
N$ by (11.4.2) it suffices to prove (11.4.3) for $t\le N/K$.  We 
follow the scheme of Section 11.3.  It suffices to be able to modify 
Proposition 11.3.3, so that when $t\le \sqrt{N}/K$ (11.3.3) can be 
replaced by

\exam{(11.4.4)}{If $P_2(L'\le a)\ge 2e^{-t^2}$, $P_2(L'\ge b)\ge 2
e^{-t^2}$, then $b-a\le K\left( t+{t^3 \over \sqrt{N}}\right)$.}

Once this is known, as in Section 11.3, we prove that
$$P(L\le a)\ge 2e^{-t^2/2}\,,~~P(L\ge b)\ge 2e^{-t^2/2}~
\text{imply}~b-a\le K\left(t+{t^3 \over \sqrt{N}}\right)\,.$$
This implies Theorem 11.4.2 since, if we set $u=t+t^3/\sqrt{N}$, for 
$u\le N/K$ we have $t\le\sqrt{N}/K$; moreover, we have $t^2\ge K^{-1}
\min (u^2,(u\sqrt{N})^{2/3})$.

The construction of $H_t$ and the proof of (11.4.4) will parallel 
the proof of Proposition 11.3.3.  In order to avoid repetition, we 
will not repeat the entire argument, but simply explain the
necessary modifications.

The construction of $H_t$ is modified as follows.  We require that 
for $k_1-1\le k\le k_0$, and each subset $S$ of ${\Cal C}_k$, with 
$\card S\le r_k$, then
$$\leqalignno{\card\{ i\le n\,;\,x_i\in\cup\{C\,;\,C\in S\}\} 
&\le KN2^{-2k}r_k+r_k\log {e2^{2k} \over r_k} &(11.4.5)\cr}$$
where we set $r_{k_1-1}=2^{2k_1}t^2/n$ and for $k\ge k_1$ we set
$$r_k=Kt2^{4k -3k_0}\,.$$
We observe that, using (11.1.12)
$$2^{-2k}r_k\ge Kt2^{2k_1-3k_0}\ge {Kt \over \sqrt{N}} 2^{2k_1-2k_0}
\ge {t^2 \over n}$$
provided $K$ is large enough.  It then follows from Proposition 
11.1.5 that imposing these extra conditions does not change the fact 
that $P_1(H_t)\ge 1-Ke^{-t^2}$.

We change the definition of the function $g(x)$ to 
$$g(x)={1 \over 2^{k_0}}(2^{k_0-\max (k_1,\ell (x))})^4\,.$$
Thus
$$\leqalignno{\Vert g\Vert_\infty &\le {1 \over 2^{k_0}} (2^{k_0-
k_1})^4\le {K \over \sqrt{n}}\left( \log {Kn \over t^2}\right)^2\le 
{K \over t} &(11.4.6)\cr}$$
and, obviously, (11.3.11) still hold.

Suppose now that we are given $a$, $b$ with
$$P_2(L'\le a)\ge 2e^{-t^2}\,,~~P_2(L'\ge b)\ge 2e^{-t^2}\,. $$

Using Proposition 11.1.5 again, we see that we can find a set $A
\subset\{ L'\le a\}$, $P_2(A)\ge e^{-t^2}$, such that whenever 
$(y_{n+1},\dots ,y_N)\in A$, we have

\exam{(11.4.7)}{For each $k_1-1\le k\le k_0$, for each subset $S$ 
of ${\Cal C}_k$, such that $\card S\le r_k$, then}
$$\card\{ n+1\le i\le N\,;\,y_i\in\cup\{C;C\in S\}\}\le KN2^{-2k}r_k
+r_k\log {e2^{2k} \over r_k}\,.$$

We then consider, using Proposition 11.1.5 again, a subset $B$ of 
$\{ L'\ge b\}$ with $P_2(B)\ge e^{-t^2}$, such that when $(z_{n+1},
\dots ,z_N)\in B$, the property similar to (11.4.7) holds.

We now appeal to Corollary 2.4.5, to find $(y_{n+1},\dots ,y_N)\in A$, 
$(z_{n+1},\dots ,z_N)\in B$ such that if $J=\{ i\,;\, n+1\le i\le N\,,
\,y_i\not= z_i\}$, then
$$\leqalignno{\sum_{i\in J} g(y_i)+g(z_i) &\le Kt\,. &(11.4.8)\cr}$$

We denote by $F$ the collection of points that consists of the points 
$(x_i)_{i\le n}$, together with the points $y_i$, $i\not\in J$.  We 
denote by $G$ the collection of points $y_i$, $i\in J$.  We have 
to show that
$$\leqalignno{\vert L(F\cup G)-L(F)\vert &\le K(t+t^3/\sqrt{n})\,. 
&(11.4.9)\cr}$$

Let us denote by $S_\ell$ the collection of squares $C\in
{\Cal C}_\ell$ that contain at least one point $y_i$, $i\in J$, $\ell 
(y_i)=\ell$.  It follows from (11.4.1) that, if $F\subset F_1\subset F
\cup G$, and if $C\in S_\ell$, we have
$$\vert L(F_1\cup (G\cap C))-L(F_1)\vert\le K2^{-\ell}\card\{ (F\cup 
G)\cap B(C,K2^{-\ell})\}\,.$$

Thereby, adding to $F_1$ all the points of $U_\ell\cap G$, where 
$U_\ell =\cup\{ C\,;\,C\in S_\ell\}$, we cannot change the value of 
$L$ by more than
$$\leqalignno{2^{-\ell}\card\{ (F\cup G)\cap B(U_\ell ,K2^{-\ell})\}
)\,. &&(11.4.10)\cr}$$
Since by definition, for $\ell (y_i)=\ell\ge k_1$, we have
$$g(y_i)\ge {1 \over 2^{k_0}} (2^{k_0-\ell})^4$$
and since $\sum\limits_{i\in J}g(y_i)\le Kt$ by (11.4.8) we see that
$$\card S_\ell\le Kt2^{4\ell -3k_0}\,.$$

Now, $B(U_\ell ,K2^{-\ell})$ is contained in a union of $\le K\card 
S_\ell$ squares $C$ of ${\Cal C}_\ell$.  Thereby, it follows from 
(11.4.5), (11.4.7) that the quantity (11.4.10) is bounded by
$$2^{-\ell}K(N2^{-2\ell}r_\ell +r_\ell \log {e2^{2\ell} \over r_\ell}
)\le Kt2^{\ell -k_0}$$
and these quantities have a sum $\le Kt$.

Now we have to control the influence of the points $y_i$ for which 
$\ell (y_i)<k_1$.

We denote by $V_\ell$ the set $\{\ell (x)=\ell\}$.  We recall that by 
Lemma 11.3.4 we have $\vert\bigcup\limits_{\ell <k_1}V_\ell\vert\le 
Kt^2/N$.  Since $V_\ell$ is union of squares of ${\Cal C}_\ell$, 
we have in particular that $V_\ell =\emptyset$ for $\ell\le k_3$, 
where $2^{-k_3}\le Kt/\sqrt{N}$.  Adding to a set $F_1$ such that 
$F\subset F_1\subset F\cup G$ the points of\break $G\cap V_\ell$, 
can, by (11.4.1), change the value of $L$ by at most
$$2^{-\ell}\card ((F\cup G)\cap B(V_\ell,K2^{-\ell}))\,.$$
Now we observe that $\vert B(V_\ell ,K2^{-\ell})\vert\le K\vert 
V_\ell\vert$.  Thus the total contribution of the points of $G\cap 
V_\ell$ is bounded by
$$\leqalignno{2^{-\ell}\card ((F\cup G)\cap V) &&(11.4.11)\cr}$$
where $\vert V\vert\le Kt^2/N$ and $V$ is a union of squares of 
${\Cal C}_{k_1-1}$.  The summation of all these quantities over 
$\ell\ge k_3$ is a most, using (11.4.7)
$$\eqalign{K2^{-k_3}\card ((F\cup G)\cap V) &\le {Kt \over \sqrt{N}}
\left( N2^{-2k_1}r_{k_1-1}+r_{k_1-1}\log{e2^{2k_1-2} \over r_{k_1-1}}
\right)\cr
&\le {Kt^3 \over \sqrt{N}}\left( 1+{2^{2k_1} \over N}\log {KN \over 
t^2}\right)\,.\cr}$$
But, using the definition of $k_1$, the last term is easily seen to 
be bounded by a constant.\rbx

\noindent {\bf 11.5.~~Simple matching}

In this section (for reasons that will become apparent later) we work 
in $[0,1]^d$ for $d\ge 2$.

A {\it matching} of a set $F$ is a decomposition of $F$ as a union 
of disjoint pairs of points (points of the same pair are matched); we 
make the convention that when $\card F$ is odd, there is exactly 
one point that is unmatched (does not belong to any pair).  A minimum 
matching is a matching that minimizes the sum of the distances of 
pairs of matched points.  We denote by $L(F)$ the length of a 
minimum matching of $F$.  For simplicity, the point to which a given 
point is matched is called its partner.

Our interest in that functional stems from the fact that it 
apparently does not have good regularity properties.  It is obvious 
that 
$$L(F\cup \{ x\})-L(F\cup\{ y\})\le \Vert x-y\Vert\,,$$
but in certain configurations this cannot be improved upon.  The 
problem is that if one tries to match $y$ to a point different from 
the partner of $x$, the partner of $x$ has to find a new partner, 
etc., and there is no apparent way to control this chain reaction.

While the behavior of $F$ is not good as far as the change of one 
point of $F$ is concerned, the situation is somewhat better when a 
significant number of points of $F$ are changed.  We set $L'(F)=\sup
\{ L(F')\,;\, F'\subset F\}$.

\proclaim{Lemma 11.5.1}  $\vert L(F)-L(G)\vert\le L'(F\triangle G)+
\sqrt{d}$.
\endproclaim

\demo{Proof}  Consider $U=F\backslash G$, $V=G\backslash F$.  Consider 
a minimal ${\Cal M}$ matching of $F$, and, for $a\in F$, denote its 
partner by $\theta (a)$.  Consider
$$H=\{ a\in F\backslash U\,;\,\theta (a)\in U\}\,.$$
When we remove $U$ from $F$, the points of $H$ lose their partners.  
Set $H'=\{\theta (a)\,;\,a\in H\}$.  Thus $H'\subset U$.  To find 
partners for the points of $V\cup H$ we consider a minimum 
matching of $V\cup H'$.  This matching induces a matching ${\Cal M}'$ 
of $V\cup H$, using the bijection $\theta$ of $H$ and $H'$.  The 
union of the trace of ${\Cal M}$ on $F\backslash (U\cup H)$ and 
${\Cal M}'$ is almost a matching of $G$, although it could happen 
that there remains an unmatched point in $V\cup H$ and one in $F
\backslash (U\cup H)$.  Then two points are then matched 
together (creating the term $\sqrt{d}$).  The matching we have 
constructed witnesses that
$$\eqalign{L(G) &\le L(F)+L(V\cup H')+\sqrt{d}\cr
&\le L(F)+L'(V\cup U)+\sqrt{d}\,. \cr}$$
To see it, it suffices to use the triangle inequality, and to observe 
that the edges  $[a,\theta (a)]$ for $a\in U'$ do disappear from 
${\Cal M}$ when $U$ is removed.\rbx
\enddemo

Here is a simple observation.

\proclaim{Lemma 11.5.2}  Consider subsets $F_1,\dots ,F_p$ of $[0,1
]^d$.  Then
$$L'\left(\bigcup\limits_{i\le p} F_i\right) \le \sum_{i\le p}L'
(F_i)+Kp^{1-1/d}$$
where, as in the rest of this section, $K$ denotes a constant that 
depends on $d$ only.
\endproclaim

\demo{Proof}  It suffices to prove this for $L$ rather than $L'$.  
The point is that if one considers an optimal matching of each $F_i$, 
their union fails to be a matching of $\bigcup\limits_{i\le p}F_i$ 
only because there could remain an unmatched point in each $F_i$, 
while we are permitted at most a single unmatched point.  Thus, it 
suffices to match all but at most one of these points, using for 
example a shortest tour through them, and matching consecutive points 
on the tour.\rbx
\enddemo

It seems an interesting question whether when $d=2$ the inequality of 
Theorem 11.2.3 would hold, at least for smaller values of $t$.  
Possibly easier is the question whether the variance of $L_N$ is 
bounded.  The best results in that direction belong to Rhee.  She 
proved that if $d=2$, $\Var L_N\le K(\log N)^2$ [R1], while if $d\ge 
3$, $\Var L_N\le KN^{1-1/d}$ [R2].  The arguments for these 
results are different.  Our methods do not allow to improve on the 
result for $d=2$, but allow significant improvement when $d\ge 3$ 
(and this is why we consider this case in this section).  Although 
this has not been checked, it seems to be an exercise to show that 
$\Var L_N\ge {1 \over K}N^{1-2/d}$ using e.g., the method of [R3].  
What we will prove is that $\Var L_N\le (\log N)^K N^{1-2/d}$.  The 
proof goes by first proving a Poissonized version of the result, and 
then using  ``dePoissonization''.  The second part of the argument is 
standard (see e.g. [R1]) and will not be given here.

The Poissonized version of the problems is the study of the r.v. 
$L_\lambda =L(\Pi_\lambda )$, where $\Pi_\lambda$ is the random 
subset of $[0,1]^d$ that is generated by a Poisson point process 
of constant intensity $\lambda$.  We consider the space $\Omega$ of 
all finite subsets of $[0,1]^d$; and on $\Omega$, we consider the 
probability $P_\lambda$ induced by $\Pi_\lambda$.  On $\Omega^2$, 
we consider the function
$$\leqalignno{f(F,G) &= L'(F\triangle G)\,. &(11.5.1)\cr}$$

For a subset $B$ of $\Omega$, we set
$$\leqalignno{f(F,B) &= \inf_{G\in B} f(F,G)=\inf_{G\in B}L'(F
\triangle G)\,. &(11.5.2)\cr}$$

We set $\gamma ={1 \over 2} -{1 \over d}$.

\proclaim{Theorem 11.5.3}  For all $\lambda\ge 3$, all subsets $B$ of 
$\Omega$, we have
$$\int_\Omega\exp{f(F,B) \over (\log\lambda )^K\lambda^\gamma} d
P_\lambda (F)\le {e \over P_\lambda (B)}\,.$$
\endproclaim

If we combine this result with Lemma 11.5.1 (and proceed as usual) we 
see that if $M_\lambda$ denotes a median of $L_\lambda$, we have
$$\int\exp\left({1 \over (\log\lambda )^K\lambda^\gamma}\vert 
L_\lambda -M_\lambda\vert\right) dP\le K$$
which certainly implies the previous claim about the variance of 
$L_\lambda$.  To prove Theorem 11.5.3, we will prove the following 
statement, which form is adapted to proof by induction.

\proclaim{Proposition 11.5.4}  There exists numbers $K_0$, $\alpha 
>1$ depending on $d$ only, such that for all $q>0$ we have, for all 
$\lambda$, $1\le\lambda\le 2^{\alpha^q}$ and all Borel subsets $B$ 
of $\Omega$,
$$\int_\Omega \exp {f(F,B) \over K^q_0\lambda^\gamma} dP_\lambda (F) 
\le {e \over P_\lambda (B)}\,.$$
\endproclaim

To see that this statement implies Theorem 11.5.3, we take for $q$ 
the smallest such that $\lambda\le 2^{\alpha^q}$, so that $\alpha^q$ 
is of order $\log\lambda$, and $K^q_0$ of order $(\log \lambda)^K$.

The proof of Proposition 11.5.4 is by induction over $q$.  For the 
case $q=1$, one uses the brutal bound 
$$f(F,G)\le K(\card F+\card G)$$
and the exponential integrability of Poisson random variables.  The 
easy details are left to the reader.

We will determine, in due time, suitable values for $K_0$ and $\alpha$ 
and we now start the proof of the induction step from $q$ to $q+1$.  
Consider $\lambda$ such that $2^{\alpha^q}\le\lambda\le 2^{\alpha^{q+
1}}$.  Consider the smallest integer $n$ such that $\lambda '=\lambda
/n^d\le 2^{\alpha^q}$.  (Thus, we can apply the induction hypothesis 
to $\lambda '$.)  By definition of $n$, we have $\lambda /(n-1)^d\ge 
2^{\alpha^q}$, so that, since $\lambda\le 2^{\alpha^{q+1}}$, we have 
$(n-1)^d\le 2^{(\alpha -1)\alpha^q}$, and thus
$$\leqalignno{n^d &\le {n^d \over (n-1)^d} 2^{(\alpha -1)\alpha^q}
\le 2^d\cdot 2^{(\alpha -1)\alpha^q}\,. &(11.5.3)\cr}$$

Also,
$$\leqalignno{\lambda ' &={\lambda \over n^d}\ge \left({n-1 \over n}
\right)^d {\lambda \over (n-1)^d}\ge\left({n-1 \over n}\right)^d 
2^{\alpha^q}\ge 2^{\alpha^q-d}\,. &(11.5.4)\cr}$$

Consider a partition of $[0,1]^d$ in $n^d$ congruent cubes 
$(C_i)_{i\le n^d}$.  From Lemma 11.5.2, we observe that
$$\leqalignno{L'(F\triangle G) &\le \sum_{i\le n^d}L'((F\triangle G)
\cap C_i) + Kn^{d-1}\,. &(11.5.5)\cr}$$

We set
$$f_i(F,G) = L'((F\triangle G)\cap C_i)\,.$$

Thus we have, from (11.5.5)
$$L'(F\triangle G)\le \sum_{i\le n^d} f_i(F,G)+Kn^{d-1}\,.$$
Thus, if we set
$$g(F,G)=\inf_{G\in B}~\sum_{i\le n^d}f_i(F,G)$$
we get by (11.5.2) that
$$\leqalignno{f(F,B) &\le g(F,B)+Kn^{d-1}\,. &(11.5.6)\cr}$$

The crucial point is that $(\Omega ,P_\lambda )$ is naturally 
isomorphic to the product of $n^d$ copies of $(\Omega ,P_{\lambda '}
)$.  To see this, let us denote by $R_i$ an affine map from $C_i$ to 
$[0,1]^d$, for $i\le n^d$.  Then the isomorphism simply associates 
$(R_i (F\cap C_i))_{i\le n^d}$ to $F$.  We observe that
$$f_i(F,G)={1 \over n}L'(R_i(F\cap C_i)\triangle R_i(G\cap C_i))$$
so that, under this isomorphism, each function $f_i$ is distributed 
like the function $h'$ on $\Omega^2$ (Provided with $P_{\lambda '}
\otimes P_{\lambda '}$), where $h'(F,G)={1 \over n}L'(F\vartriangle 
G)$.  Moreover, with the notation of Definition (2.4.1), we have 
$f_{h'}=g$.  By induction hypothesis, and taking the scaling factor
$n$ into account, we have for each Borel set $B\subset\Omega$,
$$\int_\Omega\exp (2h(F,B)) dP_{\lambda '}(F)\le {e \over P_{\lambda 
'}(B)}$$
where $h=ah'$, $a=n(2K^q_0\lambda'\,^{\gamma})^{-1}$.  It then 
follows from Theorem 2.5.1 and the definition of $g$ that
$$\forall t\le 1\,,~~\int_\Omega\exp (atg(F,B)) dP_\lambda (F)\le 
{1 \over P_\lambda (B)}\exp (3n^dt^2)$$
for each Borel set $B\subset\Omega$.  From (11.5.6), it follows 
that
$$\int_\Omega \exp (at f(F,B)) dP_\lambda (F)\le {1 \over P_\lambda 
(B)}\exp (3n^{d}t^2+Kn^{d-1}at)\,.$$
We see that if
$$\leqalignno{n^{d/2-1}a &\le K\,, &(11.5.7)\cr}$$
then, taking $t=n^{-d/2}/K$, we get
$$\leqalignno{\int_\Omega\exp\left({a \over Kn^{d/2}} f(F,B)\right) 
dP_\lambda (F) &\le {e \over P_\lambda (B)}\,. &(11.5.8)\cr}$$
Now,
$${a \over Kn^{d/2}} = {1 \over 2KK^q_0\lambda '\,^{\gamma}n^{d/2-1}} 
={1 \over 2KK^q_0\lambda^\gamma}$$
since $d\gamma =d/2-1$.  Thus, provided $K_0=2K$, (11.5.8) is exactly 
what we need to complete the induction.

Thus, it remains to check that (11.5.7) holds; but by (11.5.3), 
(11.5.4)
$$\eqalign{n^d &\le 2^d 2^{(\alpha -1)\alpha^q} \cr
a/n &\le \lambda '\,^{-\gamma}\le 2^{\gamma d} 2^{-\gamma\alpha^q}
\cr}$$
so that (11.5.7) holds for $\alpha =1+2\gamma$.\rbx
\vfil\eject

\noindent{\bf 12.~~The free energy of Spin Glasses at high 
temperature}

Consider a sequence $(\epsilon_i )_{i \le N}$ with $\epsilon_i \in 
\{-1,1\}$.  Each $\epsilon_i$ represents the two possible values of 
the spin of particule $i$.  Consider numbers $(h_{ij} )_{1 \le i <j 
\le N}$ that represent the interaction between spins.  The energy of 
a given configuration is given by $\underset {1 \le i <j \le N}
\to \sum h_{ij} \epsilon_i \epsilon_j$.  Consider a parameter $\beta 
> 0$ (that plays the role of the inverse of the temperature).  The 
so-called ``partition function" is given by
$$\leqalignno{Z_N = Z_N(h_{ij}) = 2^{-N} \sum_{(\epsilon_i) \in \{-1,
1\}^N} \exp \Big( \frac \beta {\sqrt N} \sum_{1 \le i <j \le N} h_{ij}
\epsilon_i \epsilon_j \Big) &&(12.1)\cr}$$
The role of the factor $\sqrt N$ is for normalization purposes that 
will become apparent later.

If we think to $\epsilon_i$ as a Bernoulli r.v., it is natural to 
write
$$\leqalignno{Z_N (h_{ij}) = E_\epsilon \exp \big(\frac \beta 
{\sqrt N} \sum_{1 \le i <j \le N} h_{ij} \epsilon_i \epsilon_j \big). 
&&(12.2)\cr}$$

In the model we study, the numbers $h_{ij}$ are random, and the 
sequence $(h_{ij})_{1 \le i <j \le N}$ is i.i.d.  We assume $Eh_{ij} 
= Eh^3_{ij} = 0$, and we assume for normalization purposes that $E
h^2_{ij} = 1$.  We will also assume that $E\exp \alpha|h_{ij} |< 
\infty$ for $\alpha$ small enough.  Then $EZ_N$ is well defined for 
$N$ large enough.  We are interested in the quantity $N^{-1}E \log 
Z_N$ (mean free energy per site), whose study relies ultimately on 
the study of $Z_N$.  It is proved in [A-L-R], and in [C-N] in the 
case where $h_{ij}$ is gaussian, that for $\beta<1$ the random 
variable $\log Z_N- {\beta^2 N}/ 4 $ converges in law to a 
(non-standard) normal r.v.  Equally interesting, but of a rather 
different nature is the research of tail estimates for $\log Z_N - 
{\beta^2 N}/ 4 $ that are valid for all $N$.

\proclaim {Theorem 12.1}  There exists a universal constant $K$ with 
the following property.  Assume that $E \exp \pm h_{ij}< 2$.  Then, 
for $0 < t < N/K$, $\beta <1$, 
$$\leqalignno{P\Big( |\log Z_N - \frac {\beta^2 N} 4 | \ge K \Big(t + 
\sqrt{\log \frac K {1-\beta^2}} \Big) \sqrt N \Big) \le 2e^{-t^2} 
&&(12.3)\cr}$$
\endproclaim

In particular
$$\leqalignno{-\frac K {\sqrt N} \sqrt {\log \frac 2 {1-\beta^2}} \le 
\frac 1 N E \log Z_N -\frac {\beta^2} 4 \le \frac K N .&& (12.4)\cr}$$

\remark{Comment}  In the condition $E\exp \pm h_{ij} \le 2$, the 
number $2$ can be replaced by any other (with a different constant 
$K$).  It seems reasonable to conjecture that (12.3)  is not sharp 
in the gaussian case, and that, for a given $\beta < 1$,
$$\underset {t \rightarrow\infty} \to \lim \sup_N P \big( | \log Z_N- 
\frac {\beta^2 N} 4 |\ge t \big ) = 0. $$
\endremark

It should however be pointed out that (12.3) does not hold when the 
factor $\sqrt N$ is removed from (12.3).  Indeed it would follow 
otherwise that for each $n$, $\underset N \to\sup E (4Z_N / \beta^2 
N)^n < \infty$, and it is pointed out in [A-L-R], p. 6, that this is
not the case.

The key to Theorem 12.1 will be the following deviation inequality 
$$\leqalignno{ 0<t\le 4 \sqrt N(N-1)\Rightarrow P(|\log Z_N - M_N | 
\ge t) \le 2 \exp \Big ( \frac {-t^2} {32 (N-1)}\Big)&&(12.5) \cr}$$
where $M_N$ denotes a median of $\log Z_N$.  We first show how to 
deduce this from Corollary 2.4.4.  The second crucial step will then 
be to relate $M_N$ and $\beta^2N/4$ ($\approx \log EZ_N$).

To prove (12.5), we observe that 
$$\leqalignno{|\log Z_N (h_{ij}) - \log Z_N (h'_{ij})| \le 
\frac \beta {\sqrt N} \quad \sum_{1 \le i < j \le N} | h_{ij} - 
h'_{ij} |&&(12.6)\cr}$$
as follows from the fact that 
$$|\sum_{1 \le i < j \le N} a_{ij}\epsilon_i \epsilon_j | \le \sum_{1 
\le i < j \le N} | a_{ij}|$$

We now view $\log Z_N$ as a function on ${\Bbb R}^{N(N-1)/2}$.  We 
wish to apply Corollary 2.4.4, in the case $\Omega = {\Bbb R}$, 
$h(x,y)= \frac 1 4 |x - y |$, $\mu$ the law of $h_{ij}$.  We note 
that (2.4.12) holds, since
$$\eqalign{\iint_{{\Bbb R}^2} \exp \frac 1 4 |x - y| d \mu (x) d\mu 
(y) &\le \Big( \int \exp \frac 1 4 |x| d \mu (x) \Big)^2 \cr
&\le \big( E \exp |h_{ij}| \big)^{1/2} \cr
&\le \big( E (\exp h_{ij} + \exp - h_{ij}) \big)^{1/2} \le 2.\cr}$$
Consider now $v$ and the set $A=\{\log Z_N <v\}$.  Combining (12.6) 
and (2.4.13) (used for $N(N-1)/2$ rather than $N$) we see that for 
$u>v$, we have 
$$\eqalign{u- v &\le 4\beta \sqrt N (N-1)\cr
&\Rightarrow P(\{ \log Z_N> u\} ) P ( \{ \log Z_N<v\}) \le \exp 
\big( -\frac {(u- v)^2} {32 \beta^2 (N-1)} \big)}$$
Taking successively $u=M_N$ and $v=M_N$, (12.5) follows as usual.

In order to relate $M_N$ and $\beta^2 N/4$, the key step is the 
elementary estimates
$$\leqalignno{ \frac 1 K \exp \frac {\beta^2N} 4 \le EZ_N \le K \exp 
\frac {\beta^2N} 4 &&(12.7) \cr}$$
$$\leqalignno{EZ^2_N \le \frac K {1-\beta^2}  (EZ_N)^2.&&(12.8)\cr}$$
These will be proved later.  First, we conclude the main argument.  
Consider the set $A=\{Z_N \ge \frac 1 2 EZ_N \}$.  Then
$$\eqalign {EZ_N &= E(Z_N 1_{A^c}) + E(Z_N 1_A )\cr
&\le \frac 1 2 EZ_N + E(Z^2_N)^{1/2} P(A)^{1/2}\cr}$$
so that 
$$P(A) \ge \frac 1 4 \frac {(EZ_N)^2} {EZ^2_N}$$
(a fact going back to Paley and Zigmund.)  Combining with (12.8), we 
get $P(A) \ge (1-\beta^2) /K$.  To get a lower bound for $M_N$, we 
can assume $M_N \le \log \frac 1 2 EZ_N$.  We set $t=\log (\frac 1 2 
EZ_N) - M_N$.  Since $\log Z_N \ge 0$, we have $M_N \ge 0$ and hence 
$t \le K+N/4$ by (12.7).

We certainly have
$$A \subset \{ \log Z_N \ge M_N +t \}$$
Thus, by (12.5) we have
$$\frac {(1-\beta^2)} K \le P(A) \le 2 \exp \big( - \frac {t^2} {32
(N-1)}\big)$$
so that
$$t \le K \sqrt N \big( \log \frac K {1-\beta^2} \big)^{1/2}$$
and thus
$$M_N \ge \log \big( \frac 1 2 EZ_N \big) - K \sqrt N \big( \log \frac K 
{1-\beta^2}
\big)^{1/2}.$$
We also have $M_N \le \log (2EZ_N)$.  Combining with (12.7) we get
$$|M_N - \frac {\beta^2 N } 4 | \le K \sqrt N \big( \log \frac K {1 
- \beta^2} \big)^{1/2}$$
so that (12.3) now follows from (12.5).

To prove (12.4), we first observe that the lower bound follows from 
(12.3) and a routine computation.  The upper bound follows from the 
concavity of $\log$, which implies $E \log Z_N \le \log EZ_N$, and 
(12.7).

It remains to prove (12.7), (12.8).  We start with the elementary 
inequality
$$|e^x -1 -x -\frac{x^2} 2 -\frac {x^3} {3!} | \le \frac {x^4} {4!} 
e^{|x|}$$
that is obvious on power series expansions.  Thus, for $|u| \le \frac 
1 2$, we have (since $Eh^2_{ij} = 1$, $Eh_{ij} = Eh^3_{ij} =0$),
$$\leqalignno{1 + \frac {u^2} 2 - Ku^4 \le E \exp (uh_{ij}) \le 1+ 
\frac{u^2} 2 +Ku^4 &&(12.9)\cr}$$
and thus, for $\epsilon = \pm 1$, $\beta\le 1$,
$$\exp \Big( \frac {\beta^2} {2N} - \frac {K \beta^4} {N^2} \Big) 
\le E \exp \frac {\epsilon \beta h_{ij}} {\sqrt N} \le \exp \Big( 
\frac {\beta^2} {2N} - \frac {K \beta^4} {N^2} \Big).$$

Since  
$$EZ_N = E_\epsilon \prod_{ij} E \exp \epsilon_i \epsilon_j \frac 
{\beta h_{ij}}{\sqrt N},$$
(12.7) follows.  Turning to the study of $EZ_N^2$, we have, using 
(12.9), and for $N \ge 8$, that, with obvious notations,
$$\eqalign {EZ_N^2 &= E E_\epsilon E_{\epsilon '} \exp \Big( \sum_{1 
\le i < j \le N} \frac {\beta h_{ij}} {\sqrt N} (\epsilon_i \epsilon_j 
+ \epsilon_i' \epsilon_j')\Big) \cr
&\le KE_\epsilon E_{\epsilon '} \exp \Big( \sum_{1 \le i < j \le N} 
\frac {\beta^2}{2N} (\epsilon_i \epsilon_j + \epsilon_i' \epsilon_j'
)^2 \Big) \cr}$$
Now, we have $(\epsilon_i \epsilon_j + \epsilon_i' \epsilon_j')^2 = 2 
+2 \epsilon_i \epsilon_j \epsilon_i' \epsilon_j'$.  Also, $\epsilon_i 
\epsilon_j\epsilon_i' \epsilon_j'$ is distributed like $\epsilon_i 
\epsilon_j$, so that
$$EZ_N^2 \le K \exp (\frac {\beta^2 N} 2) E_\epsilon \Big(
\frac {\beta^2}{2N} ( \sum_{1 \le i < j \le N} 2 \epsilon_i 
\epsilon_j) \Big) $$
Now,
$$\sum_{1 \le i < j \le N} 2 \epsilon_i \epsilon_j = \Big(\sum_{1 
\le i \le N} \epsilon_i\Big)^2
-N.$$
Using the subgaussian inequality
$$P_\epsilon \Big( | \sum_{i=1}^N \epsilon_i | \ge t \Big) \le 2 \exp 
\big( - \frac {t^2} {2N} \big) $$
we have
$$\eqalign {E_\epsilon \exp \Big( \frac {\beta^2} {2N} \big( 
\sum_{i=1}^N \epsilon_i\big)^2 \Big) &\le 2 \int^\infty_0 \frac d {dt} 
( \exp \frac {\beta^2t^2} {2N} ) \exp (- \frac {t^2} {2N}) dt \cr
&= 2\beta^2 / (1-\beta^2)\cr}$$
so that (12.8) follows. \hfill $\square$
\vfil\eject

\noindent{\bf 13.~~Sums of (vector valued) independent random 
variables}

The first objective of this section is to discuss the genesis of some 
key ideas of the isoperimetric approach.  This will be helped by a 
simple (but rather typical) example of application of Theorem 
3.1.1.  We will then discuss, in detail, a situation that parallels 
the situation of Chapter 8, but where the infimum over $\alpha\in
{\Cal F}$ is replaced by a supremum.  There are unexpected and subtle 
differences; this is closely connected to the fact that the conditions 
on the function $h(x,y)$ in Theorem 4.4.1 are (and must be) highly 
disymmetric in $x$ and $y$.

Consider a family ${\Cal F}$ of $N$-tuples $\alpha =(\alpha_i)_{i\le 
N}$, $\alpha_i\ge 0$; and for $y=(y_i)_{i\le N}$ set
$$\leqalignno{Z(y) &= \sup_{\alpha\in{\Cal F}}\sum_{i\le N}\alpha_i
y_i\,. &(13.1)\cr}$$

Consider a sequence $(X_i)_{i\le N}$ of positive independent r.v.  
Consider the r.v.
$$\leqalignno{Z &= Z_{\Cal F}= Z(X_1,\dots ,X_N)=\sup_{\alpha\in
{\Cal F}}~\sum_{i\le N}\alpha_iX_i\,. &(13.2)\cr}$$

We denote by $X^\ast_i$ the non-decreasing rearrangement of the 
sequence $(X_i)_{i\le N}$.  That is,
$$\leqalignno{X^\ast_i &= \sup\{ t\,;\,\card\{ j\le N\,;\,X_j\ge t\}
\ge i\}\,. &(13.3)\cr}$$

It is useful to note that
$$\sum_{i\le k}X^\ast_i = \sup\left\{\sum_{i\in I}X_i\,;\,\card I=k
\right\}\,.$$

The key motivation for Theorem 3.1.1 is the following, where $\tau =
\sup\{\alpha_i\,;\,i\le N\,,\,\alpha\in{\Cal F}\}$.

\proclaim{Proposition 13.1}  Consider $a>0$, $q,k\in\Bbb N$.  Then
$$\leqalignno{P(Z\ge qa+t) &\le {1 \over q^{k+1}P(Z\le a)^q} +P\left(
\tau\sum_{i\le k}X^\ast_i\ge t\right)\,. &(13.4)\cr}$$
\endproclaim

\remark{Comment}  To make this inequality useful, one has to estimate 
the last term, and then to choose $a$, $q$, $k$ in a efficient way.
\endremark

\demo{Proof}  Set $\Omega =\Bbb R^+$.  Consider the set
$$A=A(a)=\{ y\in\Omega^N\,;\,Z(y)\le a\}$$
where $Z(y)$ is given by (13.1).

Consider $J\subset\{ 1,\dots ,N\}$, and $\alpha\in{\Cal F}$.  A key 
observation is that, by positivity 
$$\leqalignno{\sum_{i\in J}\alpha_i y_i &\le \sum_{i\le N}a_iy_i\,. 
&(13.5)\cr}$$

Consider now $y^1,\dots ,y^q\in A$, $x\in\Omega^N$, and set
$$I=\{i\le N\,;\,x_i\not\in\{ y^1_i,\dots ,y^q_i\}\}\,.$$

Consider a partition $(J_j)_{j\le q}$ of $\{ 1,\dots ,N\}\backslash 
J$ such that if $i\in J_j$, then 
$x_i=y^j_i$.  Then, for $\alpha\in{\Cal F}$, by (13.5), since $y^j\in 
A$
$$\sum_{i\not\in I}\alpha_i x_i=\sum_{j\le q}~\sum_{i\in J_j} 
\alpha_iy^j_i\le qa$$
and thus
$$\sum_{i\le N}\alpha_ix_i\le qa+\tau\sum_{i\in I}x_i\,.$$

Hence (with the notation of Section 3.1.1) if $f(A,\dots ,A,x)\le k$, 
we have
$$\leqalignno{\sum_{i\le N}x_i &\le qa+\tau\sum_{i\le k}x^\ast_i\,. 
&(13.6)\cr}$$

If we provide the $i^{\text{th}}$ factor of $\Omega^N$ with the law 
of $X_i$, then (13.4) follows from (3.1.3) and (13.6).\rbx
\enddemo

\remark{Remark}  Certainly one can bound the term $\sum\limits_{i\in 
I}\alpha_ix_i$ in a less brutal way, by\hfil\break 
$\sup\limits_{\alpha\in{\Cal F},\card I=k}~\sum\limits_{i\in I}
\alpha_ix_i$.
\endremark

Consider now a sequence $(X_i)_{i\le N}$ of Banach space valued r.v.  
A number of classical problems of probability (in particular, laws of 
large numbers and laws of the iterated logarithm) depend 
crucially on sharp estimates of the tail probability $P(\Vert
\sum\limits_{i\le N} X_i\Vert\ge t)$.  For many years these estimates 
were found using martingales, and the results were not optimal.  One 
big obstacle is that there is no obvious substitute for the positivity 
arguments that are central to Chapter 8 and to Proposition 13.1.  
Although its importance became clear only later, a crucial 
contribution was made by M. Ledoux [L-1].  It was known at the time 
that in many situations, the tails of $\left\|\sum\limits_{i\le N}
X_i\right\|$ resemble the tails of $\left\|\sum\limits_{i\le N}g_iX_i
\right\|$, where $(g_i)_{i\le N}$ is an independent sequence of 
standard normal r.v. that is independent of the sequence $X_i$.  To 
study $\left\|\sum\limits_{i\le N} g_iX_i\right\|$, Ledoux wrote
$$\leqalignno{\left\|\sum_{i\le N}g_iX_i\right\| &=E_g\left\|\sum_{i
\le N}g_iX_i\right\|+\left(\left\|\sum_{i\le N}g_iX_i\right\|-E_g
\left\|\sum_{i\le N}g_iX_i\right\|\right) &(13.7)\cr}$$
where $E_g$ denotes conditional expectation given $(X_i)_{i\le N}$.  
The idea was that either term of the right-hand side should be easier 
to study than the term of the left-hand side.  This is particularly 
apparent for the second term, where, arguing conditionally on $X_i$, 
one can take advantage of the properties of Gaussian processes.

It turns out that the first term in the right of (13.7) has the exact 
property needed to replace the positivity used in Proposition 13.1; 
namely, if $J\subset\{ 1,\dots ,N\}$, we have
$$\leqalignno{E_g\left\|\sum_{i\in J}g_iX_i\right\| &\le E_g\left\|
\sum_{i\le N}g_iX_i\right\|\,. &(13.8)\cr}$$
The realization of the importance of positivity-like properties led 
first to the characterization of the Banach-space valued r.v. that 
satisfy the law of the iterated logarithm [L-T1].  Perhaps more 
importantly, (13.8) lead this author to the belief that some 
isoperimetric principle should be relevant, and hence to the 
theorem of [T2] (that is now superceeded by the comparable, but much 
easier to prove Theorem 3.1.1), and started the line of investigation 
that culminates in the present paper.

The author also understood that Bernoulli r.v. have regularity 
properties that almost match those of Gaussian r.v. (a crucial step 
is the comparison theorem of [T6]).  They offer the extra advantage 
that the tails of $\left\|\sum_{i\le N}\varepsilon_iX_i\right\|$ 
(where $P(\varepsilon_i =-1)=P(\varepsilon_i=1)=1/2$) always resemble 
the tails of $\left\|\sum\limits_{i\le N}X_i\right\|$.  Thus, 
rather than (13.8) one should write
$$\leqalignno{\left\|\sum_{i\le N}\varepsilon_iX_i\right\| &= 
E_\varepsilon\left\|\sum_{i\le N}
\varepsilon_iX_i\right\| +\left(\left\|\sum_{i\le N}\varepsilon_i
X_i\right\| -E_\varepsilon\left\|\sum_{i\le N}\varepsilon_iX_i\right\|
\right)\,. &(13.9)\cr}$$

To study the last term conditionally on $(X_i)_{i\le N}$, one can 
rely, in particular, upon the following result.

\proclaim{Theorem 13.2}  Consider vectors $(v_i)_{1\le i\le N}$ in a 
Banach space $W$, and set
$$\leqalignno{\sigma &= \left(\sup\left\{\sum_{i\le N}w^\ast (v_i)^2
\colon w^\ast\in W^\ast\,,\,\Vert w^\ast\Vert\le 1\right\}
\right)^{1/2}\,. &(13.10)\cr}$$
Consider a sequence $(X_i)_{i\le N}$ of independent real valued r.v. 
such that $\vert X_i\vert\le 1$.  Denote by $M$ a median of the r.v. 
$\left\|\sum\limits_{i\le N}X_iv_i\right\|$.  Then for $t>0$ we have
$$\leqalignno{P\left(\left|\left\|\sum_{i\le N}X_iv_i\right\|-M
\right|\ge t\sigma \right) &\le 4\exp \left(-{t^2 \over 16}\right)\,. 
&(13.11)\cr}$$
\endproclaim

\demo{Proof}  We observe that if we set
$${\Cal F} = \{(w^\ast (v_i))\,;\, w^\ast\in W^\ast\,,\,\Vert w^\ast
\Vert \le 1\}$$
then
$$Z=\left\|\sum_{i\le N}X_iv_i\right\|=\sup_{\alpha\in{\Cal F}}
\sum_{i\le N}\alpha_iX_i\,.$$

Thus Theorem 13.2 is a special case of Theorem 8.1.1 (using scaling).
\rbx
\enddemo

\remark{Remarks}  1)  Certainly the constant in the exponent is not 
sharp, and could be improved using (4.2.7) rather than (4.1.3), 
especially in the case of Bernoulli r.v., where the use of (4.3.8)
would yield a bound of $\exp -{1 \over 4}(t-\sqrt{\log 2})^2$ for $t
\ge\sqrt{\log 2}$.

2)  There is another bound on the tails of $\left\|\sum\limits_{i\le 
N}X_i v_i\right\|$ namely the trivial bound $\left\|\sum\limits_{i\le 
N}X_iv_i\right\|\le\sup\limits_{\Vert w^\ast\Vert\le 1}~
\sum\limits_{i\le N}\vert w^\ast (v_i)\vert$, and it is, of course, 
possible to interpolate between this bound and (13.11).  This can be 
done as follows.  For a sequence $(r_i)_{\le N}$ of real numbers, 
and $t>0$ we write
$$K_{1,2}((r_i),t)=\inf\left\{\sum_{i\le N}\vert u_i\vert +t\left(
\sum_{i\le N}w^2_i\right)^{1/2}\,;\,r_i=u_i+w_i\right\}$$
where the infimum is taken over all possible decompositions $r_i=u_i
+w_i$.  We set
$$\kappa (t)=\sup\{ K_{1,2}(w^\ast (v_i),t)\,;\,w^\ast\in W^\ast ,
\Vert w^\ast\Vert =1\}\,.$$

We observe that $\kappa (t)\le t\sigma$.  Only rather trivial 
modifications to the proof of Theorem 8.1.1 are needed to see that 
one can improve (13.11) into
$$\leqalignno{P\left(\left|\left\|\sum_{i\le N} X_iv_i\right\| -M
\right|\ge \kappa (t)\right) &\le 2\exp \left( -{t^2 \over 16}\right)
\,. &(13.12)\cr}$$
This inequality streamlines a result of [D-MS].

If one observes that $\kappa (2t)\le 2\kappa (t)$ one obtains through 
a routine computation that, for all $p\ge 1$
$$\left\|\sum_{i\le N}X_iv_i\right\|_p\le M+K\kappa (\sqrt{p})\,,$$
a rather precise form of the so called Kintchin-Kahane inequalities.  
It should also be pointed out that by (13.12), $\left\|\sum\limits_{i
\le N}X_iv_i\right\|_p\ge M-K\kappa (\sqrt{p})$, that $\left\|
\sum\limits_{i\le N}X_iv_i\right\|_p\ge M2^{-1/p}$ (obviously) and 
that $\left\|\sum\limits_{i\le N}X_i v_i\right\|_p\ge K\kappa (
\sqrt{p}/K)$.  To prove this last inequality, one reduces to the 
real-valued case; it is simple to see that this follows from 
[L-T1] lemma 4.9 (see also [M-S]).

We will not pursue the discussion of how to build upon (13.4), 
Theorems 3.1.1 and 13.2 to obtain sharp bounds for the tails of 
$\left\|\sum\limits_{i\le N}X_i\right\|$.  This is done in great detail 
in [L-T2], Theorem 6.17.  An alternative approach, that relies rather 
on Theorem 4.2.4 is developed in [T3].  The bounds obtained through 
these methods apparently are sharp in all situations (see however, 
[Ro] for a case where other ingredients are needed).  Instead, we 
continue the investigation of r.v. of the type $Z=\sup\limits_{\alpha
\in{\Cal F}}~\sum\limits_{i\le N} \alpha_iX_i$ that was started in 
Chapter 8.  In order to apply Corollary 8.2.2, we need to have 
(4.4.6), where $h(x,y)=\vert x-y\vert$, or, if $\alpha_i$ is always 
positive, $h(x,y)=(x-y)^+$.  When the variable $X_1$ is positive 
(i.e., its law $\mu$ is supported by $\Bbb R^+$), inspection of 
Theorem 4.4.1 shows that (whatever choice of $\theta$, $\xi$) no 
integrability condition on $X_1$ except boundedness, will insure 
that the conditions of this theorem hold for this choice of $h$.  We 
will now give an example that shows that this is not an artifact of 
our approach.  We will show that (13.11) cannot be essentially 
improved, even if $P(\vert X_i\vert\not= 0)$ is arbitrary small.  
This implies (by scaling) that, given any {\it finite} function $\Phi
\colon\Bbb R^+\to\Bbb R^+$, with $\Phi (0)=0$, one can find a 
real r.v. $x_1$ with $\int\Phi (\vert X_1\vert )d\lambda\le 1$, and 
vectors $(v_i)_{i\le N}$ such that (13.11) is violated.

\remark{Example 13.3}  This example is essentially a re-interpretation 
of the example presented at the end of Section 4.3.  Consider an 
independent sequence $(X_i)_{i\le N}$ of Bernoulli variables such 
that $P(X_i=1)=p$ is small.  Consider the family ${\Cal F}$ of 
$N$-tuples of form $\alpha_i=1/\sqrt{2pN}$ if $i\in I$, $\alpha_i=0$ 
otherwise, where $I$ varies over all subsets of $\{ 1,\dots ,N\}$ of 
cardinality $\le 2pN$.  Then $\sigma =1$.  Consider
$$Z=\sup_{\alpha\in{\Cal F}}~\sum_{i\le N}\alpha_iX_i = {1 \over 
\sqrt{2Np}}\sup\left\{ \sum_{i\in I}X_i\,;\,\card I\le 2Np\right\}\,.$$
We can also view $Z$ as $\left\|\sum\limits_{i\le N}X_i e_i
\right\|_{\Cal F}$, where $e_i$ is the canonical basis of $\Bbb R^N$, 
and where the norm $\Vert\cdot\Vert_{\Cal F}$ is given by $\Vert x
\Vert_{\Cal F}=\sup\limits_{\alpha\in{\Cal F}}~\sum\limits_{i\le N}
\alpha_i\vert x_i\vert$.
\endremark

The main observation is that
$$\sum_{i\le N}X_i\le 2Np\Rightarrow Z={1 \over \sqrt{2Np}}\sum_{i\le 
N}X_i\,.$$
Since the probability of the event on the left goes to $1$, as $N\to
\infty$, the r.v. $Z$ is asymptotically normal, of mean $\sqrt{Np/2}$ 
and variance $\sqrt{(1-p)/2}$; so its deviation from 
its median do not decay faster than $\exp -Kt^2$.

The conclusion to be drawn from Example 13.3 is that, in order to 
extend Theorem 13.2 to the case where $X_1$ is unbounded, we must 
require conditions of a different nature than integrability.

\proclaim{Theorem 13.4}  There exists a universal constant $L$ with 
the following property.  Consider a convex function $\psi$ on 
$\Bbb R^+$ such that $\psi (x)\le x^2$ if $x\le 1$ and $\psi (x)\ge 
x$ if $x\ge 1$.  Consider a probability measure $\mu$ on $\Bbb R$.  
Assume the following
$$\leqalignno{\forall t>0\,,\quad \mu (\{x\,;\,\vert x\vert\ge t\} ) 
&\le 2\exp (-L\psi (2t))\,.  &(13.13)\cr}$$
Given any subset $B$ of $\Bbb R$, with $\mu (B)\ge 1/2$, and any 
$t\ge 1$, we have
$$\leqalignno{\mu (\{x\,;\,\psi (\inf_{y\in B}\vert x-y\vert )\ge t
\}) &\le e^{-t}(1-\mu (B))\,.  &(13.14)\cr}$$

Consider independent real valued r.v. $(X_i)_{i\le N}$ distributed 
like $\mu$, and vectors $(v_i)_{i\le N}$ is a Banach space $W$.

Then, for all $t>0$, we have
$$\leqalignno{P\left(\left|\left\|\sum_{i\le N} X_iv_i\right\| -M
\right|\ge t\right) &\le 2\exp \left(-{1 \over L}\Psi_{\Cal F}(t)
\right) &(13.15)\cr}$$
where $M$ is a median of $\left\|\sum\limits_{i\le N}X_iv_i\right\|$, 
where
$${\Cal F} = \{(w^\ast (v_i))_{i\le N}\,;\,w^\ast\in W^\ast\,,\,\Vert 
w^\ast\Vert\le 1\}$$
and where $\Psi_{\Cal F}$ is defined in Section 8.2.
\endproclaim

\demo{Proof}  According to Corollary 8.2.2, it suffices to prove that 
the hypothesis of Theorem 4.4.1 hold when $h(x,y)=\vert x-y\vert$, in 
the case $\theta (x)=-\log x$, $w(x)=-{1 \over 2}\log x$ (so 
that $H(\xi ,w)$ holds by Proposition 2.6.1).  Only (4.4.4) has to be 
checked, since (13.14) is a rewriting of (4.4.5).

Consider $B\subset\Bbb R$ with $\mu (B)\le 1/2$.  Consider
$$\leqalignno{s &= \inf\{\vert y\vert\,;\, y\in B\}\,. &(13.16)\cr}$$
Clearly, $h(x,B)\le \vert x\vert +s$.  Thus, by convexity of $\psi$ 
we have
$$\int_{\Bbb R}\exp \psi (h(x,B))d\mu (x)\le\exp{1 \over 2}\psi (2s)
\int_{\Bbb R}\exp {1 \over 2}\psi (2x)d\mu (x)\,.$$
On the other hand, by (13.16) we have $B\cap ]-s,s[ =\emptyset$, so 
that $B\subset\{x\,;\,\vert x\vert\ge s\}$ and hence by (13.13) we 
have $\exp\psi (2s)\le (2/\mu (B))^{1/L}$.  Thus it suffices to show 
that for $L$ large enough we have
$$x\le {1 \over 2}\Rightarrow I\left({2 \over x}\right)^{1/2L}\le 
{1 \over \sqrt{x}}~~(=\exp -w(x))$$
where $I=\int_{\Bbb R}\exp {1 \over 2}\psi (2x)d\mu (x)$.  Thus it 
remains to show that under (13.13) $\lim\limits_{L\to\infty}I=0$, 
uniformly in $\psi$, an easy exercise left to the reader.\rbx
\enddemo

Theorem 13.4 can be applied to the case where $\mu$ is a measure 
$\nu_\psi$ of the type considered in Proposition 2.7.4, although in 
that case the simpler Theorem 2.7.1 will yield the same conclusion.  
There are however, situations covered by Theorem 13.4 that are not 
covered by Theorem 2.7.1, because in (13.14) we require only $t\ge 1$.  
In particular, if the law of $X$ satisfies (13.14), and if 
$\Vert Z\Vert_\infty\le 1$, the law of ${1 \over 3}(X+Z)$ satisfies 
(13.14) (it is {\it not} required that $Z$ be independent of $X$).  
(The corresponding statement for (13.13) is also true, under mild 
conditions on $\psi$, replacing if needed $1/3$ by a smaller number.)

In conclusion of this section, we want to discuss a question that 
apparently is not fully clarified by the results of the present paper.  
Consider numbers $(a_i)_{i\le N}$, and vectors $(v_i)_{i\le N}$ 
in a Banach space.  Of which order are the fluctuations of the r.v.
$Z=\left\|\sum\limits_{i\le N}a_{\rho (i)}v_i\right\|$ around its 
median $M$, when $\rho$ is seen as a random element of the 
symmetric group $S_N$, provided with the uniform probability $P$?

\proclaim{Proposition 13.5}  a)  Assuming $\vert a_i\vert\le 1$ for 
each $i$, we have
$$\leqalignno{t\ge 0 &\Rightarrow P(\vert Z-M\vert\ge t)\le 4\exp 
-{t^2 \over 16\sigma^2} &(13.17)\cr}$$
where as usual
$$\sigma^2=\sup\left\{\sum_{i\le N} w^\ast (v_i)^2\,;\,w^\ast\in 
W^\ast\,,\,\Vert w^\ast\Vert\le 1\right\}\,.$$

b)  Assuming $\Vert v_i\Vert\le 1$ for each $i$, we have
$$\leqalignno{t\ge 0 &\Rightarrow P(\vert Z-M\vert\ge t)\le 4\exp 
-{t^2 \over 16\sum\limits_{i\le N} a^2_i}\,. &(13.18)\cr}$$
\endproclaim

\remark{Remark} A first problem is to find a bound that contains 
simultaneously (13.17) and (13.18).
\endremark

\demo{Proof}  The proof follows that of Theorem 8.1.1, using now 
Theorem 5.1 rather than Theorem 4.1.1.  Thus, we indicate only the 
key points.

To prove a), one notes that for $\rho,\tau\in S_N$, and $I=\{i\le N
\,;\,\rho (i)\not=\tau (i)\}$, then
$$\left|\sum_{i\le N}w^\ast (v_i)a_{\rho (i)}-\sum_{i\le N}w^\ast 
(v_i)a_{\tau (i)}\right|\le \sum_{i\in I}\vert w^\ast (v_i)\vert\,.$$

To prove b), one observes that $Z$ has the same distribution as 
$\left\|\sum\limits_{i\le N}a_i v_{\rho (i)}\right\|$; and, with 
the notations above, one now has
$$\eqalignno{\left|\sum_{i\le N}w^\ast (v_{\rho (i)})a_i-\sum_{i\le 
N}w^\ast (v_{\tau (i)})a_i\right| &\le\sum_{i\in I}\vert a_i\vert\,. 
&\bx\cr}$$
\enddemo

It should be pointed out that it is most likely that a phenomenon 
similar to that of Example 13.3 occurs in case a, and that (13.16) 
cannot be improved even if a large majority of the numbers $a_i$ 
are equal to zero.
\vfil\eject

\baselineskip 14pt plus 2pt

\centerline{\bf References}

\exam{[A-L-R]}{M. Aizenman, J. L. Lebowitz, D. Ruelle, Some rigorous 
results on the  Sherrington-Kirkpatrick spin glass model, Commun. 
Math. Phys. 112, 1987, 3-20.}

\exam{[A-S]}{N. Alon, J. Spencer, The Probabilistic Method, Wiley, 
1991.}

\exam{[B1]}{B. Bollob\'as, The chromatic number of random graphs, 
Combinatorica 8, 1988, 49-55.}

\exam{[B2]}{B. Bollab\'as, Random graphs revisited, Proceeding of 
Symposia on Applied Mathematics, Vol. 44, 1991, 81-98.}

\exam{[B-B]}{B. Bollob\'as, G. Brightwell, The height of a random 
partial order:  Concentration of Measure, Annals of Applied Probab., 
2, 1992, 1009-1018.}

\exam{[C-N]}{F. Comets, J. Neveu, The Sherrington-Kirkpatrick Model 
of Spin Classes and Stochastic Calculus:  the high temperature case.}

\exam{[D-MS]}{S. Dilworth, S. Montgomery-Smith,}

\exam{[F]}{A. M. Frieze, On the length of the longest monotone 
subsequence in a random permutation, Ann. Appl. Prob. 1, 1991, 
301-305.}

\exam{[H]}{W. Hoeffding, Probability inequalities for sums of bounded 
random variables, J. Amer. Statist. Assoc. 58, 1963, 13-30.}

\exam{[J]}{S. Janson, Poisson approximation for large derivations, 
Random Structures and Algorithms 1, 221-290.}

\exam{[Ka]}{R. M. Karp, An upper bound on the expected cost of an 
optimal assignment; In {\it Discrete Algorithm and Complexity:  
Proceedings of the Japan-US joint Seminar}, Academic Press, 1987, 
1-4.}

\exam{[K1]}{H. Kesten, Aspects of first-passage percolation, Ecole 
d' Et\'e de Probabilit\'e de Saint-Flour XIV, Lecture Notes in Math 
1180, 125-264, Springer, New York.}

\exam{[K2]}{H. Kesten, On the speed of convergence in first passage 
percolation, Ann. Applied Probab. 3, 1993, 296-338.}

\exam{[K-S]}{R. M. Karp, J. M. Steele, Probabilistic analysis of 
heuristics, in The Traveling Salesman Problem, John Wiley and Sons, 
1985, 181-205.}

\exam{[Lea]}{J. Leader, Discrete Isoperimetric inequalities, 
Proceeding of Symposia on Applied Mathematics, Vol. 44, 1991, 57-80.}

\exam{[L]}{M. Ledoux, Gaussian randomization and the law of the 
iterated logarithm in type 2 Banach spaces, Unpublished manuscript, 
1985.}

\exam{[L-T1]}{M. Ledoux, M. Talagrand, Characterization of the law of 
the iterated logarithm in Banach spaces, Ann. Probab. 16, 1988, 
1242-1264.}

\exam{[L-T2]}{M. Ledoux, M. Talagrand, Probability in Banach Spaces, 
Springer Verlag, 1991.}

\exam{[Lu]}{T. Luczak, The chromatic number of Random graphs, 
Combinatorica 11, 1991, 45-54.}

\exam{[M1]}{B. Maurey, Construction de suites sym\'etriques, Comptes 
Rendus Acad. Sci. Paris 288, 1979, 679-681.}

\exam{[M2]}{B. Maurey, Some deviation inequalities, Geometric and 
Functional Analysis 1, 1991.}

\exam{[McD]}{C. McDiarmid, On the method of bounded differences, in 
``Survey in Combinatorics (J. Simons, Ed.) London Mathematical 
Society Lecture Notes, Vol. 141, Cambridge Univ. Press, London/New 
York, 1989, 148-188.}

\exam{[M-H]}{C. McDiarmid, Ryan Hayward, Strong concentration for 
Quicksort, Proceedings of the Third Annual ACM-SIAM Symposium on 
Discrete Algorithms (SODA), 1992, 414-421.}

\exam{[M-S]}{V. Milman, G. Schechtman, Asymptotic theory of finite 
dimensional normed spaces, Lecture Notes in Math 1200, Springer 
Verlag, 1986.}

\exam{[P]}{G. Pisier, Probabilistic methods in the geometry of Banach 
spaces.  Probability and Analysis.  Varena (Italy) 1985.  Lecture 
Notes in Math. 1206, Springer Verlag, 1986, 167-241.}

\exam{[R1]}{W. Rhee, On the fluctuations of the stochastic traveling 
salesperson problem, Math. of Operation Research 16, 1991, 482-489.}

\exam{[R2]}{W. Rhee, A matching problem and subadditive Euclidean 
functionals, Ann. Applied. Probab. 3, 1993, 794-801.}

\exam{[R3]}{W. Rhee, On the fluctuations of simple matching, 
Manuscript, 1992.}

\exam{[Ro]}{J. Rosinski, Remarks on a Strong Exponential 
Integrability of Vector Valued Random Series and Triangular Arrays,}

\exam{[R-T]}{W. Rhee, M. Talagrand, A sharp deviation inequality for 
the stochastic traveling salesman problem, Ann. Probab. 17, 1989, 
1-8.}

\exam{[S-S]}{E. Shamir, J. Spencer, Sharp concentration of the 
chromatic number of random graphs $G_{n,p}$, Combinatorica 7, 
121-129.}

\exam{[T1]}{M. Talagrand, An isoperimetric theorem on the cube and 
the Kintchine Kahane inequalities, Proc. Amer. Math. Soc. 104, 1988, 
905-909.}

\exam{[T2]}{M. Talagrand, Isoperimetry and integrability of the sum 
of independent Banach space valued random variables, Ann. Probab. 
17, 1989, 1546-1570.}

\exam{[T3]}{M. Talagrand, A new isoperimetric inequality for product 
measure, and the tails of sums of independent random variables, 
Geometric and Functional analysis 1, 1991, p. 211-223.}

\exam{[T4]}{M. Talagrand, A new isoperimetric inequality for product 
measure, and the concentration of measure phenomenon, Israel Seminar 
(GAFA), Springer Verlag Lecture Notes in Math. 1469, 1991, p. 94-124.}

\exam{[T5]}{M. Talagrand, Some isoperimetric Inequalities and their 
applications, Proceedings of the International Congress of 
Mathematicians, Kyoto 1990, Springer Verlag, 1991, p. 1011-1029.}

\exam{[T6]}{M. Talagrand, Regularity of infinitely divisible 
processes, Ann. Probab. 21, 1993, 362-432.}

\exam{[T7]}{M. Talagrand, Supremum of some canonical processes, Amer. 
J. Math., to appear.}

\exam{[W]}{D. W. Walkup, On the expected value of a random assignment 
problem, SIAM J. Comput. 8, 1979, 440-422.}

\exam{[Y]}{V. V. Yurinskii, Exponential bounds for large deviations.  
Theor. Prob. Appl. 19, 1974, 154-155.}
\bigskip

\settabs 3\columns

\+Equipe d'Analyse - Tour 48 &&Department of Mathematics\cr
\+U.A. au C.N.R.S. n$^\circ$ 754 &&The Ohio State University\cr
\+Universit\'e Paris VI &and &231 West 18th Avenue\cr
\+4 Pl Jussieu &&Columbus, Ohio  43210\cr
\+75230 Paris Cedex 05 &&USA\cr

\end